\documentclass[preprint, 12pt, authoyear]{elsarticle}

\usepackage{amssymb}
\usepackage[all,arc]{xy}
\usepackage{enumerate}
\usepackage{mathrsfs}
\usepackage{amsbsy}
\usepackage{amsfonts}
\usepackage{amsmath}
\usepackage{amssymb}
\usepackage{amsthm}

\usepackage{bbm} 
\usepackage{caption}
\usepackage{color}

\usepackage[mathscr]{eucal}

\usepackage{float}

\usepackage{graphicx}  
\usepackage{morefloats}

\usepackage{hyperref} 

\usepackage{latexsym}  
\usepackage{listings} 
\usepackage{lscape}  

\usepackage{multirow,array}
\newcolumntype{P}[1]{>{\centering\arraybackslash}m{#1}}

\usepackage{rotating} 
\usepackage{stmaryrd}
\usepackage{subcaption, graphicx}

\usepackage{url}
\usepackage{verbatim}
\usepackage{wrapfig}  
\usepackage[nocenter]{qtree}

\usepackage{epstopdf}
\DeclareGraphicsExtensions{.eps}

\newtheorem{thm}{Theorem}[section]
\newtheorem{cor}[thm]{Corollary}
\newtheorem{prop}[thm]{Proposition}
\newtheorem{lem}[thm]{Lemma}

\theoremstyle{definition}
\newtheorem{defn}[thm]{Definition}

\theoremstyle{remark}

\newcommand{\Z}{\mathbb{Z}}
\newcommand{\R}{\mathbb{R}}

\newcommand{\myz}{\omega}

\allowdisplaybreaks

\makeatletter
\let\c@equation\c@thm
\makeatother
\numberwithin{equation}{section}

\makeatletter
\def\ps@pprintTitle{%
 \let\@oddhead\@empty
 \let\@evenhead\@empty
 \def\@oddfoot{}%
 \let\@evenfoot\@oddfoot}
\makeatother

\begin{document}

\begin{frontmatter}

\title{Discrete Directional Gabor Frames}

\author{Wojciech Czaja, Benjamin Manning, James M. Murphy, Kevin Stubbs}

\begin{abstract}We develop a theory of discrete directional Gabor frames for functions defined on the $d$-dimensional Euclidean space.  Our construction incorporates the concept of ridge functions into the theory of isotropic Gabor systems, in order to develop an anisotropic Gabor system with strong directional sensitivity. We present sufficient conditions on a window function $g$ and a sampling set $\Lambda_{\omega}$ for the corresponding directional Gabor system $\{g_{m,t,u}\}_{(m,t,u)\in\Lambda_{\omega}}$ to form a discrete frame.  Explicit estimates on the frame bounds are developed. A numerical implementation of our scheme is also presented, and is shown to perform competitively in compression and denoising schemes against state-of-the-art multiscale and anisotropic methods, particularly for images with significant texture components.  
\end{abstract}

\end{frontmatter}

\section{Introduction}\label{sec1}

Given a square integrable function $g \in L^2(\R)$ and constants
$a,b>0$, the associated {\em Gabor system} (also known as {\em Weyl--Heisenberg system}) generated by $g$ and the lattice $a\Z \times b\Z$, $\mathcal{G} (g, a,b) = \{ g_{m,n} \}_{m,n \in \Z}$, is defined by
\[
g_{m,n} (x) = e^{2 \pi i am x} g(x - bn).
\]
In 1946, Dennis Gabor proposed to study such systems for their usefulness in the analysis of information conveyed by communication channels \cite{Gabor46}. The resulting theory led to many applications ranging from auditory signal processing, to pseudodifferential operator analysis, to uncertainty principles. The edited volumes by Benedetto and Frazier \cite{benedetto1993} and by Feichtinger and Strohmer \cite{FS1, FS2}, as well as Gr\"ochenig's treatise \cite{Grochenig}, provide detailed treatments of various aspects of this rich and beautiful theory.

Among many interesting recent developments in time-frequency analysis, we want to consider the set of ideas which expands the notion of traditional Gabor systems by including directional information. Directional content is especially important in many applications, such as those in remote sensing or medical imaging, where interesting features often propagate in a particular direction. Thus, it is not surprising that this avenue of research achieved significant success in the closely related field of wavelet analysis. Contourlets \cite{DoV03}, curvelets and ridgelets \cite{acha_curvelets1, acha_curvelets2, CG02}, bandlets \cite{LePM}, wedgelets \cite{Don99}, and shearlets \cite{GKL06, LLKW, shearlet_acha, shearlet_edge}, are just a few examples of directional multiscale constructions introduced to address the problems of identifying edges and directions in various forms of imagery.

Counter to a common belief that in image analysis wavelet-based techniques are superior to other approaches, Gabor methods play a significant role in several important areas, including fingerprint recognition, texture analysis, and vasculature detection. These examples of data are not covered by typical assumptions about existence of a sparse, curvilinear set of singularities, and thus present a new type of difficulties. To overcome these difficulties in the context of Gabor analysis, an effort similar to that in the directional wavelet theory has been undertaken. Directional Gabor filters have been proposed as a model for multichannel neuronal behavior by Daugman \cite{D1, D2} and Watson \cite{Watson}, and to solve problems in texture analysis by Porat and Zeevi \cite{PZ}. Wang et al. proposed to use directional Gabor filters in character recognition \cite{Wang}. Hong, Jain, and Wan proposed to use Gabor ridge filter banks for fingerprint enhancement \cite{Hong}. Gabor methods have also been employed to achieve state-of-the-art performance in image denoising, compared to multiscale anisotropic methods, such as shearlets and curvelets \cite{Gabor_Framelets}. Directional Fourier methods have also been successfully deployed for edge detection \cite{Czaja+Wickerhauser, CW2}.

On the theoretical side, Grafakos and Sansing proposed the concept of directional time-frequency analysis in \cite{Grafakos_Sansing}. In this paper, they introduce directional sensitivity in the time-frequency setting by considering projections onto the elements of the unit sphere. The Radon transform arises naturally in this context and it enables continuous and semi-continuous reconstruction formulas.  Giv proposed a related variant of short-time Fourier transform \cite{Giv} and established its orthogonality relations, as well as some operator theoretic properties. The relationship with the classical short-time Fourier transform and the quasi-shift invariance of the new scheme has been analysed in \cite{AceskaGiv}. Results from \cite{Grafakos_Sansing} have been extended to more general continuous representations in Sobolev spaces in \cite{CFM}. Gabor systems have also been analyzed from the viewpoint of shearlet theory \cite{gabor_shearlets}.  

In the present paper, we develop a theory of discrete directional Gabor frames, related to the concept of Gabor ridge systems.  In Section \ref{sec2} we review some of the necessary features of time frequency analysis.  In Section \ref{sec3} we characterize certain minimal constraints on the spaces of functions that can be efficiently represented by directional Gabor systems. An elementary example is provided in Section \ref{sec4}, and Section \ref{sec5} is concerned with the derivation of sufficient conditions for the existence of discrete directional Gabor systems. We close with some numerical examples in Section \ref{sec6}.

\section{Background}\label{sec2}

We begin with a review of some relevant background in frame theory and Gabor analysis.

\begin{defn}\label{discrete_frame}

Let $\mathcal{H}$ be a Hilbert space.  A discrete set $\{\phi_{i}\}_{i\in I}\subset \mathcal{H}$ is called a \emph{(discrete) frame for $\mathcal{H}$} if there exist constants $0<A\le B<\infty$ such that:

\begin{align*} \forall f\in \mathcal{H}, \ A\|f\|_{\mathcal{H}}^{2}\le \sum_{i\in I}|\langle f, \phi_{i}\rangle_{\mathcal{H}}|^{2} \le B\|f\|_{\mathcal{H}}^{2}.\end{align*}The optimal choices of $A,B$ are the \emph{frame bounds}.  A frame is said to be \emph{tight} if $A=B$.  It is \emph{Parseval} if $A=B=1$.  

\end{defn}

\begin{defn}\label{continuous_frame}

Let $\mathcal{H}$ be a separable Hilbert space, $X$ a locally compact Hausdorff space equipped with a positive Radon measure $\mu$ such that supp($\mu)= X$.  A set $\{\psi_{x}\}_{x\in X}$ is a \emph{continuous frame with respect to $\mu$} for $\mathcal{H}$ if there exist constants $0<A\le B<\infty$ such that:

\begin{align*}\forall f\in \mathcal{H}, \ A\|f\|_{\mathcal{H}}^2\le \int_{X}|\langle f, \psi_{x}\rangle_{\mathcal{H}}|^2 \; d\mu(x) \le B\|f\|_{\mathcal{H}}^2.\end{align*}
\end{defn}

Definition \ref{continuous_frame} generalizes Definition \ref{discrete_frame} by allowing for a continuous indexing of the frame elements.  In the study of continuous frames, the \emph{discretization problem} arises naturally: when can a continuous frame be discretized to acquire a discrete frame?  That is, when can the continuous indexing set be replaced with a discrete indexing set, while still retaining the frame property?  A prominent abstract approach to the discretization problem involves the coorbit space theory of Feichtinger and Gr\"{o}chenig \cite{atomic_decomp_FG}, \cite{Fornasier_Rauhut}.  This method considers the representations induced by the group structure of the indexing set $X$.  Thus, coorbit space theory is most effective when the continuous frame is parametrized by a group.  This is the case for continuous wavelet systems and continuous shearlet systems \cite{shearlet_coorbit}, which can be parametrized by the affine and shearlet groups, respectively.

There is also substantial interplay between discrete and continuous representations in the context of Gabor theory, that is, between discrete Gabor frames and the short-time Fourier transform.  We are interested in studying this interplay in the context of Gabor systems that have a degree of directional sensitivity incorporated.  The anisotropic generalization of Gabor systems developed by Grafakos and Sansing is continuous, and the theory they develop proves certain systems are continuous representations for particular function spaces.  However, a discrete theory is desired, both out of pure mathematical interest and for applications purposes.  Some partial results towards this exist, such as Theorem \ref{semi_discrete_thm}.  However, a regime for full discretization has not yet been proved.  This article develops sufficient conditions for full discretization; see Theorem \ref{main_result}.  

We now consider some relevant background on classical Gabor theory and the continuous theory of directional Gabor systems.  Unless stated otherwise, $\langle \cdot, \cdot\rangle$ denotes the usual $L^{2}$ inner product.

\begin{defn}For $g\in L^{2}(\mathbb{R}^{d})$ and $m,t\in\mathbb{R}^{d}$, $g^{m,t}\in L^{2}(\mathbb{R}^{d})$ is the function defined by

\begin{align*}g^{m,t}(x):=e^{2\pi i m\cdot x}g(x-t).\end{align*} 

\end{defn}

\begin{defn}Let $f \in L^{1}(\mathbb{R}^{d})$.  The \emph{Fourier transform} of $f$ is the function $\hat{f}\in L^{\infty}(\mathbb{R}^{d})$ given by

\begin{align*}\hat{f}(\gamma)=\int_{\mathbb{R}^{d}}f(x)e^{-2\pi i \gamma \cdot x}dx.\end{align*}

\end{defn}

Note that the Fourier transform is uniquely extendable to a unitary operator on $L^{2}(\mathbb{R}^{d})$; this result is Plancherel's Theorem \cite{BenedettoCzaja}.  The inverse Fourier transform of $f$ is denoted $\check{f}$.

\begin{defn}  Let $g\in L^{2}(\mathbb{R}^{d})$.  The \emph{short-time Fourier transform} of a function $f\in L^{2}(\mathbb{R}^{d})$ with respect to the window function $g$ is the function $V_{g}(f):\mathbb{R}^{2d}\rightarrow\mathbb{C}$ given by 

\begin{align*}V_{g}(f)(m,t):=\int_{\mathbb{R}^{d}}f(x)\overline{g(x-t)}e^{-2\pi i m\cdot x} \; dx, \ t,m\in\mathbb{R}^{d}.\end{align*}

\end{defn}

The short-time Fourier transform is also called a ``sliding window Fourier transform," \cite{Grochenig} or ``voice transform," for its connections to speech processing \cite{STFT_speech}.  The window function $g$ may be chosen with desired regularity and localization properties.  An important property of the short-time Fourier transform is its invertibility.  Indeed, the original function, $f$, is recovered from $V_{g}(f)$ by integrating against translations and modulations of any $\psi\in L^{2}(\mathbb{R}^{d})$ such that $\langle \psi, g\rangle \neq 0$: for such a $\psi$,

\begin{align}\label{STFT_reproducing}f = \frac{1}{\langle \psi, g\rangle}\int_{\mathbb{R}^{d}}\int_{\mathbb{R}^{d}}V_{g}(f)(m,t)\psi^{m,t}\; dm\; dt,\end{align}where convergence is pointwise almost everywhere.  In particular, if $\psi = g$ , we see the short-time Fourier transform is a continuous frame. 

The \emph{Radon transform}, in particular the Fourier slice theorem, shall play a crucial role in the theory of discrete directional Gabor frames.  Recall that $\mathcal{S}(\mathbb{R}^{d})$ denotes the space of Schwartz functions on $\mathbb{R}^{d}$ \cite{Hormander}.

\begin{defn}

Let $f\in\mathcal{S}(\mathbb{R}^{d})$.  The \emph{Radon transform} of $f$ is the function $R(f): S^{d-1}\times\mathbb{R}\rightarrow\mathbb{C}$ given by the formula

\begin{align*} R(f)(u,s):=\int_{u\cdot x=s}f(x) \; dx, \ u\in S^{d-1}, \ s\in\mathbb{R}.\end{align*}

\end{defn}

The notation $R_{u}(f):=R(f)(u, \cdot)$ shall be used in what follows.  The operator $R: f\mapsto R_{u}(f)$ can be extended to a continuous operator from $L^{1}(\mathbb{R}^{d})$ to $L^{1}(\mathbb{R})$ uniformly in $u\in S^{d-1}$.  Since the Fourier transform is naturally defined on $L^{1}(\mathbb{R}^{d})$, it is reasonable to consider its relation to the Radon transform.  The Fourier slice theorem describes this relationship.

\begin{thm} (Fourier slice theorem)  Let $f\in L^{1}(\mathbb{R}^{d})$.  Then the Fourier transform of $f$ and $R(f)$ are related in the following way:

\begin{align*}\widehat{R_{u}(f)}(\gamma) =\hat{f}(\gamma u).\end{align*}

\end{thm}

The continuous theory of Grafakos and Sansing involves weighting the frame elements in the frequency domain.  When the weights are polynomials, this corresponds to differentiation in the time domain.  Indeed, this motivates the following definition.  

\begin{defn}For $\alpha>0$, the \emph{differential operator $D^{\alpha}$ of order $\alpha$} for functions $h:\mathbb{R}\rightarrow\mathbb{C}$ is given by

\begin{align*}D^{\alpha}(h)=(\hat{h}(\gamma)|\gamma|^{\alpha})^{\vee}.\end{align*}

\end{defn}

This definition naturally extends to multi-indices $\alpha=(\alpha_{1},...,\alpha_{d})\in\mathbb{N}$ for functions $h:\mathbb{R}^{d}\rightarrow\mathbb{C}$. 

We now come to the main definition of \cite{Grafakos_Sansing}.

\begin{defn}Let $g\in\mathcal{S}(\mathbb{R})$ be a real-valued, not identically zero window function.  For $m,t\in\mathbb{R}$, $d\in\mathbb{Z}^{+}$, we define

\begin{align*}G^{m,t}(s)=D^{\frac{d-1}{2}}(g^{m,t})(s)=(\widehat{g^{m,t}}(\gamma)|\gamma|^{\frac{d-1}{2}})^{\vee}(s), \ s\in\mathbb{R}.\end{align*}The \emph{weighted Gabor ridge functions generated by $g$} are functions $G_{m,t,u}:\mathbb{R}^{d}\rightarrow\mathbb{C}$ given by

\begin{align*}G_{m,t,u}(x)=G^{m,t}(u\cdot x)=(\widehat{g^{m,t}}(\gamma)|\gamma|^{\frac{d-1}{2}})^{\vee}(u\cdot x),  \ x\in\mathbb{R}^{d}, u\in S^{d-1}.\end{align*}
\end{defn}
The convention of using lower case letters to refer to the generator and capital letters to refer to the weighted Gabor ridge function was established in \cite{Grafakos_Sansing}, and shall be maintained throughout this section.  We note that ridge functions have seen significant application in the theories of function approximation and neural networks \cite{ridge_Candes}, \cite{ridge_Chui}.  

It is clear that the weighted Gabor ridge function $G_{m,t,u}$ is constant along any hyperplane $\{x\in\mathbb{R}^{d} \ | \  u\cdot x=C\}$ for constant $C$.  Along the direction $u$, these functions modulate like a one dimensional Gabor function. The weighting in the Fourier domain permits reconstruction of a signal $f$ from the coefficients  
\begin{align*}\{\langle f, G_{m,t,u}\rangle\}_{(m,t,u)\in\mathbb{R}\times\mathbb{R}\times S^{d-1}}.\end{align*}
Weighted Gabor ridge functions yield a continuous reproducing formula that is analogous to (\ref{STFT_reproducing}), with an additional integral to account for the directional character of the ridge function construction \cite{Grafakos_Sansing}.

\begin{thm}Let $g, \psi\in \mathcal{S}(\mathbb{R})$ be two window functions such that $\langle g, \psi\rangle\neq 0$.  Suppose $f\in L^{1}(\mathbb{R}^{d})$ and $\hat{f}\in L^{1}(\mathbb{R}^{d})$.  Then,

\begin{align*}f=\frac{1}{2\langle g, \psi\rangle}\int_{S^{d-1}}\int_{\mathbb{R}}\int_{\mathbb{R}}\langle f, G_{m,t,u}\rangle\Psi_{m,t,u}\; dm\; dt\; du.\end{align*}
\end{thm}

Note that without the weights in the Fourier domain, the reconstruction formula does not hold.  Indeed, let $g_{m,t,u}(x)=g^{m,t}(u\cdot x)$, and $R^{*}$ the adjoint of the Radon transform:
\begin{align*}R^{*}g(x)=\int_{S^{d-1}}g(u,u\cdot x)du.\end{align*} We have the following result \cite{Grafakos_Sansing}.

\begin{thm} Let $f\in L^{1}(\mathbb{R}^{d}), \hat{f}\in L^{p}(\mathbb{R}^{d})$ for some $1<p<d$.  Then for $g, \psi\in \mathcal{S}(\mathbb{R})$ with $\langle g, \psi\rangle\neq 0$, the following identity holds:

\begin{align*} R^{*}(R(f))=\frac{1}{\langle g, \psi\rangle}\int_{S^{d-1}}\int_{\mathbb{R}}\int_{\mathbb{R}}\langle f, g_{m,t,u}\rangle \psi_{m,t,u}\; dm \; dt \; du.\end{align*}
\end{thm}
Thus, instead of reconstructing $f$ itself, we reconstruct the \emph{weighted back-projection} $R^{*}(R(f))$.  Indeed, one may compute that for a suitable constant $C_{d}$,  

\begin{align*}(-\Delta)^{\frac{d-1}{2}}R^{*}(R(f))=C_{d}f,\end{align*}so that we do not recover $f$ perfectly in general.  Here, the fractional power of the Laplacian $\Delta$ is understood in the sense of a pseudodifferential operator.

There is also a Parseval-type formula for Gabor ridge functions, Corollary 1 in \cite{Grafakos_Sansing}:

\begin{thm}\label{ridge_Parseval}

Let $f\in L^{1}(\mathbb{R}^{d})\cap L^{2}(\mathbb{R}^{d})$ and let $g\in\mathcal{S}(\mathbb{R})$ be not identically $0$.  There exists a constant $C_{g}$, depending only on $g$, such that

\begin{align*}\|f\|_{2}^{2}=C_{g}\int_{S^{d-1}}\int_{\mathbb{R}}\int_{\mathbb{R}}|\langle f, G_{m,t,u}\rangle|^{2}\; dm \; dt \; du.\end{align*}
\end{thm}

Note that Theorem \ref{ridge_Parseval} expresses the fact that the Gabor ridge system
\begin{align*}\{G_{m,t,u}\}_{(m,t,u)\in\mathbb{R}\times\mathbb{R}\times S^{d-1}}\end{align*}forms a continuous representation for $f\in L^{1}(\mathbb{R}^{d})\cap L^{2}(\mathbb{R}^{d})$.

In the case of Gabor ridge functions, the continuous system is indexed by the set $\mathbb{R}\times\mathbb{R}\times S^{d-1}$.  In general, this set does not admit a non-trivial group structure.  Indeed, the only spheres which admit non-trivial Lie structures are $S^{0}, S^{1}$ and $S^{3}$ \cite{Hofmann}.  Thus, the use of co-orbit theory to prove a discretization appears challenging.  We shall prove sufficient conditions for the discretization of Gabor ridge systems, but using techniques other than coorbit theory.  

As a starting point, consider the following \emph{semi-discrete} representing formula for Gabor ridge systems, which is Theorem 5 in \cite{Grafakos_Sansing}.

\begin{thm}\label{semi_discrete_thm}There exist $g, \psi\in\mathcal{S}(\mathbb{R})$ and $\alpha,\beta>0$ such that for all $f\in L^{1}(\mathbb{R}^{d})\cap L^{2}(\mathbb{R}^{d})$, we have have:

\begin{align}\label{semi_discrete_rep}A\|f\|_{2}^{2}\le \int_{S^{d-1}}\sum_{m\in\mathbb{Z}}\sum_{t\in\mathbb{Z}}|\langle f, G_{\alpha m, \beta t, u}\rangle|^{2}du\le B\|f\|_{2}^{2},\end{align}where the constants $A,B$ depend only on $g,\alpha,\beta$.  Moreover, for this choice of $g,\psi$, we have

\begin{align*}f=\frac{1}{2}\int_{S^{d-1}}\sum_{m\in\mathbb{Z}}\sum_{t\in\mathbb{Z}}\langle f, G_{\alpha m, \beta t, u}\rangle \Psi_{\alpha m, \beta t, u}du.\end{align*}
\end{thm}

The representation (\ref{semi_discrete_rep}) is called semi-discrete because there are two discrete sums, but also an integral over $S^{d-1}$.  We are interested in a full discretization of the Gabor ridge system.  Grafakos and Sansing did not present such a theory in \cite{Grafakos_Sansing}, though they suggested the co-orbit theory of Feichtinger and Gr\"{o}chenig might be used.  We pursue a different approach, based on the work of Hern\'{a}ndez, Labate, and Weiss \cite{HLW}.

\section{Finding a Space of Functions to Represent}\label{sec3}

We begin by noting that, once discretized, the frequency weights used in the weighted Gabor ridge construction are not needed.  Indeed, the frequency weights are needed in the continuous regime to perform a continuous change of variables to address the fact that $R^{*}(R(f))\neq f$ in general.  Moreover, when considering the discrete implementation of our construction, the weights will have the effect of imposing additional averaging of the experimental data.  Our target image class is textures, which will be adversely affected by such averaging.  Hence, there is also numerical significance attached to not considering the frequency weights at present.  Thus, we consider frame elements of the form 

\begin{align*}g_{m,t,u}(x):=g^{m,t}( u\cdot x), m,t\in\mathbb{R}, \ u\in S^{d-1}, x\in\mathbb{R}^{d}\end{align*} for a function $g:\mathbb{R}\rightarrow\mathbb{C}$.  We seek a discrete system of the form

\begin{align*}\{g_{m,t,u}\}_{(m,t,u)\in\Lambda},\end{align*}along with a space of functions for which this set will be a discrete frame.  The indexing set $\Lambda\subset\mathbb{R}\times\mathbb{R}\times S^{d-1}$ must be discrete with respect to the natural topology of $\mathbb{R}\times\mathbb{R}\times S^{d-1}$.  The nomenclature used to describe such a system shall be \emph{discrete directional Gabor frame}.

A first, somewhat naive, approach to developing a discrete directional Gabor frame is to start with the semi-discrete representation, Theorem \ref{semi_discrete_thm}, and choose a fixed, discrete set of points on the unit circle at which to sample.  This would essentially be replacing the integral in (\ref{semi_discrete_rep}) with a finite sum.  To investigate this approach, we consider systems of the form

\begin{align*} \{g_{m, t, u}\}_{(m,t,u)\in\Delta\times Q},\end{align*}where $\Delta\subset\mathbb{R}\times\mathbb{R}$ is an arbitrary uniformly discrete set, and $Q\subset S^{d-1}$ is an arbitrary finite set.  Note that since $S^{d-1}$ is compact, $Q$ is finite if and only if it is discrete.  We show that such a system cannot be a frame for any reasonable space of functions.  This proof is based on ideas found in Chapter 4 of \cite{Lax_Zalcman}.  

\begin{prop}\label{prop1}Let $Q\subset S^{d-1}$ be finite and let $\epsilon>0$.  Then there exists $f\in \mathcal{S}(\mathbb{R}^{d})$ with support in $B_{\epsilon}(0)$ such that $\|f\|_{2}=1$ and

\begin{align*}R_{u_{i}}(f)\equiv0, \ \forall u_{i}\in Q.\end{align*}

\begin{proof}Note that 

\begin{align*}\|R_{u_{i}}(f)\|_{2}^{2}=&\int_{\mathbb{R}}\left | \widehat{R_{u_{i}}(f)(\gamma)}\right|^{2}\;d\gamma \\=& \int_{\mathbb{R}}|\hat{f}(\gamma u_{i})|^{2}\;d\gamma.\end{align*}These equalities follow from Parseval's identity and the Fourier slice theorem, respectively.  Let $\ell_{1},...,\ell_{q}$ be the lines $\{\gamma u_{i}\}_{u_{i}\in Q, \gamma\in\mathbb{R}}$.  Let $P$ be a polynomial vanishing on these lines, and let $g\in\mathcal{S}(\mathbb{R}^{d})$ have support in $B_{\epsilon}(0)$.  Then $P(\xi)\hat{g}(\xi)\in\mathcal{S}(\mathbb{R}^{d})$ vanishes on each of the lines $\ell_{1},...,\ell_{q}$; in particular, 

\begin{align*}\int_{\ell_{i}}|P(\xi)\hat{g}(\xi)|^{2}\;d\xi=0, \ i=1,...,q.\end{align*}Let $D=(-i\frac{\partial}{\partial x_{1}},-i\frac{\partial}{\partial x_{2}},...,-i\frac{\partial}{\partial x_{d}})$ and define $f=C\cdot P(D)g$, where $C$ is a normalization constant so that $\|f\|_{2}=1$.  By construction, $f$ is supported in $B_{\epsilon}(0)$, and $\hat{f}=C'\cdot P(\xi)\hat{g}(\xi)$, so that 

\begin{align*}0=\int_{\ell_{i}}|P(\xi)\hat{g}(\xi)|^{2}\; d\xi=\int_{\ell_{i}}|\hat{f}(\xi)|^{2}\;d\xi=\int_{\mathbb{R}}|\hat{f}(\gamma u_{i})|^{2}\; d\gamma=\|R_{u_{i}}(f)\|_{2}^{2},\end{align*}where $\ell_{i}=\{\gamma u_{i}\}_{\gamma\in\mathbb{R}}$.  It follows that $R_{u_{i}}(f)\equiv 0$, as desired.

\end{proof}
\end{prop}

\begin{cor}Let $g\in L^{2}(\mathbb{R})$, let $\Delta\subset\mathbb{R}\times \mathbb{R}$ be an arbitrary discrete set, and let $Q\subset S^{d-1}$ be a fixed finite set.  Let $\epsilon >0$ be fixed.  Then the system

\begin{align*}\{g_{m, t, u}\}_{(m,t,u)\in \Delta\times Q}\end{align*}is not a frame for any subspace of $L^{2}(\mathbb{R}^{d})$ which contains $\{ f\in\mathcal{S}(\mathbb{R}^{d}) \ | \  \text{supp}(f)\subset B_{\epsilon}(0)\}$.  In fact, it cannot even be Bessel for such a space.    

\begin{proof}By Lemma 1 in \cite{Grafakos_Sansing}, $\langle f, g_{m, n, u}\rangle = \langle R_{u}(f), g^{m,n}\rangle$.  Then:

\begin{align*}&\sum_{u\in Q}\sum_{(m,n)\in\Delta}|\langle f, g_{m, n, u}\rangle|^{2}\\=&\sum_{u\in Q}\sum_{(m,n)\in\Delta}|\langle R_{u}(f), g^{m,n}\rangle|^{2}.\end{align*}  By Proposition \ref{prop1}, we may find $\tilde{f}\in\mathcal{S}(\mathbb{R}^{d})$, $\text{supp}(\tilde{f})\subset B_{\epsilon}(0), \|\tilde{f}\|_{2}=1$, such that 

\begin{align*}R_{u}(\tilde{f})\equiv0, \forall u\in Q.\end{align*}Noting that this forces 

\begin{align*}\sum_{u\in Q}\sum_{(m,n)\in\Delta}|\langle \tilde{f}, g_{m, n, u}\rangle|^{2}=0,\end{align*}the result is proved. 
\end{proof}

\end{cor}

Thus, infinitely many directions from $S^{d-1}$ must be included in any discrete directional Gabor frame representing Schwartz functions supported on some open set about the origin.  

This leads us to consider what space of functions can be represented with a discrete directional Gabor frame.  The following result demonstrates that we cannot have a discrete directional Gabor frame for any function space that contains $\mathcal{S}(\mathbb{R}^{d})$, regardless of sampling set.  

\begin{thm}\label{strongnegative}Let $\{g_{m,t,u}\}_{(m,t,u)\in \Lambda}$ be any discrete directional Gabor system.  Then there exists a sequence $\{\phi_{n}\}_{n=1}^{\infty}\subset\mathcal{S}(\mathbb{R}^{d})$ such that

\begin{enumerate}

\item $\lim\limits_{n\rightarrow\infty}\|\phi_{n}\|_{2}= 0.$

\item $\sum\limits_{(m,t,u)\in\Lambda}|\langle \phi_{n}, g_{m,t,u}\rangle|^{2}\ge 1, \forall n.$
\end{enumerate}

In particular, $\{g_{m,t,u}\}_{(m,t,u)\in\Lambda}$ cannot be a frame for any function space containing $\mathcal{S}(\mathbb{R}^{d})$; it cannot even be Bessel.  

\begin{proof}

We will show that $\phi_{n}$ may be chosen with the property that 

\begin{align}\label{line1}|\langle \phi_{n}, g_{m,t,u}\rangle|^{2}\ge 1,\end{align}for a fixed $(m,t,u)\in\Lambda$.  By rotating $\phi_{n}$ as needed, we may assume without loss of generality that $u=(1,0,...,0)\in S^{d-1}$.  Then we compute:

\begin{align*}|\langle \phi_{n}, g_{m,t,u}\rangle|^{2}=&|\langle R_{u}(\phi_{n}), g^{m,t}\rangle|\\=&\left|\int_{\mathbb{R}}\hat{\phi}_{n}(\gamma u)\overline{\widehat{g^{m,t}}(\gamma)}d\gamma\right|.\end{align*}Now, let us choose $\phi_{n}$ such that

\begin{align*}\hat{\phi_{n}}(\lambda)=\hat{\phi_{n}}(\lambda_{1},\lambda_{2},...,\lambda_{d})=\hat{\psi_{n}}(\lambda_{1})\cdot\hat{\eta_{n}}(\lambda_{2},...,\lambda_{d}),\end{align*}where $\hat{\psi}_{n}\in\mathcal{C}_{c}^{\infty}(\mathbb{R})$ is chosen such that $\|\hat{\psi_{n}}\|_{\infty}$ is uniformly bounded in $n$, and such that

\begin{align}\label{line2}\left|\int_{\mathbb{R}}\hat{\phi}_{n}(\gamma u)\overline{\widehat{g^{m,t}}(\gamma)}|d\gamma\right|=|\hat{\eta_{n}}(0,...,0)|\left|\int_{\mathbb{R}}\hat{\psi}_{n}(\gamma)\overline{\widehat{g^{m,t}}(\gamma)}d\gamma\right|=|\hat{\eta_{n}}(0,...,0)|.\end{align}The $\hat{\eta}_{n}\in\mathcal{C}_{c}^{\infty}(\mathbb{R}^{d-1})$ are chosen such that $|\hat{\eta_{n}}(0,...,0)|=1,$ for all $n$, and $$\lim_{n\rightarrow\infty}\|\hat{\eta}_{n}\|_{2}=0.$$  Note that with this choice, (\ref{line2}) implies (\ref{line1}), and that

\begin{align*}\lim_{n\rightarrow\infty}\|\phi_{n}\|_{2}=&\lim_{n\rightarrow\infty}\|\hat{\phi}_{n}\|_{2}\\=&\lim_{n\rightarrow\infty}\|\hat{\psi}_{n}\cdot\hat{\eta}_{n}\|_{2}\\\le&\lim_{n\rightarrow\infty}\|\hat{\psi}_{n}\|_{\infty}\|\hat{\eta}_{n}\|_{2}= 0,\end{align*}which gives the desired result.\end{proof}

\end{thm}

So, it is clear that in order to discretize a Gabor ridge system to acquire a discrete directional Gabor frame, infinitely many directions on the sphere must be incorporated, and our space of functions must be chosen very specifically.  Indeed, we must shrink the space of functions not to include all of $\mathcal{S}(\mathbb{R}^{d})$.  In the spirit of classical sampling, we will consider functions supported in a fixed compact domain.

\section{An Elementary Example}\label{sec4}

We show it is possible to produce a discrete directional Gabor system representing a space of functions with restricted support.  This example is based on manipulating directional Gabor systems to act like Fourier series.    

\begin{thm}\label{toy_example}Let $g(x):=\chi_{[-\frac{1}{2}, \frac{1}{2}]}(x), \ x\in\mathbb{R}$, the indicator function on $[-\frac{1}{2}, \frac{1}{2}]$.  Let $\Gamma\subset S^{d-1}\times\mathbb{R}$ be such that the mapping $\psi:\Gamma\rightarrow\mathbb{Z}^{d}$ given by $(u,m)\mapsto m u$ is a bijection.  Set 

\begin{align*}\Lambda=\{(m,n,u) \ | \ (u,m)\in\Gamma, n\in\mathbb{Z}\}.\end{align*}  Then we have

\begin{align*}\sum_{(m,n,u)\in\Lambda}|\langle f, g_{m,n,u}\rangle|^{2}=\|f\|_{2}^{2},\end{align*}for all $f \in L^2(\R^d)$ where ${\rm supp}(f) \subset B_{\frac{1}{2}}(0)$.

\begin{proof}

We compute:

\begin{align*}\sum_{(m,n,u)\in\Lambda}|\langle f, g_{m,n,u}\rangle|^{2}=&\sum_{(u,m)\in\Gamma}\sum_{n\in\mathbb{Z}}|\langle f, g_{m,n,u}\rangle|^{2}\\=& \sum_{(u,m)\in\Gamma}|\langle f, g_{m,0,u}\rangle|^{2}\\=&\sum_{(u, m)\in\Gamma}\left|\int_{\mathbb{R}^{d}}f(x)\chi_{[-\frac{1}{2}, \frac{1}{2}]}(u\cdot x)e^{-2\pi i m u\cdot x}\; dx\right|^{2}\\=&\sum_{(u, m)\in\Gamma}\left|\int_{B_{\frac{1}{2}}(0)}f(x)e^{-2\pi i m u\cdot x}\; dx\right|^{2}\\=&\sum_{k\in\mathbb{Z}^{d}}|\hat{f}(k)|^{2}\\=&\|f\|_{2}^{2}.\end{align*}

\end{proof}

\end{thm}

We remark that the penultimate line is simply the sum of the Fourier coefficients of $f$.  This example shows the existence of a tight discrete directional Gabor frame for the subspace of $L^2(\R^d)$ consisting of functions supported in $B_{\frac{1}{2}}(0)$.  Note that in this proof, we make no use of the translations in the system $\{g_{m,n,u}\}_{(m,n,u)\in\Lambda}$.  Indeed, the window $g$ is supported on $[-\frac{1}{2}, \frac{1}{2}]$, and we consider functions $f$ supported only on the ball $B_{\frac{1}{2}}(0)$.  It is of interest if this parameter set can be made non-trivial, to produce more interesting examples of discrete directional Gabor frames.

\section{Sufficient Conditions for a Discrete System}\label{sec5}

In this section, we develop sufficient conditions on the window function $g$ and sampling set $\Lambda\subset \mathbb{R}\times\mathbb{R}\times S^{d-1}$ for the discrete directional Gabor system 

\begin{align*}\{g_{m,n,u}\}_{(m,n,u)\in\Lambda}\end{align*}to be a frame.  We will consider a particular class of indexing sets, parametrized by real numbers $\omega>0$ that determine the localization of the window function $g$.  We will thus refer to our index sets as $\Lambda_{\omega}$, to denote this parametrization.

\begin{defn}Let $\mathcal{U}=\mathcal{U}_{d}$ be the subset of $\mathcal{S}(\mathbb{R}^{d})$ with support in $[-\frac{1}{2}, \frac{1}{2}]^{d}$.  \end{defn}

We prove sufficient conditions for the existence of discrete directional Gabor frames using tools developed by Hern\'{a}ndez, Labate, and Weiss \cite{HLW}.  This theory uses absolutely convergent almost periodic Fourier series to provide necessary and sufficient conditions for certain functions to generate discrete frames.  Among the classes of functions they characterize are classical discrete Gabor and discrete wavelet systems.  In the process of developing sufficient conditions for discrete directional Gabor frames, we show their methods are partially extensible to certain anisotropic systems.

Recall that $g_{m,n,u}(x):=g^{m,n}(u\cdot x), m,n\in\mathbb{R}, \ u\in S^{d-1}, x\in\mathbb{R}^{d}$.

\begin{lem}\label{Lemma1}Let $f\in\mathcal{S}(\mathbb{R}^{d})$, $g\in L^{2}(\mathbb{R})$.  Then for a fixed $(u,m)\in S^{d-1}\times\mathbb{R},$

\begin{align*}\sum_{n\in \myz \mathbb{Z}}|\langle f, g_{m,  n,u}\rangle|^{2}= \frac{1}{\myz} \int_{\mathbb{R}}\sum_{k\in  \mathbb{Z}/\myz}\hat{f}(\gamma u)\overline{\hat{g}(\gamma-m)}\overline{\hat{f}\left(\left(\gamma+k\right)u\right)}\hat{g}\left(\gamma+k-m\right)d\gamma.\end{align*}

\begin{proof}By Lemma 1 in \cite{Grafakos_Sansing} and the Fourier slice theorem, 

\begin{align*}\langle f, g_{m,n,u}\rangle=&\langle R_{u}f, g^{m,n}\rangle\\=&\int_{\mathbb{R}}\hat{f}(\gamma u)\overline{\widehat{g^{m,n}(\gamma)}}\; d\gamma\\=&\int_{\mathbb{R}}\hat{f}(\gamma u)\overline{\hat{g}(\gamma-m)}e^{2\pi in\gamma}\; d\gamma.\end{align*}Periodizing this integral, we see

\begin{align*}\langle f, g_{m,n,u}\rangle=&\sum_{k\in \mathbb{Z}/\myz}\int_{0}^{1/\myz}\hat{f}\left(\left(\gamma+k\right)u\right)\overline{\hat{g}\left(\gamma-m+k\right)}e^{2\pi in\gamma}\; d\gamma\\=&\int_{0}^{1/\myz}\sum_{k\in \mathbb{Z}/\myz}\hat{f}\left(\left(\gamma+k\right)u\right)\overline{\hat{g}\left(\gamma-m+k\right)}e^{2\pi in\gamma}\; d\gamma.\end{align*}Thus,

\begin{align*}&\sum_{n\in \myz \mathbb{Z}}|\langle f, g_{m,n,u}\rangle|^{2}\\=&\sum_{n\in \myz \mathbb{Z}}\left|\int_{0}^{1/\myz}\left (\sum_{k\in  \mathbb{Z}/\myz}\hat{f}\left(\left(\gamma+k\right)u\right)\overline{\hat{g}\left(\gamma-m+k\right)}\right)e^{2\pi in\gamma}\; d\gamma \right|^{2}\\=& \frac{1}{\myz} \int_{0}^{1/\myz}\left|\sum_{k\in \mathbb{Z}/ \myz}\hat{f}\left(\left(\gamma+k\right)u\right)\overline{\hat{g}\left(\gamma-m+k\right)}\right|^{2}\; d\gamma.\end{align*}To transition from the penultimate to the ultimate line, we applied Parseval's formula to the function

\begin{align*}\gamma\mapsto\sum_{k\in \mathbb{Z}/\myz}\hat{f}\left(\left(\gamma+k\right)u\right)\overline{\hat{g}\left(\gamma-m+k\right)}.\end{align*}

Now,
\begin{align*} & \int_{0}^{1/\myz}\left|\sum_{n\in \mathbb{Z}/\myz}\hat{f}\left(\left(\gamma+n\right)u\right)\overline{\hat{g}(\gamma+n-m)}\right|^{2}d\gamma\\=&\int_{0}^{1/\myz}\sum_{j,\ell \in \mathbb{Z}/\myz}\hat{f}\left(\left(\gamma+j\right)u\right)\overline{\hat{g}\left(\gamma+j-m\right)}\overline{\hat{f}\left(\left(\gamma+\ell\right)u\right)}\hat{g}\left(\gamma+\ell-m\right)d\gamma\\=&\int_{0}^{1/\myz}\sum_{k,j\in \mathbb{Z}/\myz}\hat{f}\left(\left(\gamma+j\right)u\right)\overline{\hat{g}\left(\gamma+j-m\right)}\overline{\hat{f}\left(\left(\gamma+j+k\right)u\right)}\hat{g}\left(\gamma+j+k-m\right)d\gamma\\=&\int_{\mathbb{R}}\sum_{k\in  \mathbb{Z}/\omega}\hat{f}(\gamma u)\overline{\hat{g}(\gamma-m)}\overline{\hat{f}\left(\left(\gamma+k\right)u\right)}\hat{g}\left(\gamma+k-m\right)d\gamma.\end{align*}

To go from the penultimate to ultimate line, we sum over $j$ and exploit the periodization.  Since $f$ is Schwartz class, we may interchange the integral and sum.  This gives the result.

\end{proof}

\end{lem}

Let $\Gamma\subset S^{d-1}\times\mathbb{R}$ be as in Theorem \ref{toy_example}, i.e., let $\Gamma\subset S^{d-1}\times\mathbb{R}$ be such that the mapping $\psi:\Gamma\rightarrow\mathbb{Z}^{d}$ given by $(u,m)\mapsto m u$ is a bijection.  Let $\Lambda_\myz =\{(m,n,u)\ | \ (u,m)\in\Gamma, n\in \myz \mathbb{Z}\}$ where $\myz > 0$.  While many such sets exists, one simple construction is 

\begin{align*} \Gamma=\left\{\left(\frac{a}{\|a\|_{2}}, \|a\|_{2}\right)\right\}_{a\in\mathbb{Z}^{d}\setminus\{0\}}\cup\{(u_{0},0)\},\end{align*}for arbitrary $u_{0}\in S^{d-1}$.  Note that by construction, $\Gamma \neq A\times B$ for $A\subset S^{d-1}, B\subset\mathbb{R}$, which implies $\Lambda_{\omega}$ also does not split as a direct product.  Moreover, with $\Gamma$ defined generally as above, the sampling set $\Lambda_{\myz}$ is indeed discrete. 

\begin{lem}
The set $\Lambda_{\myz}$ is a discrete subset of $\mathbb{R}\times\mathbb{R}\times S^{d-1}$.

\begin{proof}
It suffices to show that $\Gamma$ is a discrete subset of $S^{d-1} \times \mathbb{R}$. Let $(u, m) \in \Gamma$ and $\alpha_u$ be the smallest positive real number so that $\alpha_u u \in \mathbb{Z}^{d}$. First suppose $m \ne 0$. Let $N_m = \{ v \in S^{d-1} : \alpha_v < |m| +1\}$. Notice that $N_m$ is finite and $u \in N_m$. So we can choose a neighborhood $u\in U \subset S^{d-1}$ so that $U \cap N_m = \{ u \}$. Let $V = U \times (m-1, m+1 )$. Now, suppose $(v,l) \in V \cap \Gamma$. This implies $\alpha_v < |m|+1$ and so $v \in U \cap N_m$. Thus $v = u$. Since $|m-l| \ge 1$ when $m \ne l$, then we must have $l = m$. Hence $V \cap \Gamma = \{ (u,m) \}$.

It remains to prove the case when $m = 0$. Since $(u,m) \mapsto mu$ is a bijection between $\Gamma$ and $\Z^d$, there is precisely one $u_{0} \in S^{d-1}$ so that $(u_{0},0)  \in \Gamma$.  We set $V = S^{d-1} \times (-1,1)$. Since every $(v,\ell) \in \Gamma$ satisfies $|\ell| \ge 1$ except when $(u,\ell) = (u_{0},0)$, it follows that $V \cap \Gamma = \{(u_{0},0)\}$. Therefore, $\Gamma$ is discrete.

\end{proof}

\end{lem}

We are now prepared to give a sufficient condition for a discrete directional Gabor frame.  A key ingredient in the proof is a multidimensional version of Kadec's $\frac{1}{4}$-Theorem \cite{2D_Kadec}.

\begin{thm}\label{Kadec1/4}Let $\Omega=\{\lambda_{k}=(\lambda_{k_{1}},...,\lambda_{k_{d}})\ | \ k\in\mathbb{Z}^{d}\}\subset\mathbb{R}^{d}$ and suppose there exists $L<\frac{1}{4}$ such for all $k\in\mathbb{Z}^{d}$, the associated $\lambda_{k}$ satisfies 
\begin{align*}\|k-\lambda_{k}\|_{\infty} \le L.\end{align*}Then the collection $\{e^{-2\pi i \lambda\cdot x}\}_{\lambda\in\Omega}$ is a Riesz basis for $L^{2}([-\frac{1}{2}, \frac{1}{2}]^{d})$, with frame bounds $A=(2\pi)^{d}(\cos(\pi L)-\sin(\pi L))^{2d}, B=(2\pi)^{d}(2-\cos(\pi L)+\sin(\pi L))^{2d}$.\end{thm}

We also require a result on perturbing Fourier frames to acquire Bessel systems.  The following result appears in \cite{Duffin_Schaeffer} in the case $d=1$; we require the result for arbitrary $d$.  The result is probably known, but we give a proof for completeness. 

\begin{thm}\label{Bessel} Suppose $\{\lambda_{k}\}_{k\in\mathbb{Z}^{d}}\subset\mathbb{R}^{d}$ is such that there exists $0<A\le B<\infty$ such that for all $f\in\mathcal{U}$,

\begin{align*}A\|f\|_{2}^{2}\le\sum_{k\in\mathbb{Z}^{d}}|\hat{f}(\lambda_{k})|^{2}\le B\|f\|_{2}^{2}.\end{align*}Let $\{\mu_{k}\}_{k\in\mathbb{Z}^{d}}\subset\mathbb{R}^{d}$ and $M>0$ be such that $\|\lambda_{k}-\mu_{k}\|_{\infty}<M, \forall k\in\mathbb{Z}^{d}$.  Then 

\begin{align*}\sum_{k\in\mathbb{Z}^{d}}|\hat{f}(\mu_{k})|^{2}\le B(1+\sqrt{B'})^{2}\|f\|_{2}^{2},\end{align*}where $B'=\frac{B}{A}\left(e^{M^{2}\rho^2d}-1\right)\left(e^{\frac{\pi^{2}d}{\rho^{2}}}-1\right)$ for any $\rho>0$.  In particular, if $A=B=\rho=1$, then $B'=\left(e^{M^{2}d}-1\right)\left(e^{\pi^{2}d}-1\right)$.  

\begin{proof}Since $f\in\mathcal{U}$, $f$ is supported in $[-\frac{1}{2}, \frac{1}{2}]^{d}$.  So, by Schwartz's Paley-Wiener Theorem \cite{Hormander}, $\hat{f}$ is analytic on $\mathbb{C}^{d}$.  Consider its power series expansion about a fixed $\lambda_{k}:$

\begin{align*}\hat{f}(\gamma)=\sum_{\alpha\in\mathbb{N}^{d}}\frac{\partial^{\alpha}\hat{f}(\lambda_{k})}{\alpha!}(\gamma-\lambda_{k})^{\alpha},\end{align*}where $\alpha=(\alpha_{1},...,\alpha_{d})\in\mathbb{N}^{d}$ is a multi-index.  Then for any $\rho>0$, an application of the Cauchy-Schwarz inequality yields:

\begin{align*}|\hat{f}(\mu_{k})-\hat{f}(\lambda_{k})|^{2}=&\left|\sum_{\alpha\neq(0,...,0)}\frac{\partial^{\alpha}\hat{f}(\lambda_{k})}{\alpha!}(\mu_{k}-\lambda_{k})^{\alpha}\right|^{2}\\ \le &\sum_{\alpha\neq(0,...,0)}\frac{|\partial^{\alpha}\hat{f}(\lambda_{k})|^{2}}{\rho^{2|\alpha|}\alpha!}\cdot \sum_{\alpha\neq(0,...,0)}\frac{M^{2|\alpha|}\rho^{2|\alpha|}}{\alpha!}\\=&\left(e^{M^{2}\rho^{2}d}-1\right)\cdot\sum_{\alpha\neq(0,...,0)}\frac{|\partial^{\alpha}\hat{f}(\lambda_{k})|^{2}}{\rho^{2|\alpha|}\alpha!}.\end{align*}

Recalling the correspondence between polynomial multiplication and the Fourier transform, $\partial^{\alpha}\hat{f}(\gamma)=\left(\frac{2\pi}{i}\right)^{|\alpha|}\left(x^{\alpha}f(x)\right)^{\hat{}}(\gamma)$, and thus $(\partial^{\alpha}\hat{f})^{\check{}}\in\mathcal{U}$.  We may thus sum over $k\in\mathbb{Z}^{d}$ as follows: \begin{align*}\sum_{k\in\mathbb{Z}^{d}}|\hat{f}(\mu_{k})-\hat{f}(\lambda_{k})|^{2}\le&\sum_{k\in\mathbb{Z}^{d}}\left(e^{M^{2}\rho^{2}d}-1\right)\cdot\sum_{\alpha\neq(0,...,0)}\frac{|\partial^{\alpha}\hat{f}(\lambda_{k})|^{2}}{\rho^{2|\alpha|}\alpha!}\\=&\left(e^{M^{2}\rho^{2}d}-1\right)\sum_{\alpha\neq(0,...,0)}\frac{1}{\rho^{2|\alpha|}\alpha!}\sum_{k\in\mathbb{Z}^{d}}|\partial^{\alpha}\hat{f}(\lambda_{k})|^{2}\\ \le& B\left(e^{M^{2}\rho^{2}d}-1\right)\sum_{\alpha\neq(0,...,0)}\frac{1}{\rho^{2|\alpha|}\alpha!}\|\partial^{\alpha}\hat{f}\|_{2}^{2}\\ =&B\left(e^{M^{2}\rho^{2}d}-1\right)\sum_{\alpha\neq(0,...,0)}\frac{(2\pi)^{2|\alpha|}}{\rho^{2|\alpha|}\alpha!}\int_{\mathbb{R}^{d}}|x|^{2\alpha}|f(x)|^{2}dx \\ \le&B\left(e^{M^{2}\rho^{2}d}-1\right)\sum_{\alpha\neq(0,...,0)}\frac{(2\pi)^{2|\alpha|}}{\rho^{2|\alpha|}\alpha!}\frac{1}{2^{2|\alpha|}}\|f\|_{2}^{2}\\ \le&\frac{B}{A}\left(e^{M^{2}\rho^{2}d}-1\right)\left(e^{\frac{\pi^{2}d}{\rho^{2}}}-1\right)\sum_{k\in\mathbb{Z}^{d}}|\hat{f}(\lambda_{k})|^{2}.\end{align*}

Thus, 

\begin{align*}\sum_{k\in\mathbb{Z}^{d}}|\hat{f}(\mu_{k})-\hat{f}(\lambda_{k})|^{2}\le B'\sum_{k\in\mathbb{Z}^{d}}|\hat{f}(\lambda_{k})|^{2},\end{align*}so that

\begin{align*}\sum_{k\in\mathbb{Z}^{d}}|\hat{f}(\mu_{k})|^{2}\le(1+\sqrt{B'})^{2}\sum_{k\in\mathbb{Z}^{d}}|\hat{f}(\lambda_{k})|^{2}\le B(1+\sqrt{B'})^{2}\|f\|_{2}^{2},\end{align*}as desired.

\end{proof}

\end{thm}

\begin{thm}\label{main_result}Let $g\in L^{2}(\mathbb{R})$ be such that $\hat{g}$ is compactly supported and such that for $k \in (\Z/\myz)\setminus \{0\}, \  \hat{g}(\gamma)\hat{g}(\gamma+k) = 0$ almost everywhere.  Furthermore, suppose it is not the case that $\hat{g}$ is zero almost everywhere on the interval  $[-1/4,1/4]$. Then with $\Lambda_\myz$ as above, $\{g_{m,n,u}\}_{(m,n,u)\in\Lambda_\myz}$ is a discrete frame for the subspace of functions in $L^2(\mathbb{R}^d)$ with support contained in $[-1/2,1/2]^d$.

\begin{proof}

Suppose $ g \in L^2(\R)$ satisfies the above hypotheses. By Lemma \ref{Lemma1}, if $f \in \mathcal{U}$, then
\begin{align*}
& \sum_{(u,m) \in \Gamma} \sum_{n \in \myz \Z} | \langle f,  g_{m,n,u} \rangle |^2 \\ &= \frac{1}{\myz} \sum_{(u,m) \in \Gamma} \int_{\mathbb{R}}\sum_{k\in \mathbb{Z}/\myz}\hat{f}(\gamma u)\overline{\hat{g}(\gamma-m)}\overline{\hat{f}\left(\left(\gamma+k\right)u\right)}\hat{g}\left(\gamma+k-m\right) \; d\gamma \\ &= \frac{1}{\myz} \sum_{(u,m) \in \Gamma} \int_{\mathbb{R}} |\hat{f}((\gamma  + m) u)|^2 |\hat{g}(\gamma) |^2  \; d\gamma 
\end{align*}

Let $a(\gamma) = (2\pi)^{2d}(\cos(\pi|\gamma|)-\sin(\pi|\gamma|))^{2d}$. Since $a(\gamma) > 0$ for $\gamma \in (-1/4,1/4)$, then 

\begin{align*}A = \frac{1}{\myz} \int_{-\frac{1}{4}}^{\frac{1}{4}} a(\gamma) |\hat{g}(\gamma)|^2 \; d\gamma > 0.\end{align*} For each $(u,m) \in \Gamma$, $\| (\gamma + m)u - m u\|_2 = |\gamma|$ and thus for $\gamma \in (-1/4,1/4)$, we can apply Theorem \ref{Kadec1/4} to obtain
\[
A\| f \|_2^2 = \frac{1}{\myz} \int_{-\frac{1}{4}}^{\frac{1}{4}} a(\gamma) \| f \|_2^2 |\hat{g}(\gamma)|^2\;d\gamma\le \frac{1}{\myz} \sum_{(u,m) \in \Gamma} \int_{\mathbb{R}} |\hat{f}((\gamma  + m) u)|^2 |\hat{g}(\gamma) |^2  \; d\gamma.
\]
This establishes the lower frame bound.

To establish the upper frame bound, for each $\gamma$ we apply Theorem \ref{Bessel} to the sequence of points $\{ (\gamma + m) u \}_{(u,m) \in \Gamma}$. Set 
\[
b(\gamma) = (1+\sqrt{c(\gamma)})^2 \quad \mbox{where} \quad c(\gamma) =\left(e^{\gamma^{2}d}-1\right)\left(e^{\pi^{2}d}-1\right).
\]
Since $\hat{g}$ is compactly supported, we may choose $M > 0$ so that the support of $\hat{g}$ is contained in $[-M,M]$. We let
\[
B = \frac{1}{\myz} \int_{-M}^{M}  b(\gamma)|\hat{g}(\gamma)|^2 \; d \gamma.
\]
Then by Theorem \ref{Bessel}, we have
 \begin{align*}& \frac{1}{\myz} \sum_{(u,m) \in \Gamma} \int_{\mathbb{R}} |\hat{f}((\gamma  + m) u)|^2 |\hat{g}(\gamma) |^2  \; d\gamma \\ =& \frac{1}{\myz} \int_{\mathbb{R}} \sum_{(u,m) \in \Gamma}|\hat{f}((\gamma  + m) u)|^2 |\hat{g}(\gamma) |^2  \; d\gamma \\ \le& \frac{1}{\myz} \int_{-M}^{M}b(\gamma)|\hat{g}(\gamma)|^{2}\|f\|_{2}^{2}d\gamma\\ =& B \| f \|_2^2.\end{align*}
 
This establishes both frame bounds for $\{g_{m,n,u}\}_{(m,n,u)\in\Lambda_\myz}$ for functions $f \in \mathcal{U}$. The result follows since $\mathcal{U}$ is dense in the set of $L^2$ functions supported in $[-\frac{1}{2},\frac{1}{2}]^d$.
\end{proof}

\end{thm}

We remark that the condition that $\hat{g}$ be compactly supported is not strictly necessary.  As long as $\hat{g}$ decays fast enough so that 

\begin{align*}\int_{-\infty}^{\infty}b(\gamma)|\hat{g}(\gamma)|^{2}\; d\gamma<\infty,\end{align*}the upper bound still holds.  Of course, $b(\gamma)$ grows like $e^{d\gamma^{2}}$, so $\hat{g}(\gamma)$ would need to decay extremely rapidly.  

Theorem \ref{toy_example} and Theorem \ref{main_result} suggest that the windows $g(x)=\chi_{[-\frac{1}{2}, \frac{1}{2}]}$, and $g(x)=\text{sinc}\left(\frac{x}{16}\right)^4$ can both generate discrete directional Gabor systems, with appropriately chosen sampling set $\Lambda_{\omega}$.  This is interesting since one is compactly supported in time, the other in frequency.  In our numerical experiments, we shall consider only a window that is compactly supported in frequency.  This is in order to match the typical types of window functions used in state-of-the-art anisotropic frames, such as shearlets and curvelets.  

We note that it is perhaps possible to acquire tight frames, and to extend the class of functions for which directional Gabor systems are a frame, by considering frequency weightings, as in \cite{Grafakos_Sansing}.  This is not pursued in the present article, particularly for the detrimental averaging effect such weighting may have on the numerical implementation, especially for textural images.  However, it remains an interesting question.

\section{Numerical Implementation}\label{sec6}

In this section, we consider a prototype MATLAB implementation of our discrete directional Gabor frame in two dimensions.  This algorithm decomposes an image into its directional Gabor coefficients, as determined by a selected window function $g$ and sampling set $\Lambda_{\omega}$.  In our case, the inverse frame operator is approximated with MATLAB's standard preconditioned conjugate gradient routine, with error tolerance $10^{-6}$.  We perform experiments with the window function:

\begin{align*} g(x)=\text{sinc}\left(\frac{x}{16}\right)^4,\end{align*} which has compact Fourier support $K=[-\frac{1}{4}, \frac{1}{4}]$.  This is suggested by Theorem \ref{main_result}.  Other windows could be used in order to optimize numerical performance, but we do not perform such an optimization in the present article.  

In addition, to create a numerical implementation of our discrete directional Gabor frame, we must choose an appropriate finite 

\begin{align*}\tilde{\Lambda} \subset \Lambda_\omega =\{(m,n,u)\ | \ (u,m)\in\Gamma, n\in \myz \mathbb{Z}\},\end{align*}with $\Gamma$ as in Theorem \ref{toy_example}.  While the choice of the $n$ coordinate is flexible, it is less clear how to choose an appropriate finite collection $\{(u,m) : (u,m) \subset \mathbb{S}^1 \times \mathbb{R}\}$ such that $(u,m) \mapsto mu$ is onto some subset of $\mathbb{Z}^2$. 

From an implementation standpoint, it is more useful to pick some $A \subset \mathbb{Z}^2$ finite and choose $(u,m)$ such that $(u,m) \mapsto mu$ is onto $A$. In the experiments which follow, for a $2M \times 2N$ image we choose 

\begin{align*}&A = \left([-M, M] \times [-N, N]\right) \cap \mathbb{Z}^2, \\ &u \in \left\{ \frac{(x_1, x_2)}{\|(x_1, x_2)\|_{2}} \ : \ (x_1, x_2) \in A, \ \gcd(x_1, x_2) = 1 \right \}.\end{align*}Note that $m$ is fixed by the choices of $A$ and $u$.  In our implementation we store the collection 

\begin{align*}D = \{ (x_1, x_2) : (x_1, x_2) \in A, \ \gcd(x_1, x_2) = 1 \}.\end{align*} instead of normalizing. Because the components of the vectors in $D$ are coprime, the choices for $m$ fixed by our choice of $A$ are all integers.

In testing our algorithm, we consider two choices of parameter set $\Lambda_{\omega}$.  The first case takes $\Gamma$ as above, and uses only $n=0$ for the translation parameter.  This produces a non-redundant Gabor system, which is essentially a directionally weighted Fourier transform.  We also consider a true, redundant directional Gabor system, in which the translation parameter runs over the set $n=4k, |k|\le 5$ for integers $k$.  This corresponds to taking $\omega=4$.  

We also consider experiments using several methods of state-of-the-art fast anisotropic representations: shearlets \cite{shearlet_acha, shearlet_edge} and curvelets \cite{acha_curvelets1, acha_curvelets2}.  These methods are known to be theoretically optimal for sparse representation among a large class of image-like signals \cite{shearlets_book}.  We consider both the ShearLab implementation of shearlets \cite{ShearLab} and the FFST implementation of shearlets \cite{FFST}.  We consider the CurveLab implementation of curvelets \cite{CurveLab}.  

We note that our redundant directional Gabor frame has redundancy on the order of the redundancy of shearlets and curvelets.  Indeed, the redundant transform we have implemented has approximately 11-fold redundancy, compared to approximately 3 for curvelets, exactly 29 for FFST, and 49 for Shearlab.  In order to also make a comparison with a non-redundant system, experiments with the undecimated wavelet transform with Haar wavelet are also presented.  Other wavelets, including biorthogonal and Daubechies wavelets with many vanishing moments, were considered as well.  The use of such wavelets did not in general improve on the Haar wavelet for the present experiments.

\subsection{Experimental Data}

\begin{figure}[H]
\centering
\includegraphics[width=.25\textwidth]{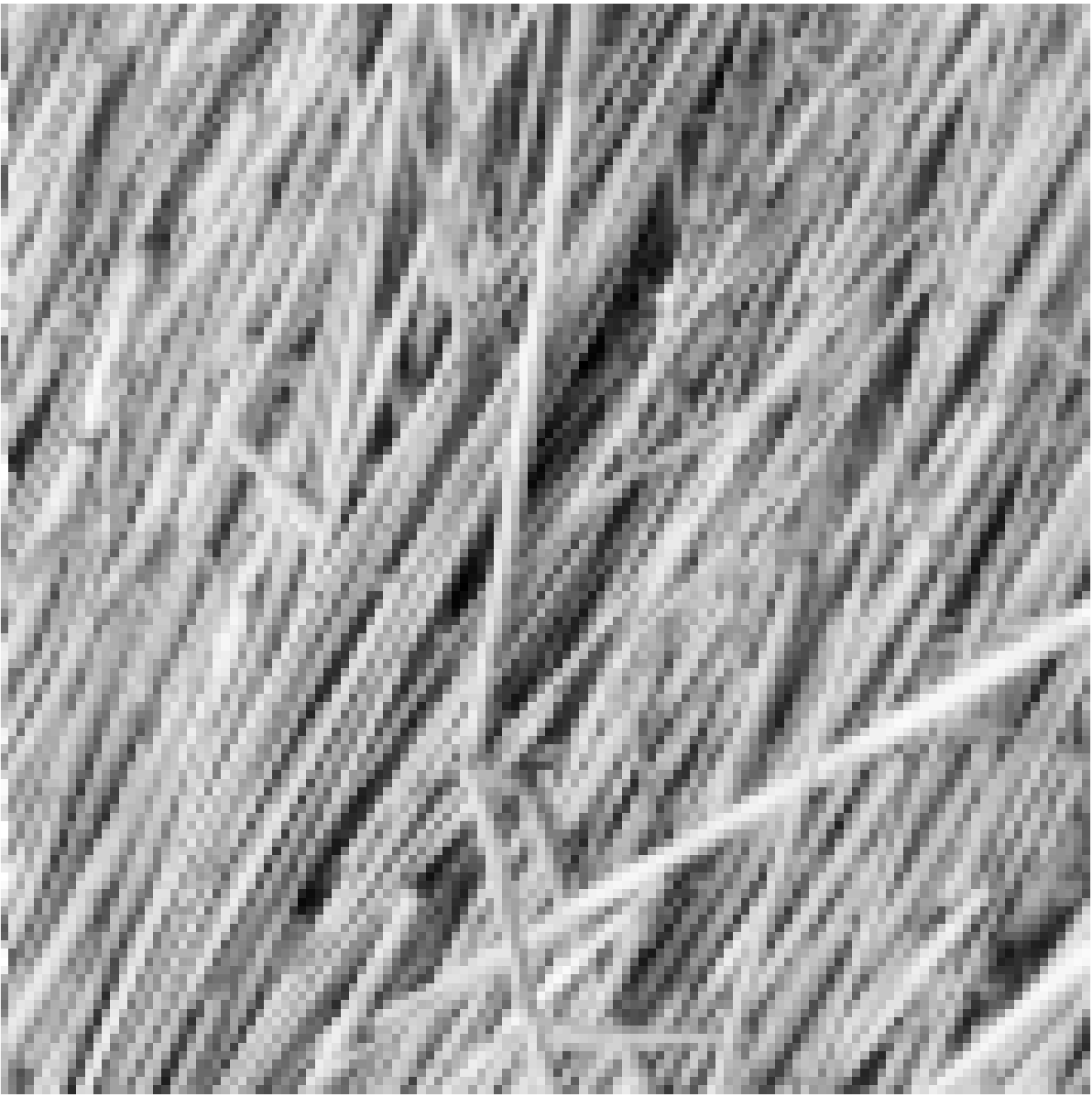}
\includegraphics[width=.25\textwidth]{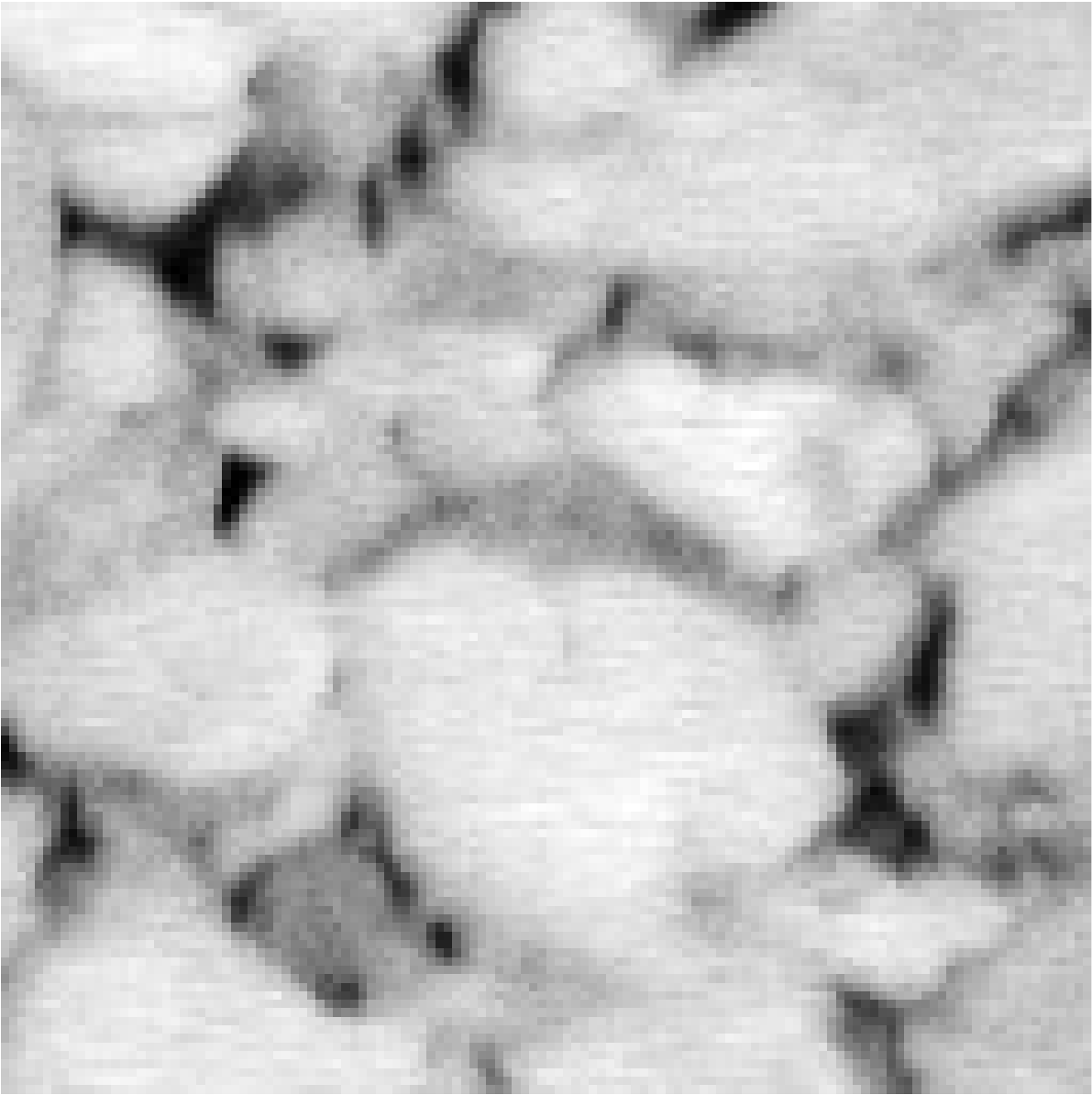}
\includegraphics[width=.25\textwidth]{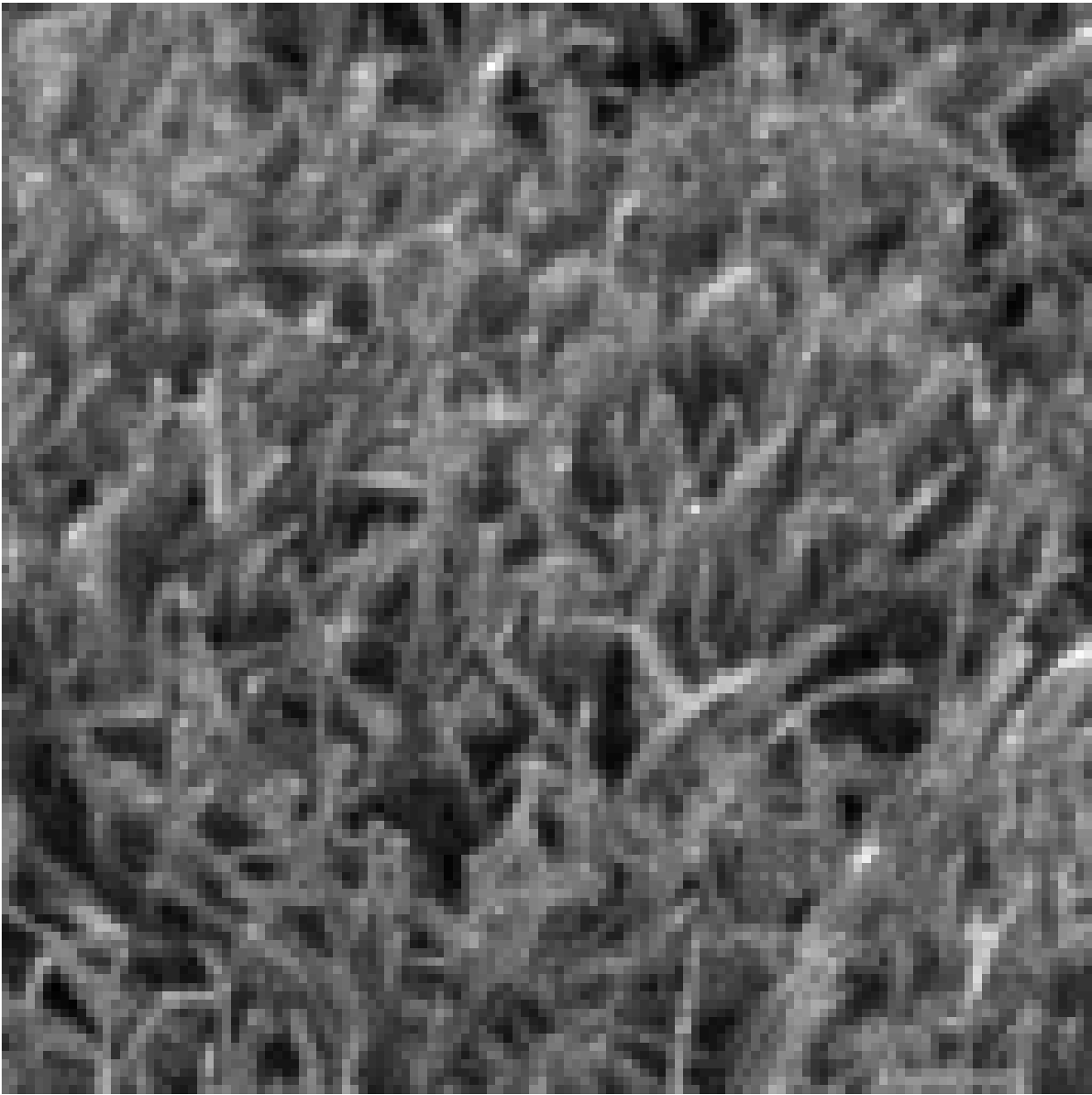}
\includegraphics[width=.25\textwidth]{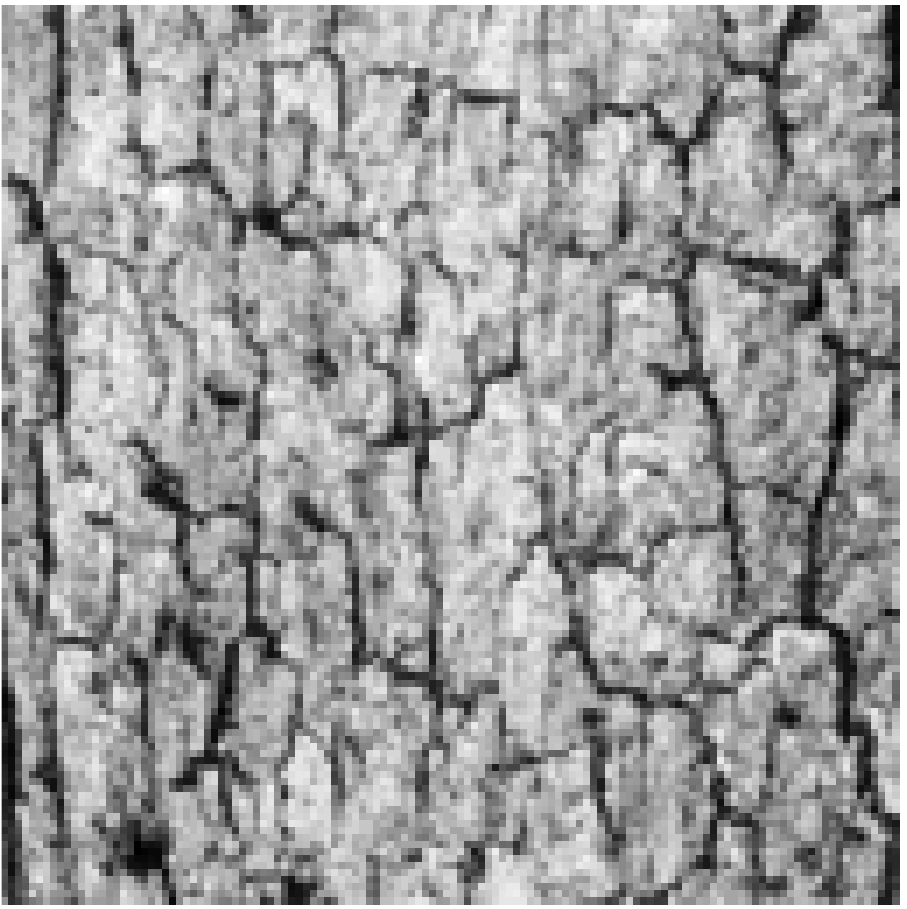}
\includegraphics[width=.25\textwidth]{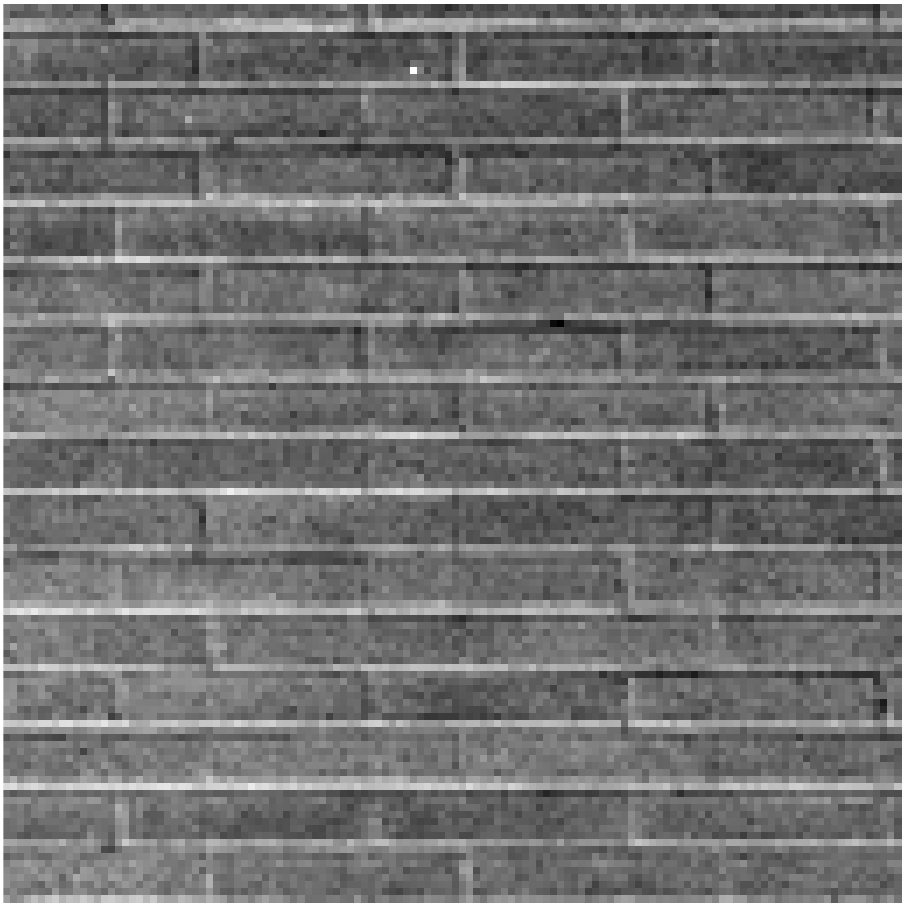}
\includegraphics[width=.25\textwidth]{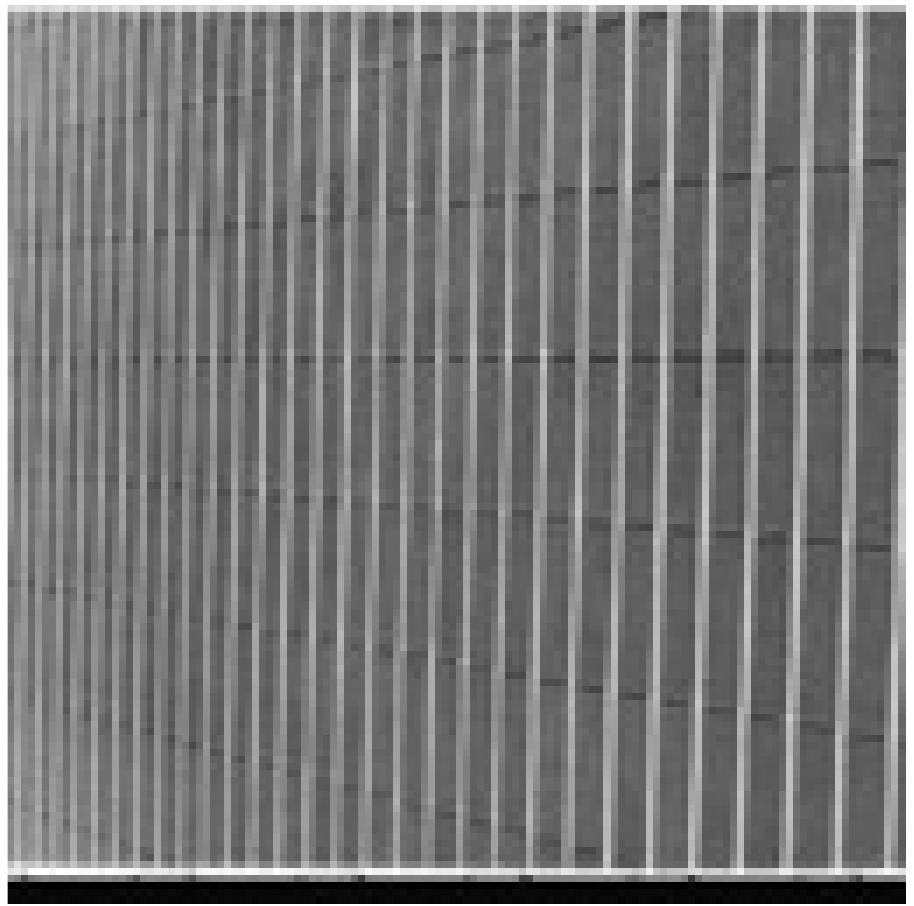}
\includegraphics[width=.25\textwidth]{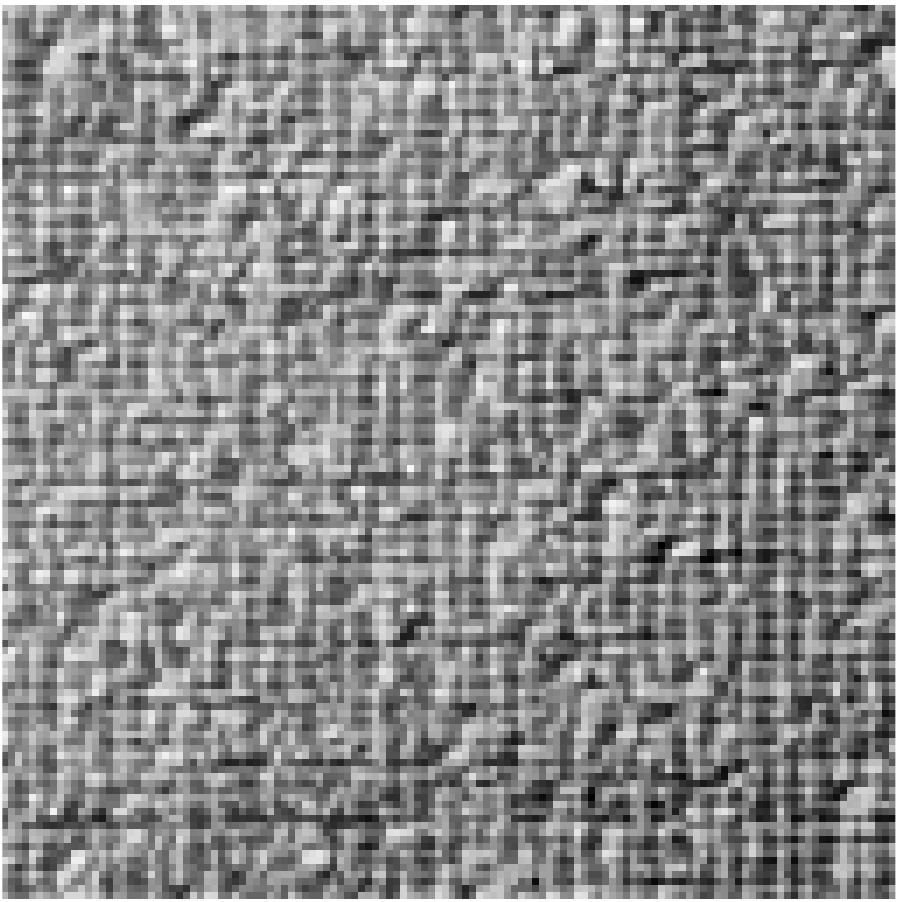}
\includegraphics[width=.25\textwidth]{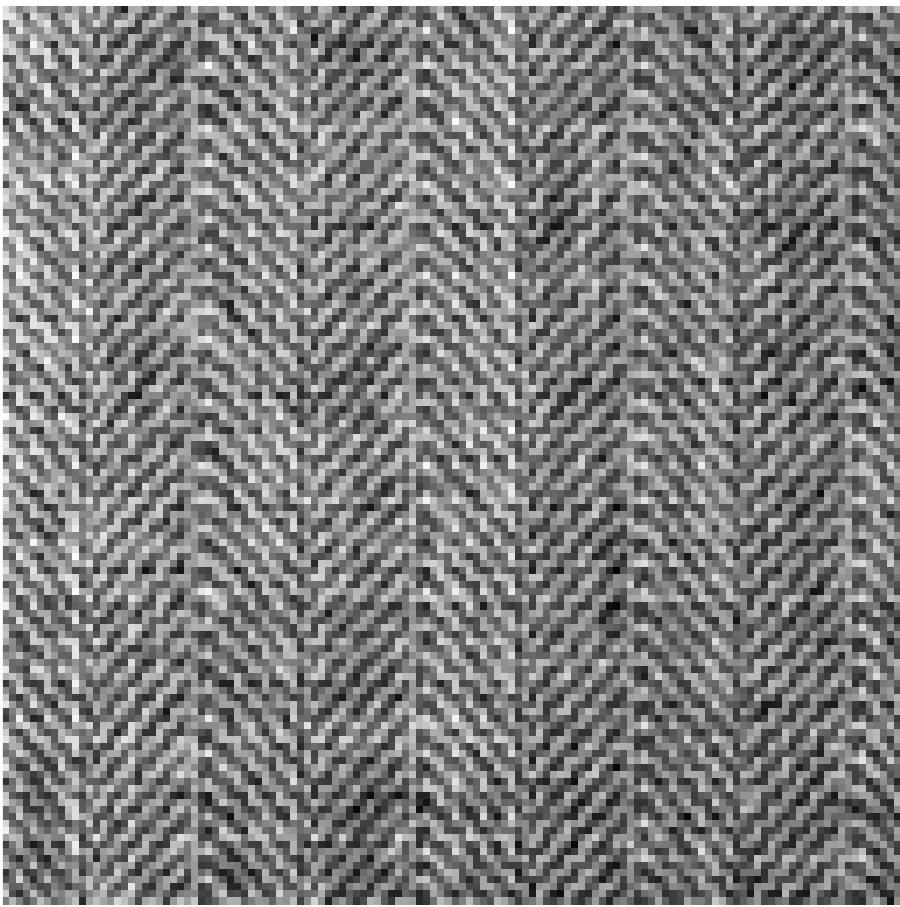}
\includegraphics[width=.25\textwidth]{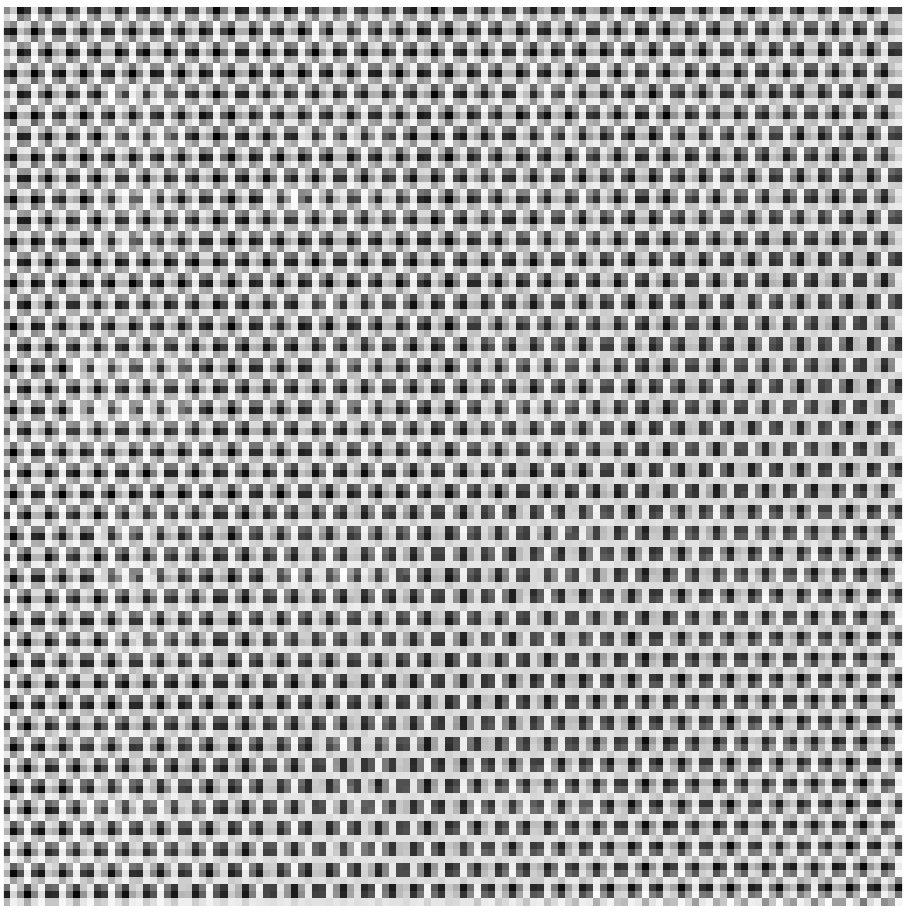}
\caption{Texture images for numerical experiments.  Top row, from left to right: straw, rocks, grass.  Middle row, from left to right: cracked mud, bricks, bars.  Bottom row, from left to right: fabric, grate, honeycomb.}
\label{texture_images}
\end{figure}

\begin{figure}[H]
\centering
\includegraphics[width=.25\textwidth]{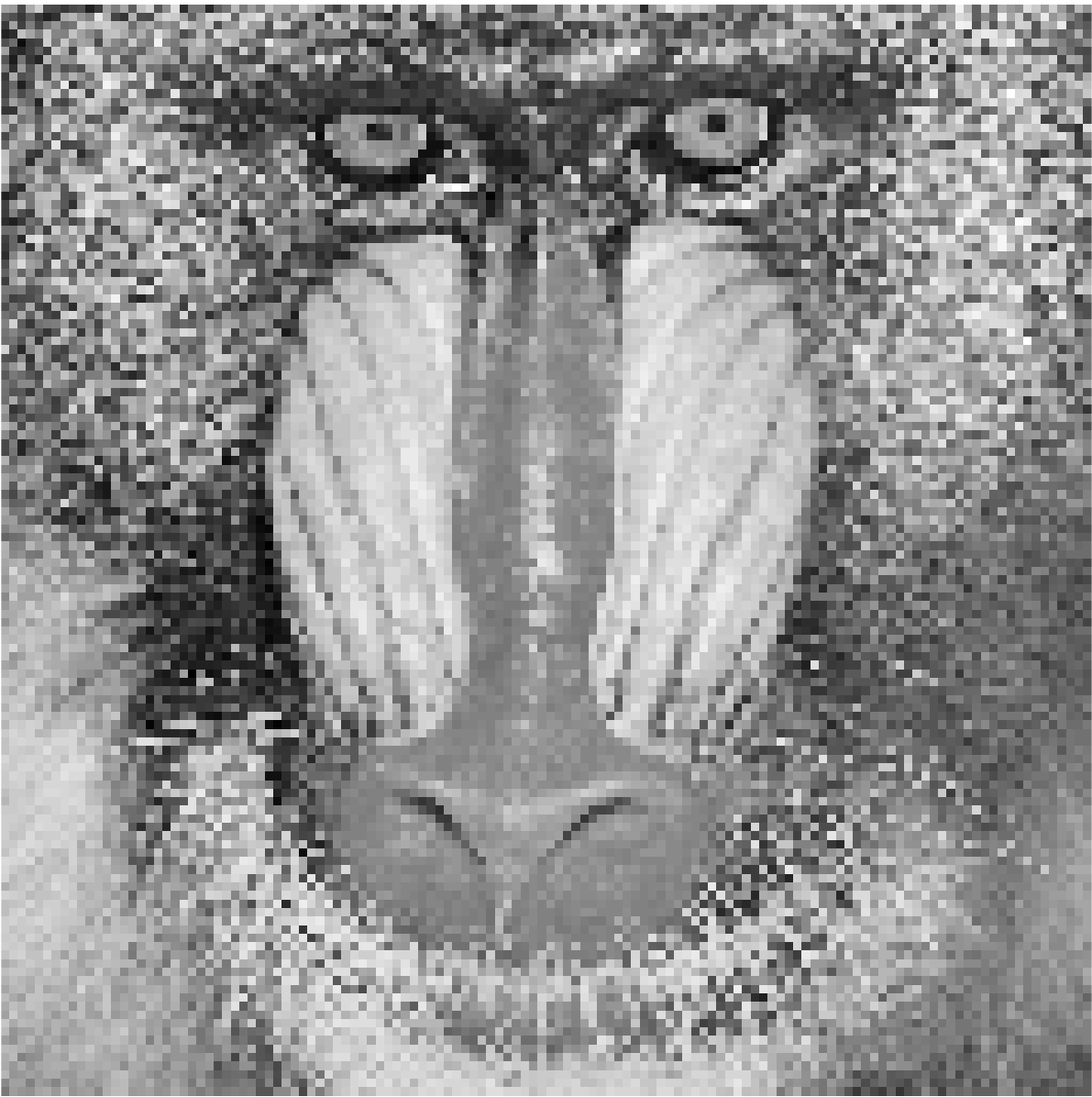}
\includegraphics[width=.25\textwidth]{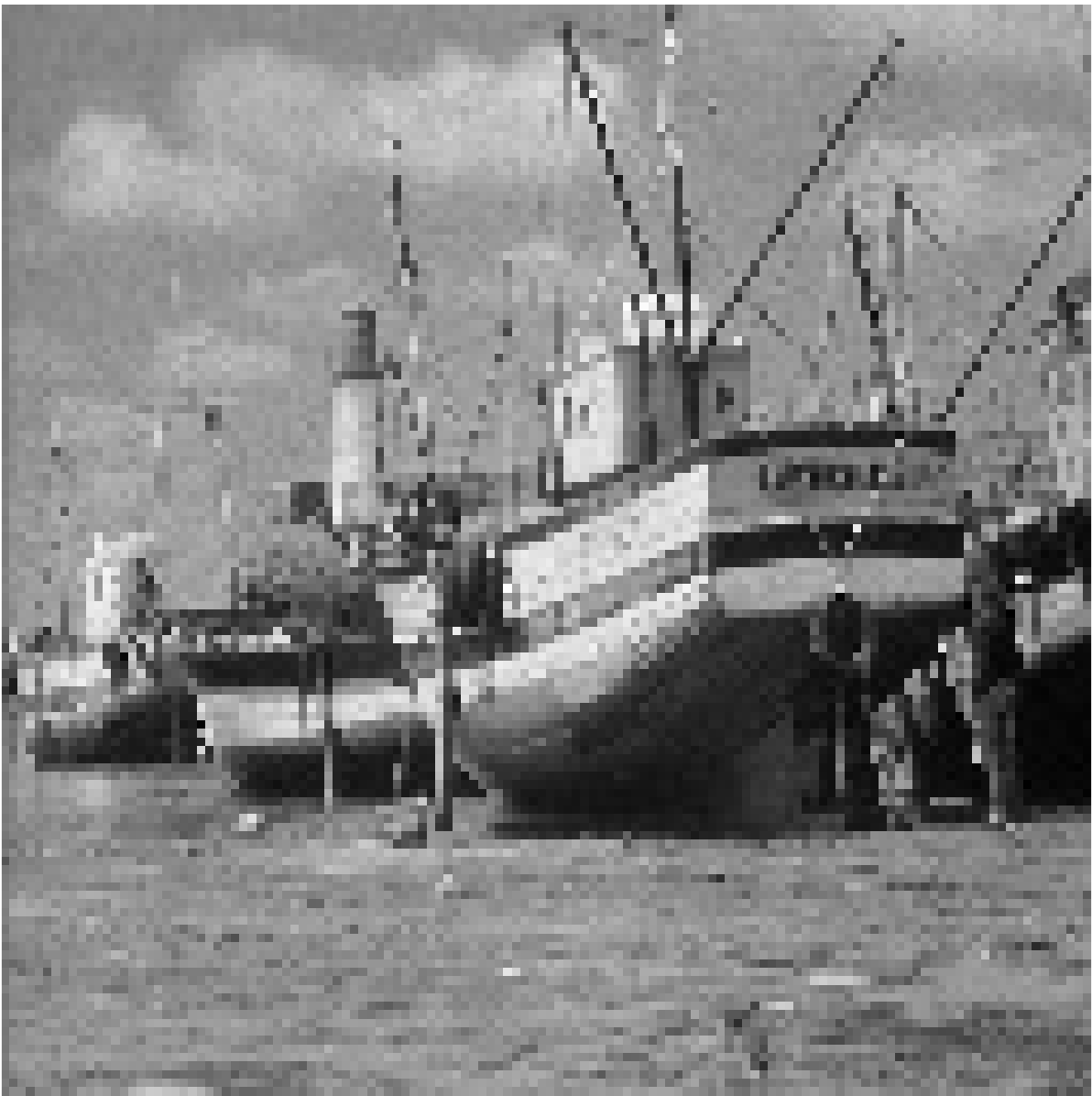}
\includegraphics[width=.25\textwidth]{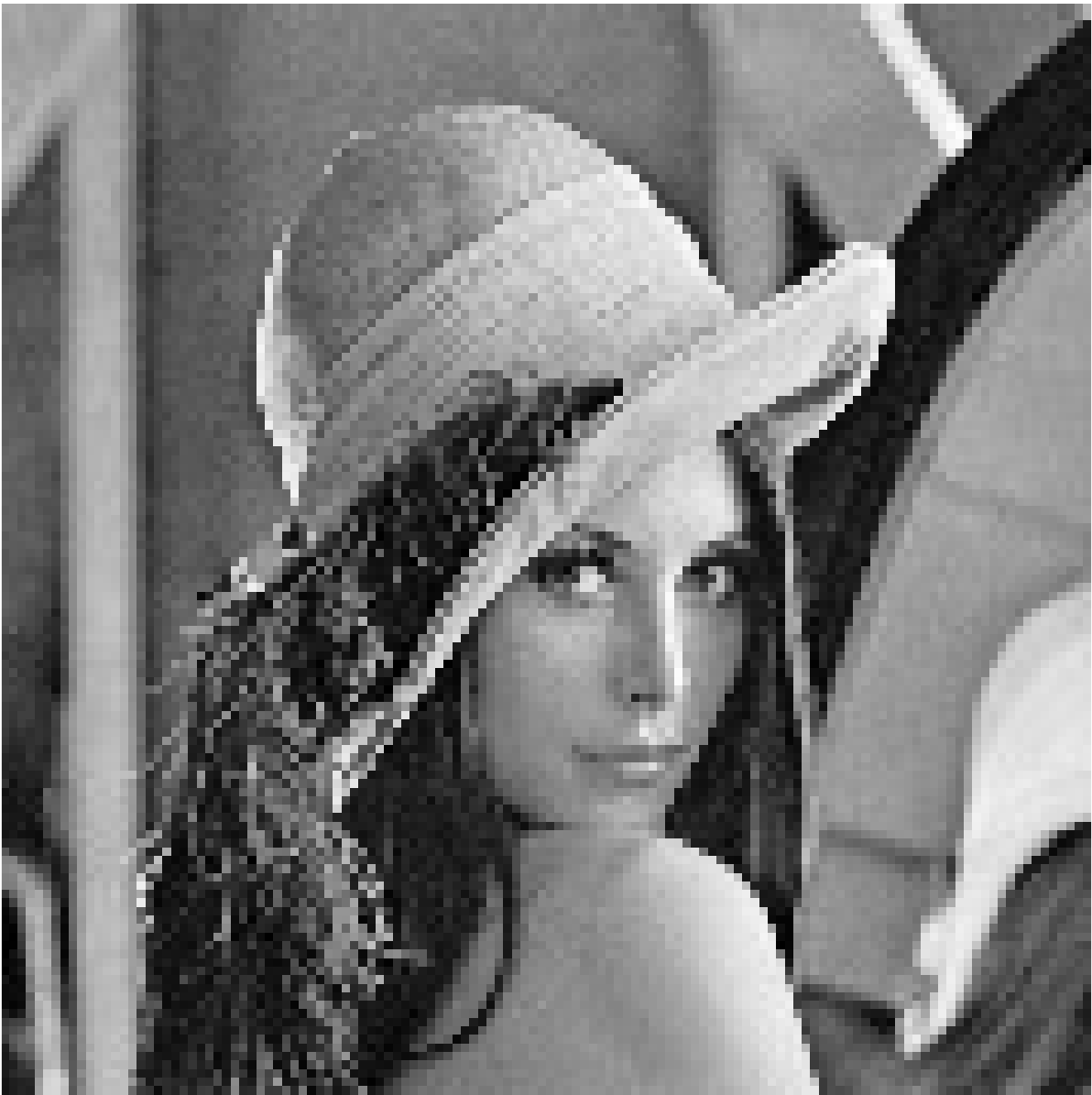}
\caption{Real images for numerical experiments, from left to right: mandrill, boat, Lena.}
\label{IP_images}
\end{figure}

Experiments for compression and denoising are considered. Motivated by extensive existing literature which utilizes Gabor-type expansions for texture analysis, we perform our experiments on nine texture samples of size $128\times 128$, cf., Figure \ref{texture_images}. We compare and contrast these results with experiments on three natural images, cf., Figure \ref{IP_images}. All images are taken from the USC-SIPI Image Database.  We hypothesized that our method would offer strong performance on textures, where Gabor methods are known to do well, and multiscale anisotropic methods are not known to be near optimal. For natural images, we hypothesized mixed results. On the one hand, these images do contain textures, but on the other, they also have smooth edges, where shearlets and curvelets and known to be theoretically near optimal. We note that our experimental images are downsampled to size $128\times 128$, for computational purposes.

\subsection{Experiments with Compression}

To evaluate a frame's efficiency in capturing directional information, we first perform compression experiments with global hard thresholding on the coefficients.  That is, we set a certain percentage of the smallest coefficients to 0, then reconstruct by computing the inverse frame operator.  The threshold is determined based on all coefficients in the images, independent of scale in the case of the redundant methods.  Note that this means many more coefficients are used in the redundant methods than in the non-redundant methods.  

\begin{figure}[H]    
\centering
\begin{minipage}[t]{0.15\textwidth}
\includegraphics[width=\linewidth]{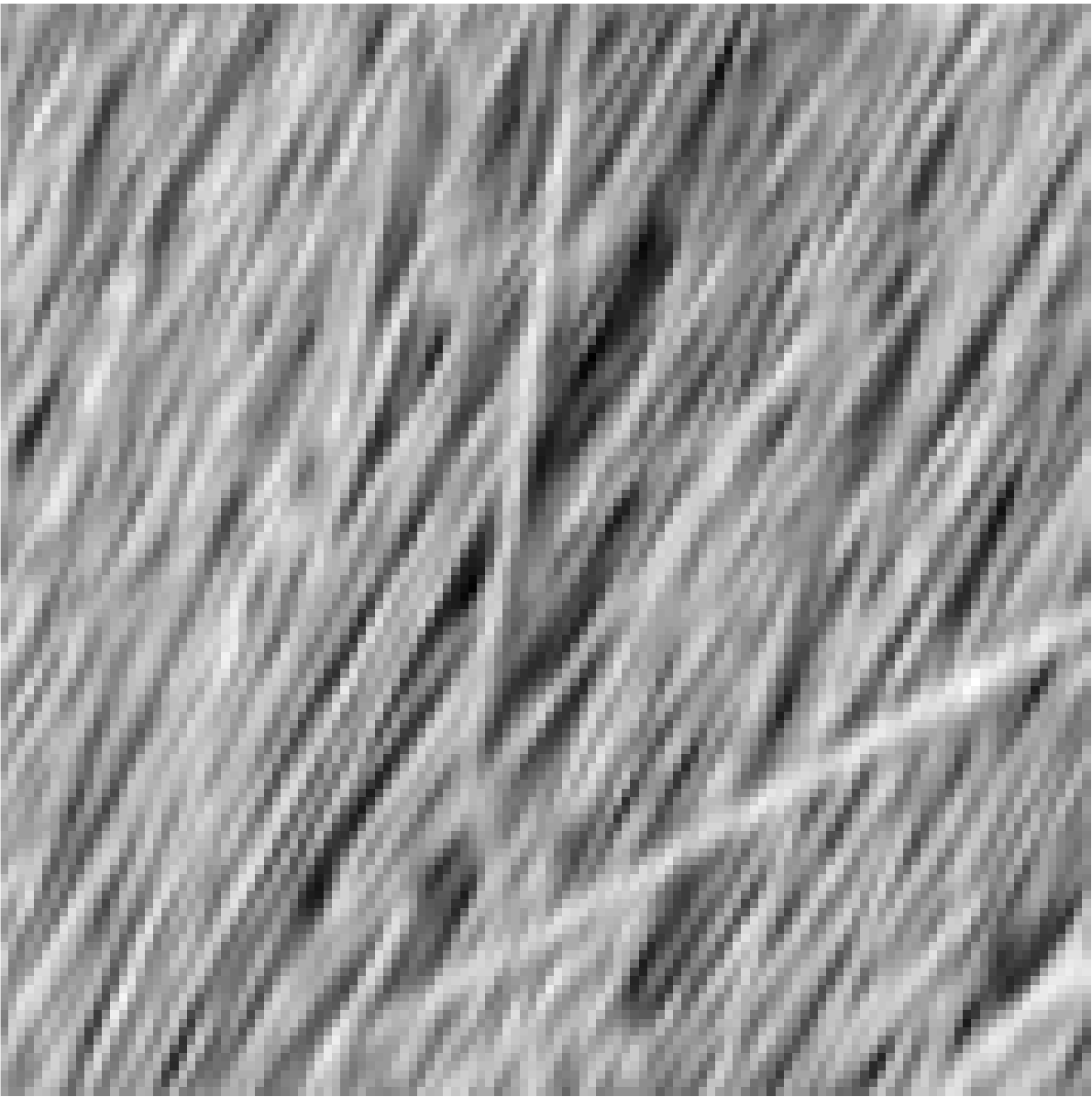}
\end{minipage}
\begin{minipage}[t]{0.15\textwidth}
\includegraphics[width=\linewidth]{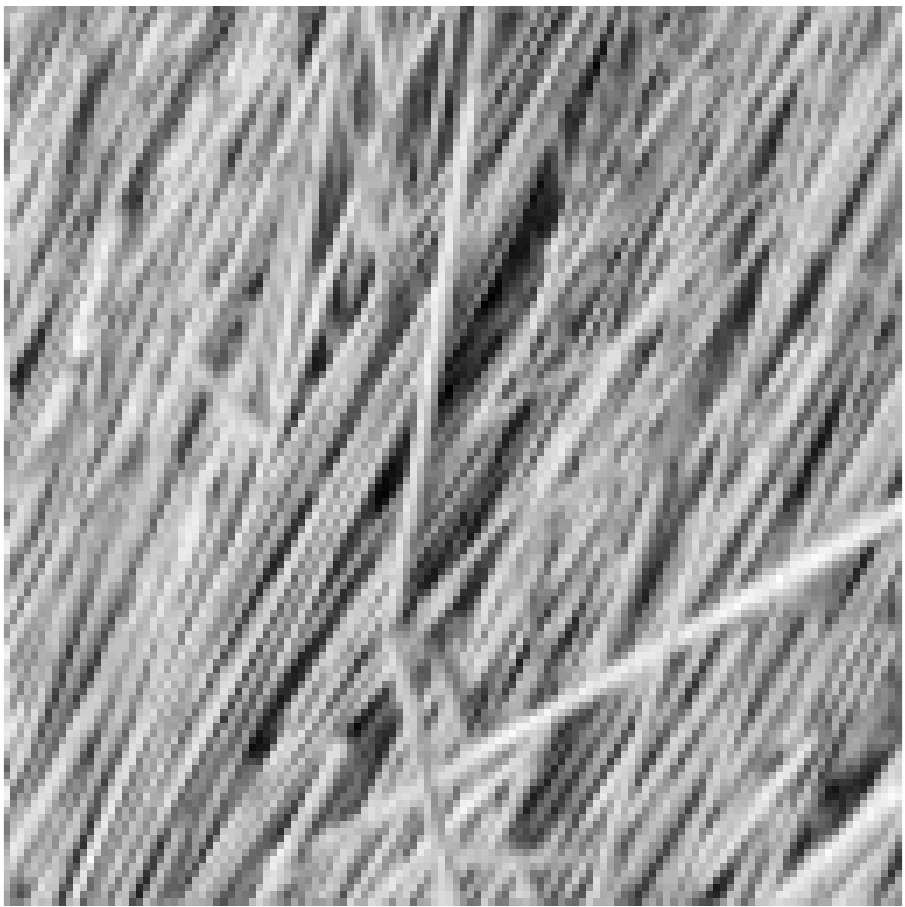}
\end{minipage}
\hspace{\fill}
\begin{minipage}[t]{0.15\textwidth}
\includegraphics[width=\linewidth]{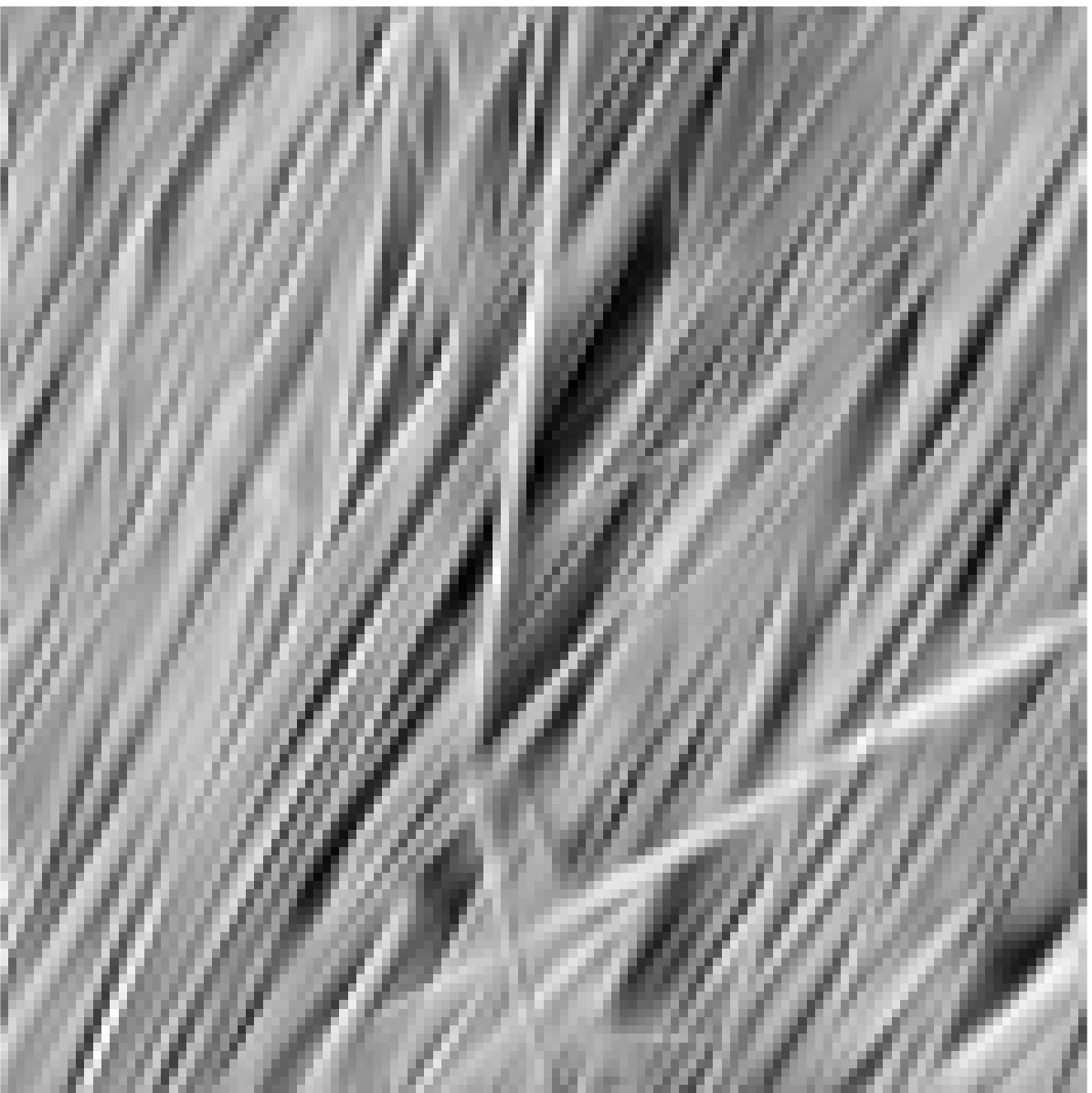}
\end{minipage}
\hspace{\fill}
\begin{minipage}[t]{0.15\textwidth}
\includegraphics[width=\linewidth]{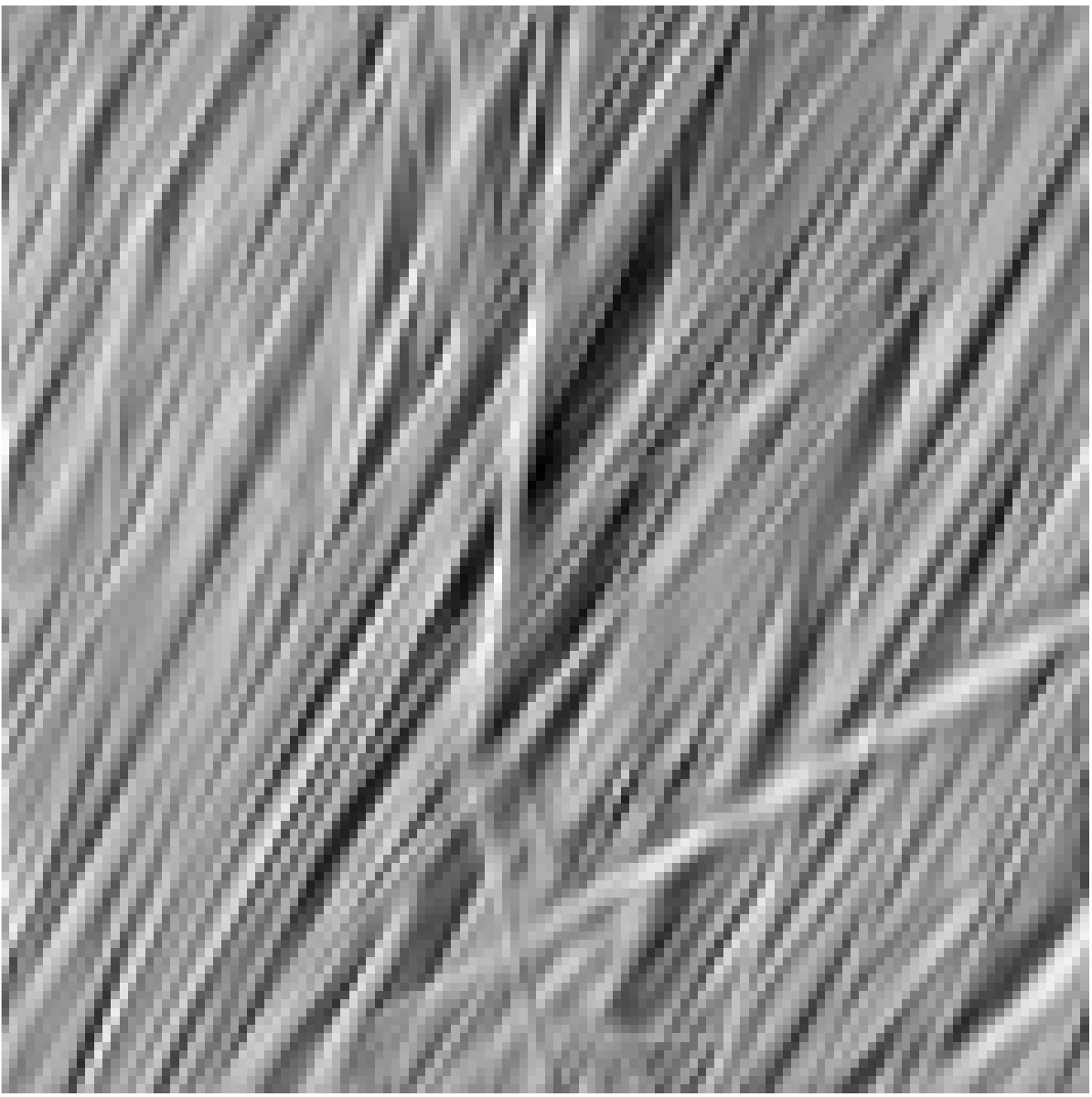}
\end{minipage}
\hspace{\fill}
\begin{minipage}[t]{0.15\textwidth}
\includegraphics[width=\linewidth]{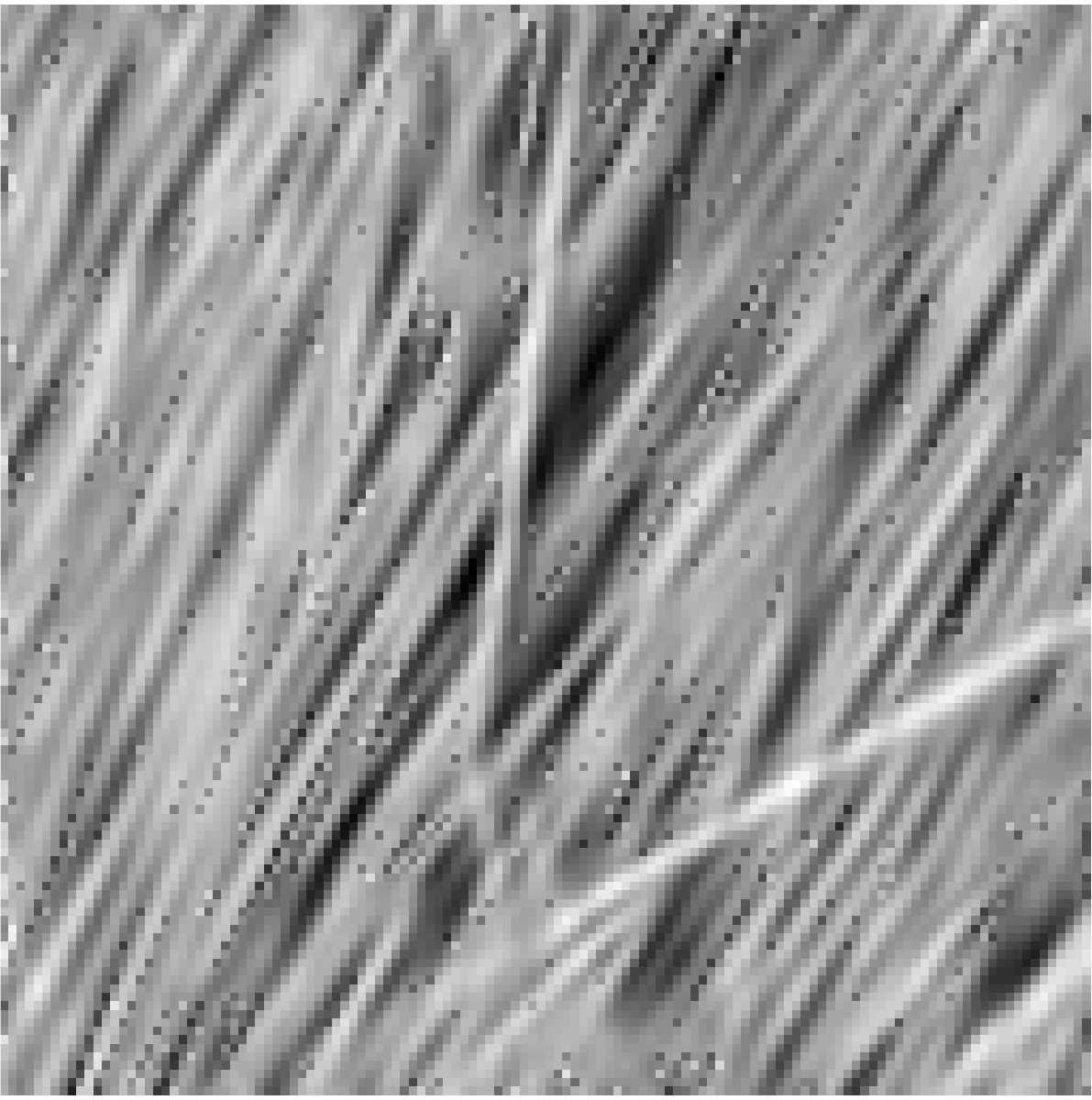}
\end{minipage}
\hspace{\fill}
\begin{minipage}[t]{0.15\textwidth}
\includegraphics[width=\linewidth]{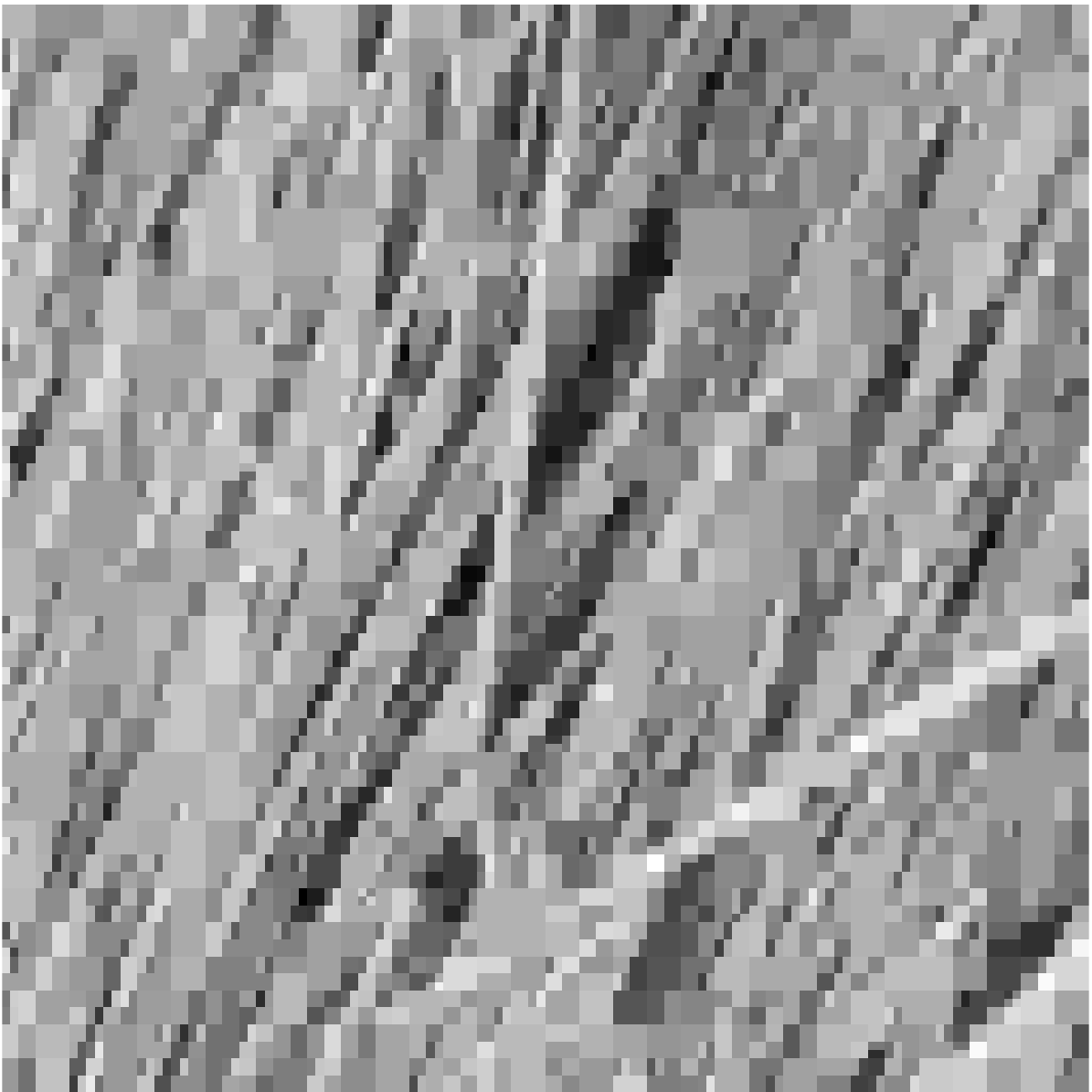}
\end{minipage}

\vspace*{0.5cm}

\begin{minipage}[t]{0.15\textwidth}
\includegraphics[width=\linewidth]{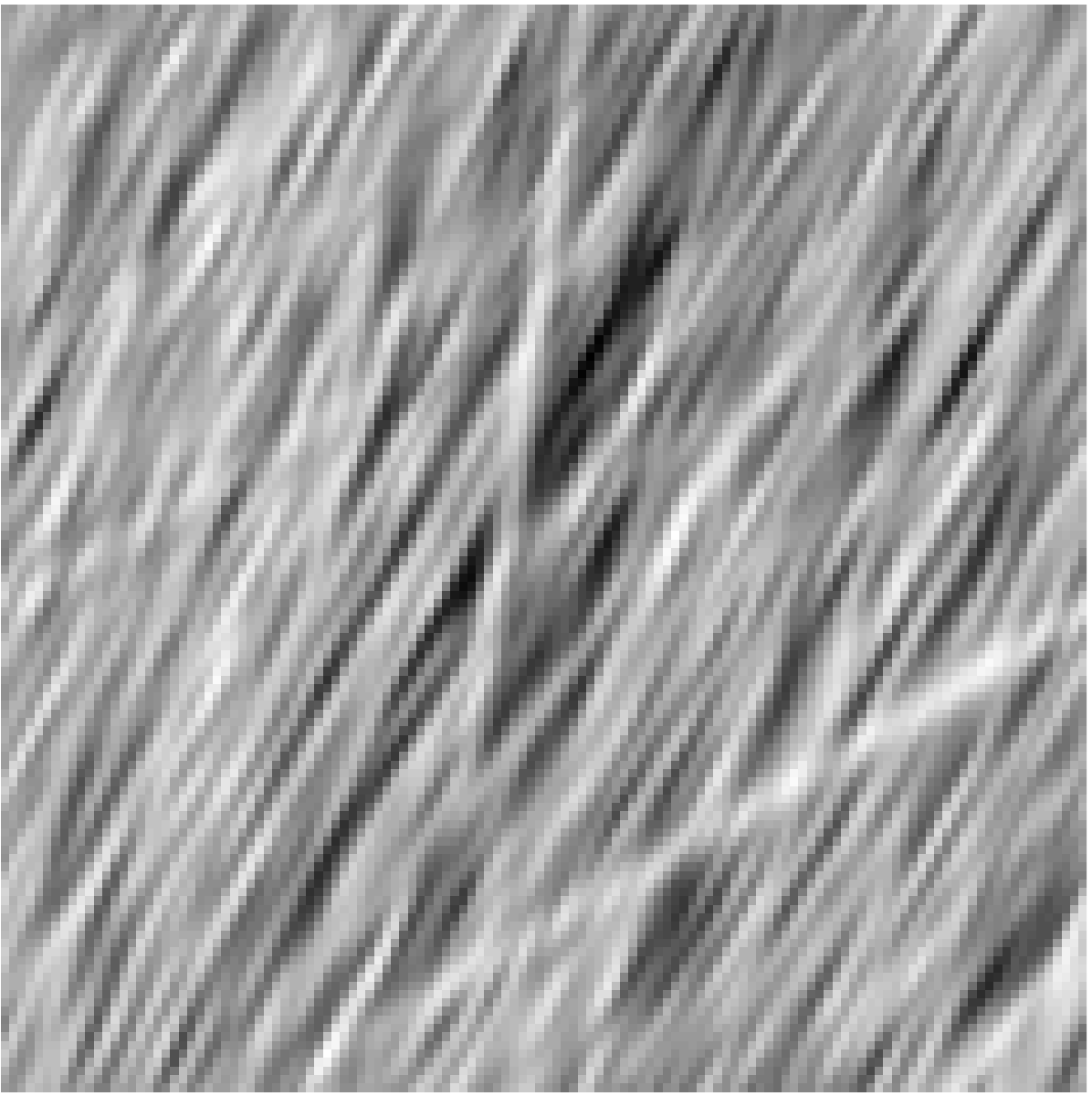}
\end{minipage}
\begin{minipage}[t]{0.15\textwidth}
\includegraphics[width=\linewidth]{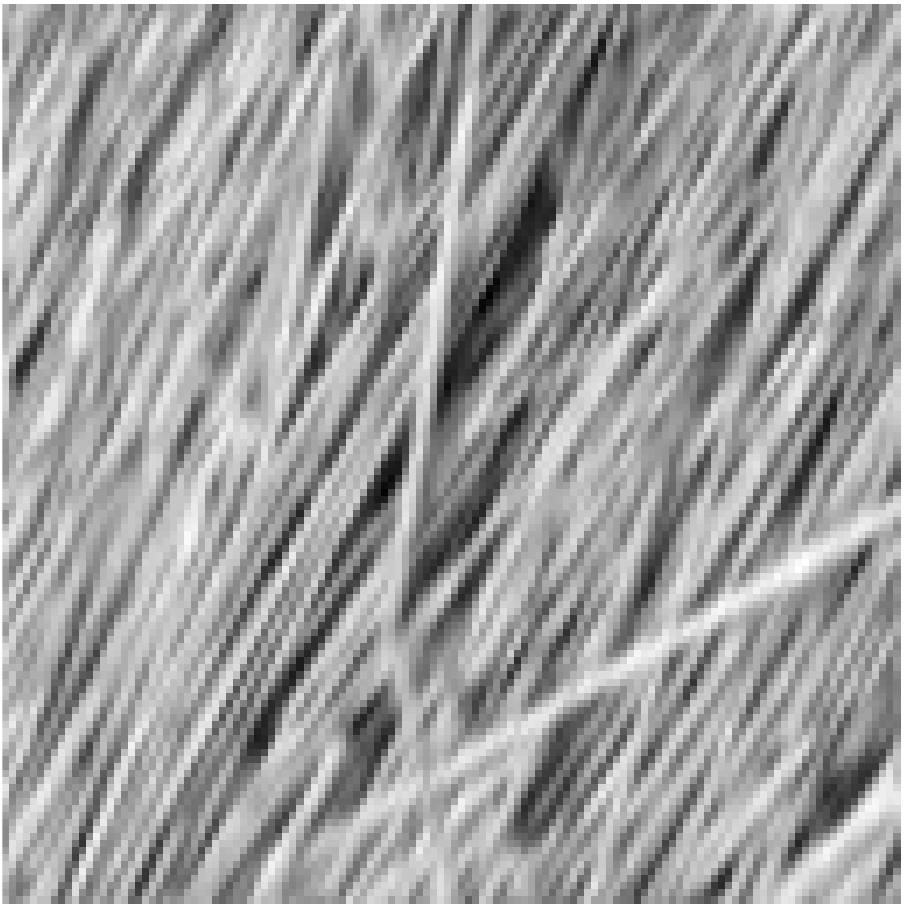}
\end{minipage}
\hspace{\fill}
\begin{minipage}[t]{0.15\textwidth}
\includegraphics[width=\linewidth]{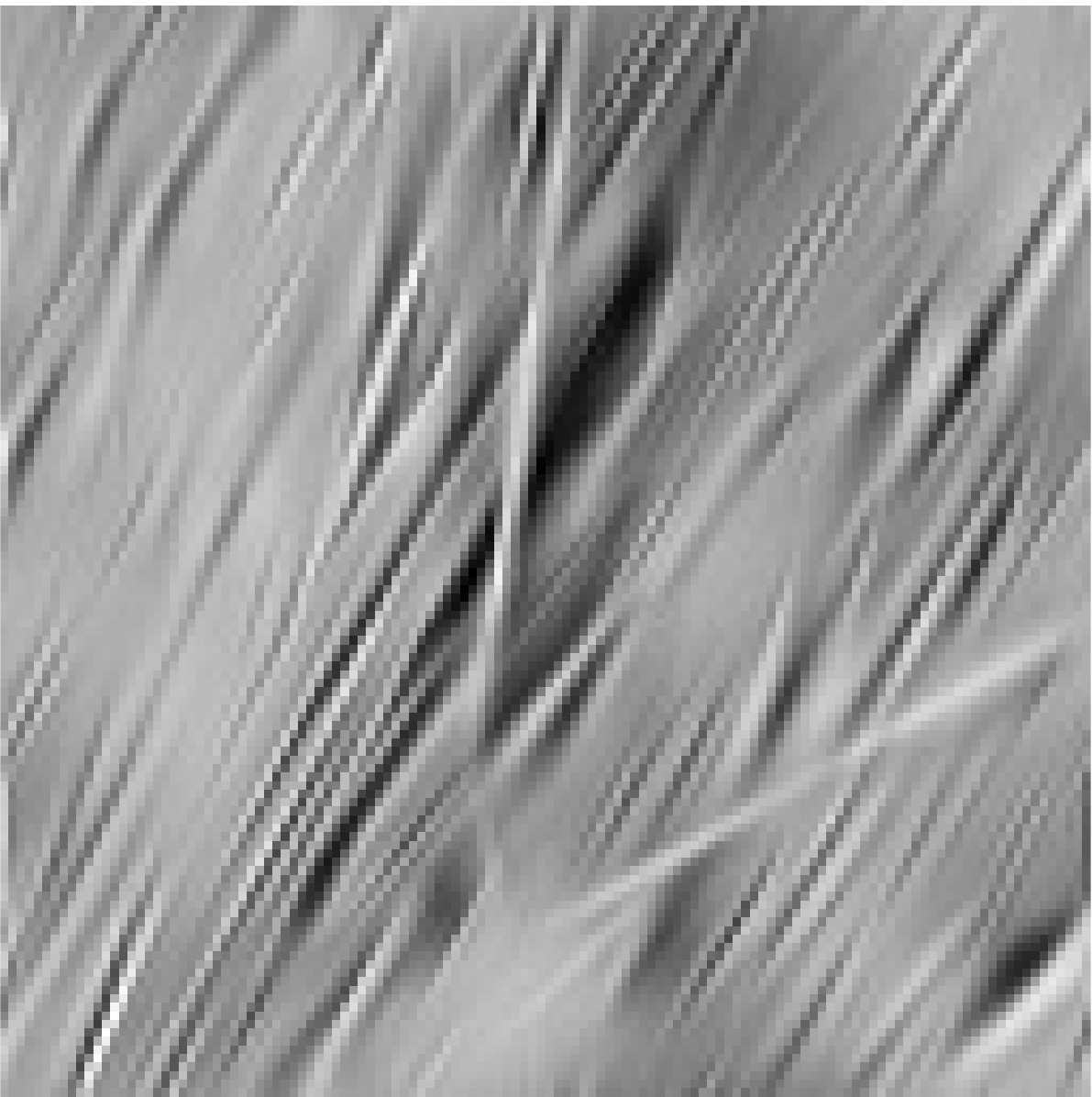}
\end{minipage}
\hspace{\fill}
\begin{minipage}[t]{0.15\textwidth}
\includegraphics[width=\linewidth]{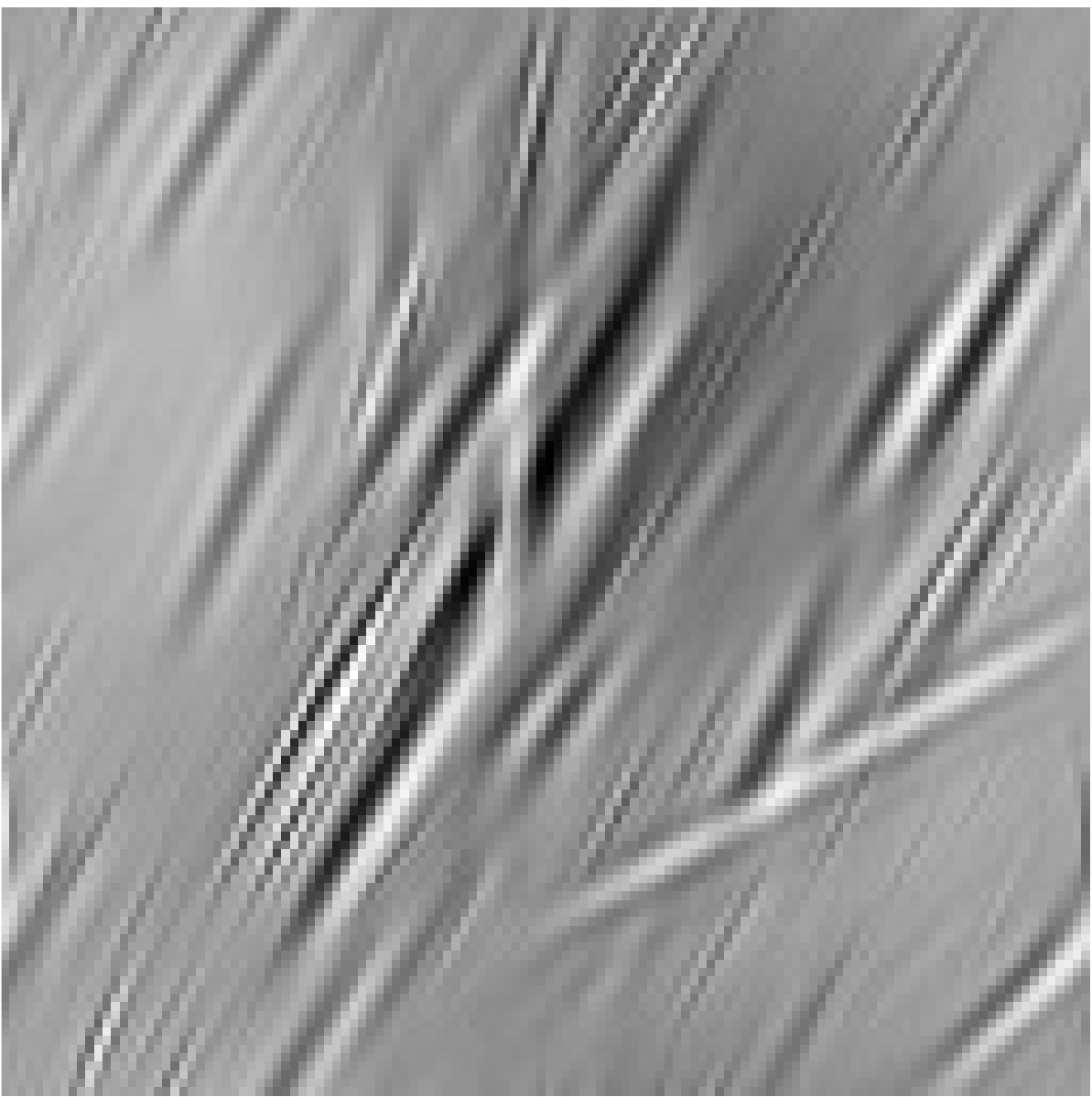}
\end{minipage}
\hspace{\fill}
\begin{minipage}[t]{0.15\textwidth}
\includegraphics[width=\linewidth]{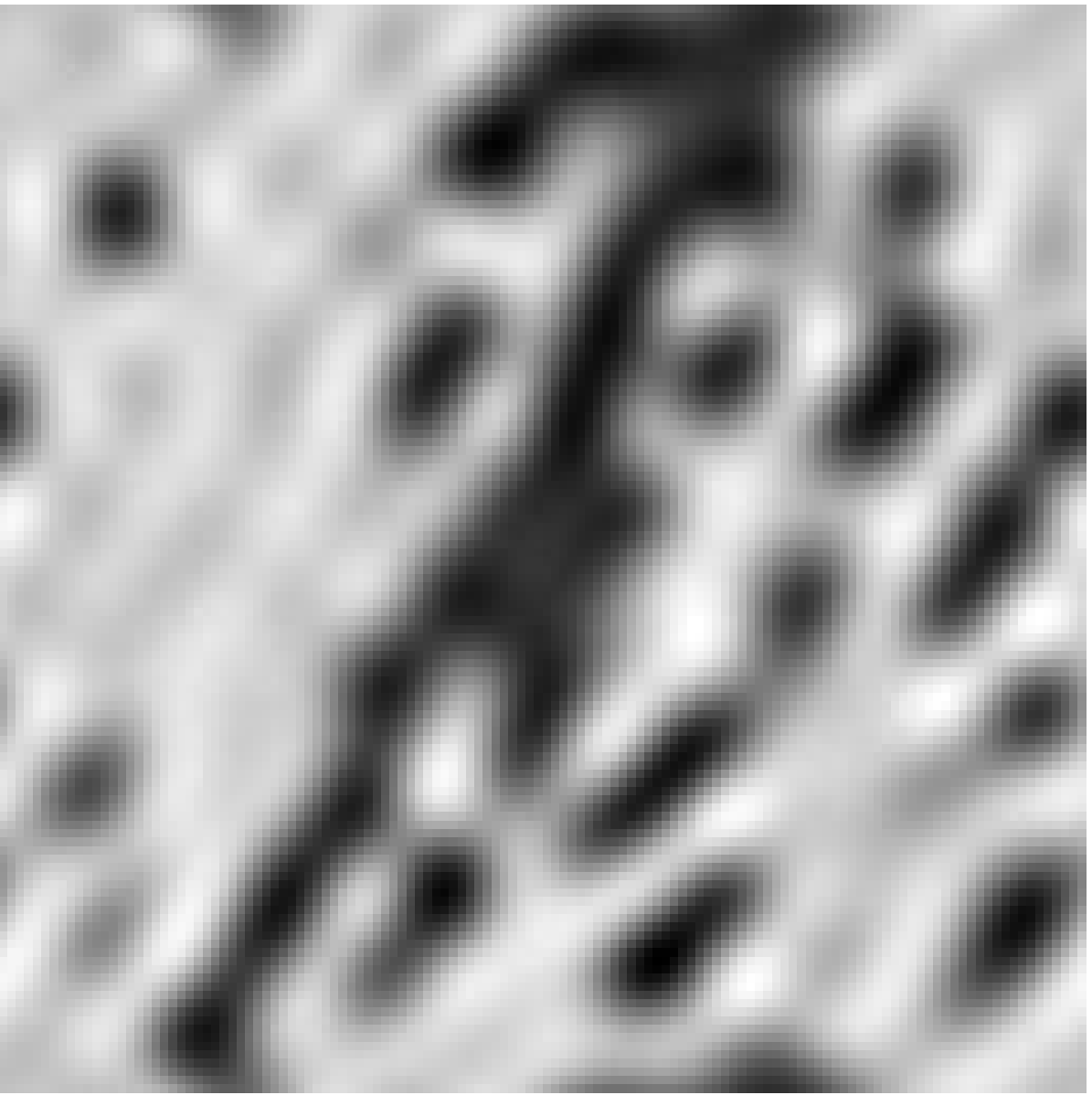}
\end{minipage}
\hspace{\fill}
\begin{minipage}[t]{0.15\textwidth}
\includegraphics[width=\linewidth]{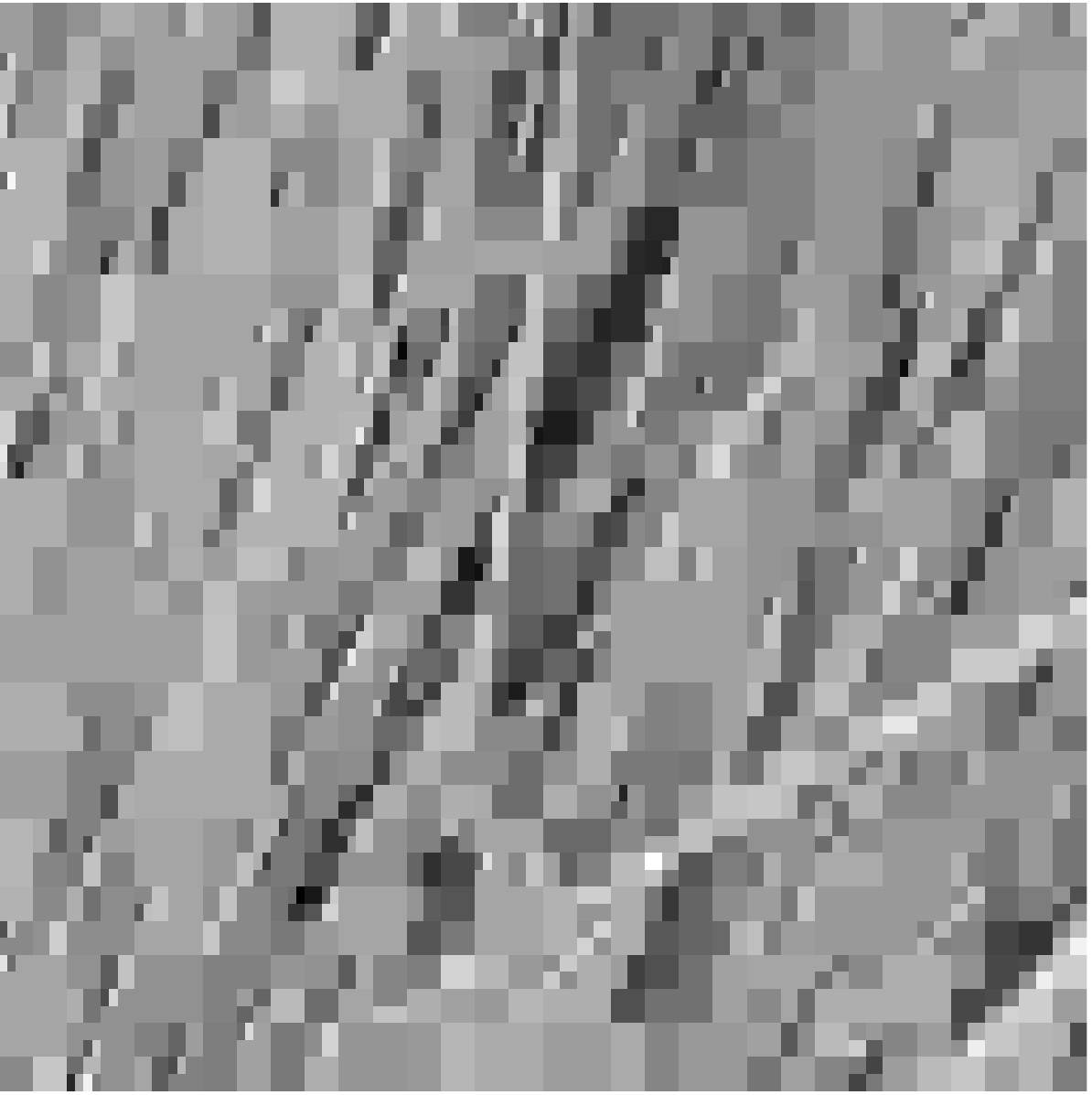}
\end{minipage}

\vspace*{0.5cm}

\begin{minipage}[t]{0.15\textwidth}
\includegraphics[width=\linewidth]{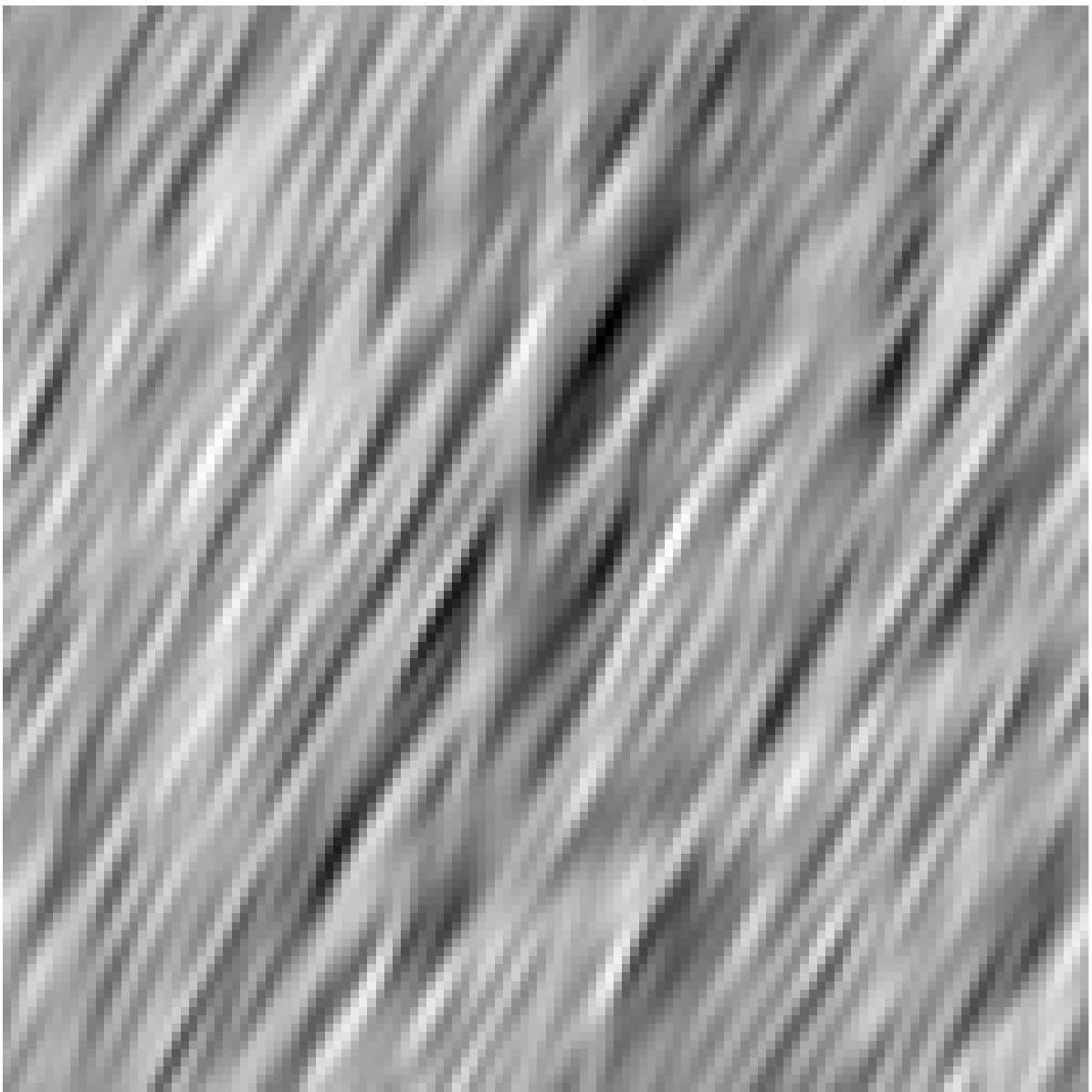}
\end{minipage}
\begin{minipage}[t]{0.15\textwidth}
\includegraphics[width=\linewidth]{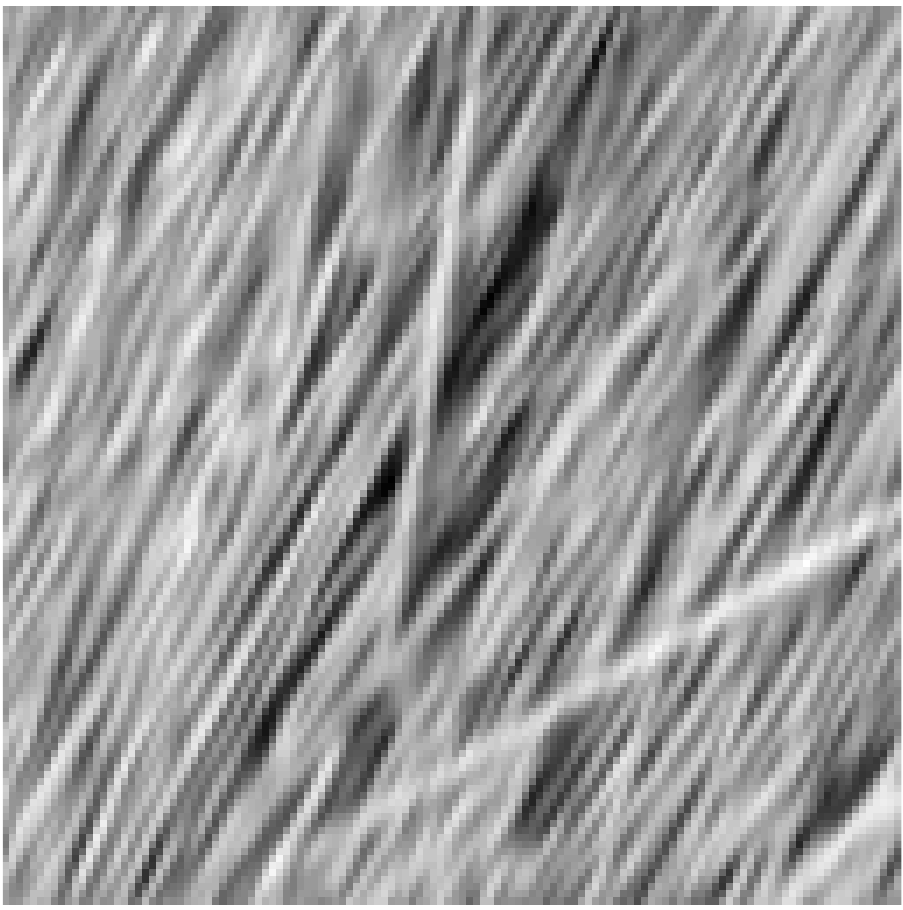}
\end{minipage}
\hspace{\fill}
\begin{minipage}[t]{0.15\textwidth}
\includegraphics[width=\linewidth]{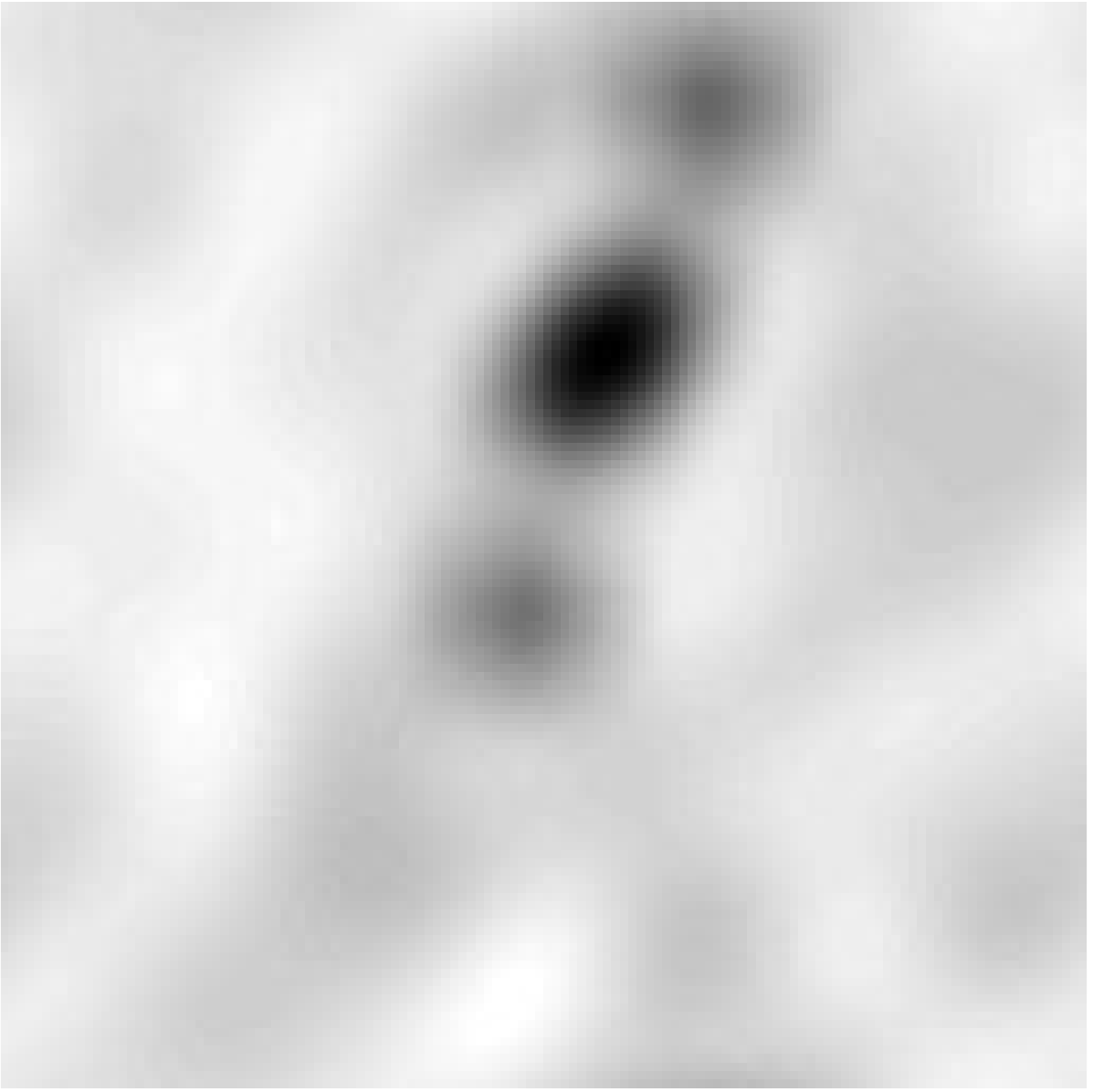}
\end{minipage}
\hspace{\fill}
\begin{minipage}[t]{0.15\textwidth}
\includegraphics[width=\linewidth]{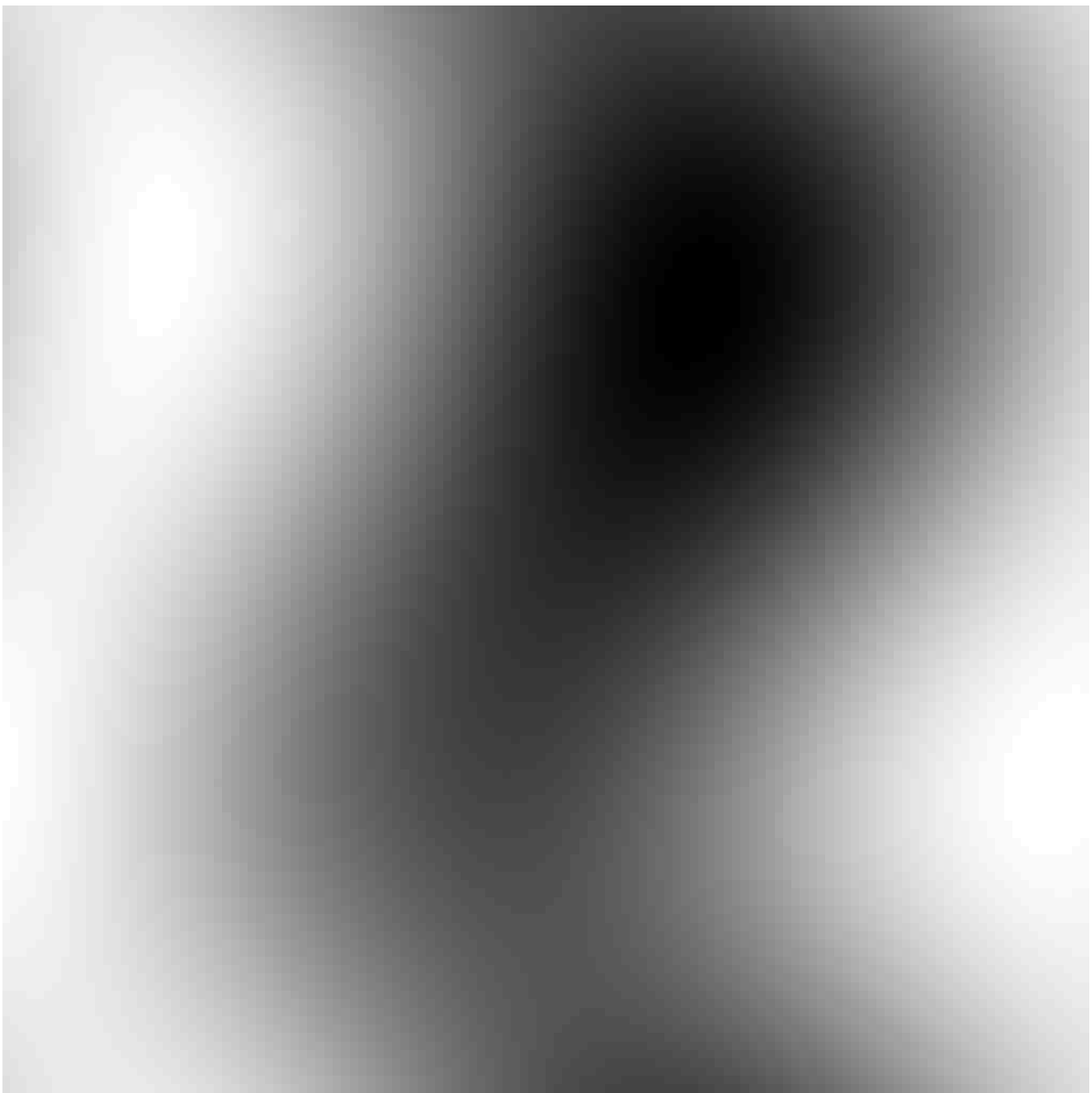}
\end{minipage}
\hspace{\fill}
\begin{minipage}[t]{0.15\textwidth}
\includegraphics[width=\linewidth]{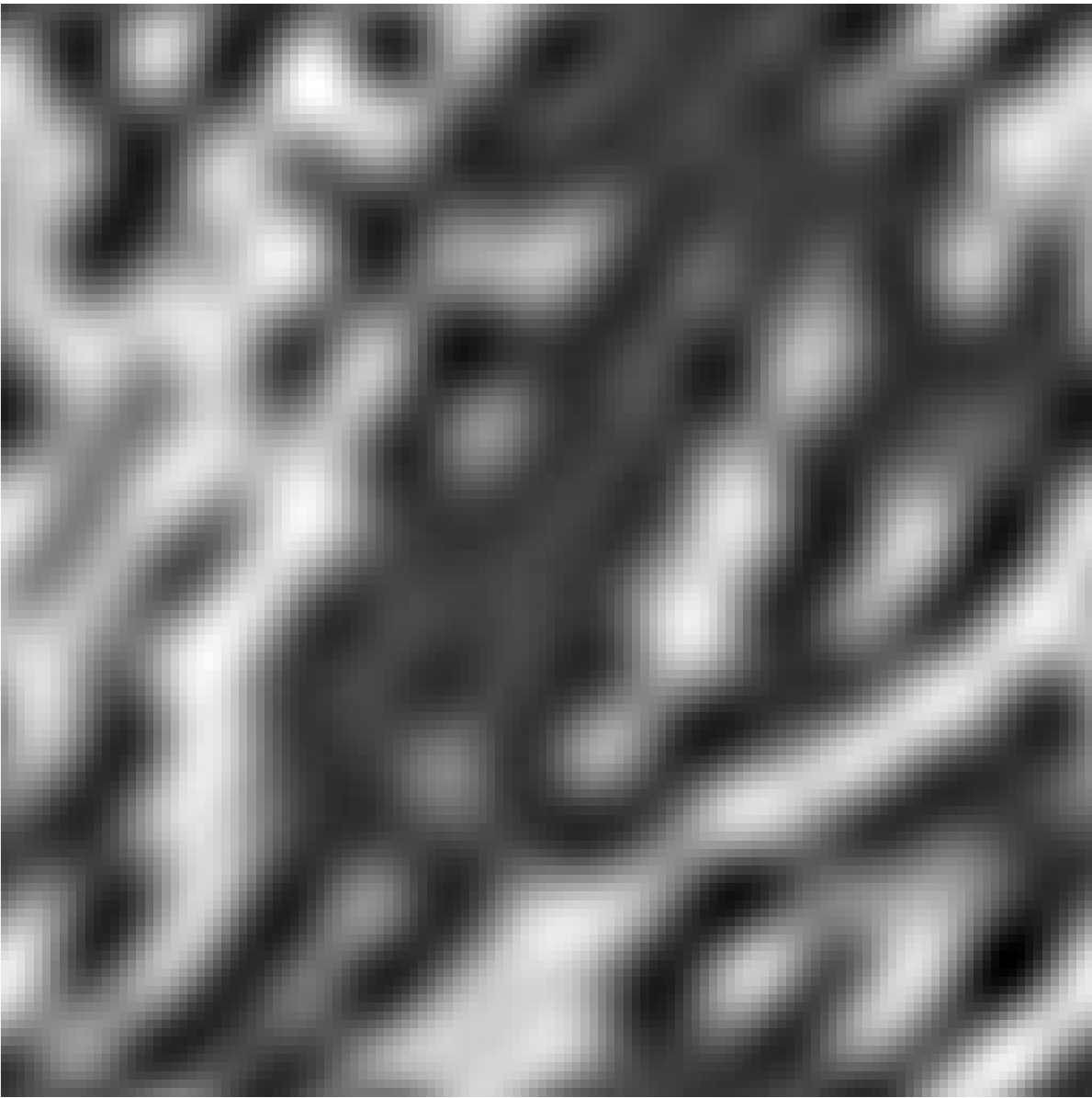}
\end{minipage}
\hspace{\fill}
\begin{minipage}[t]{0.15\textwidth}
\includegraphics[width=\linewidth]{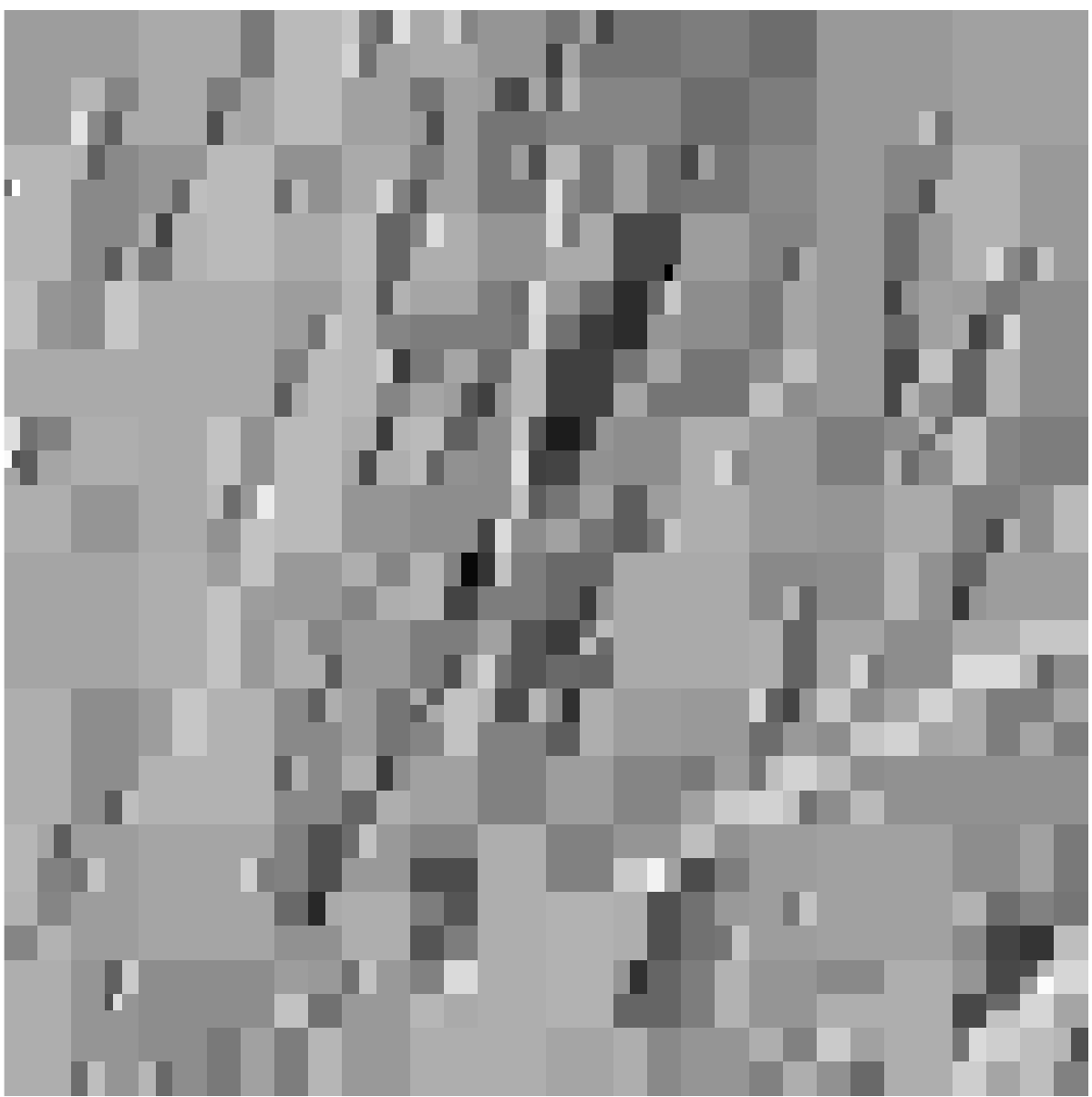}
\end{minipage}

\vspace*{0.5cm}

\begin{minipage}[t]{0.15\textwidth}
\includegraphics[width=\linewidth]{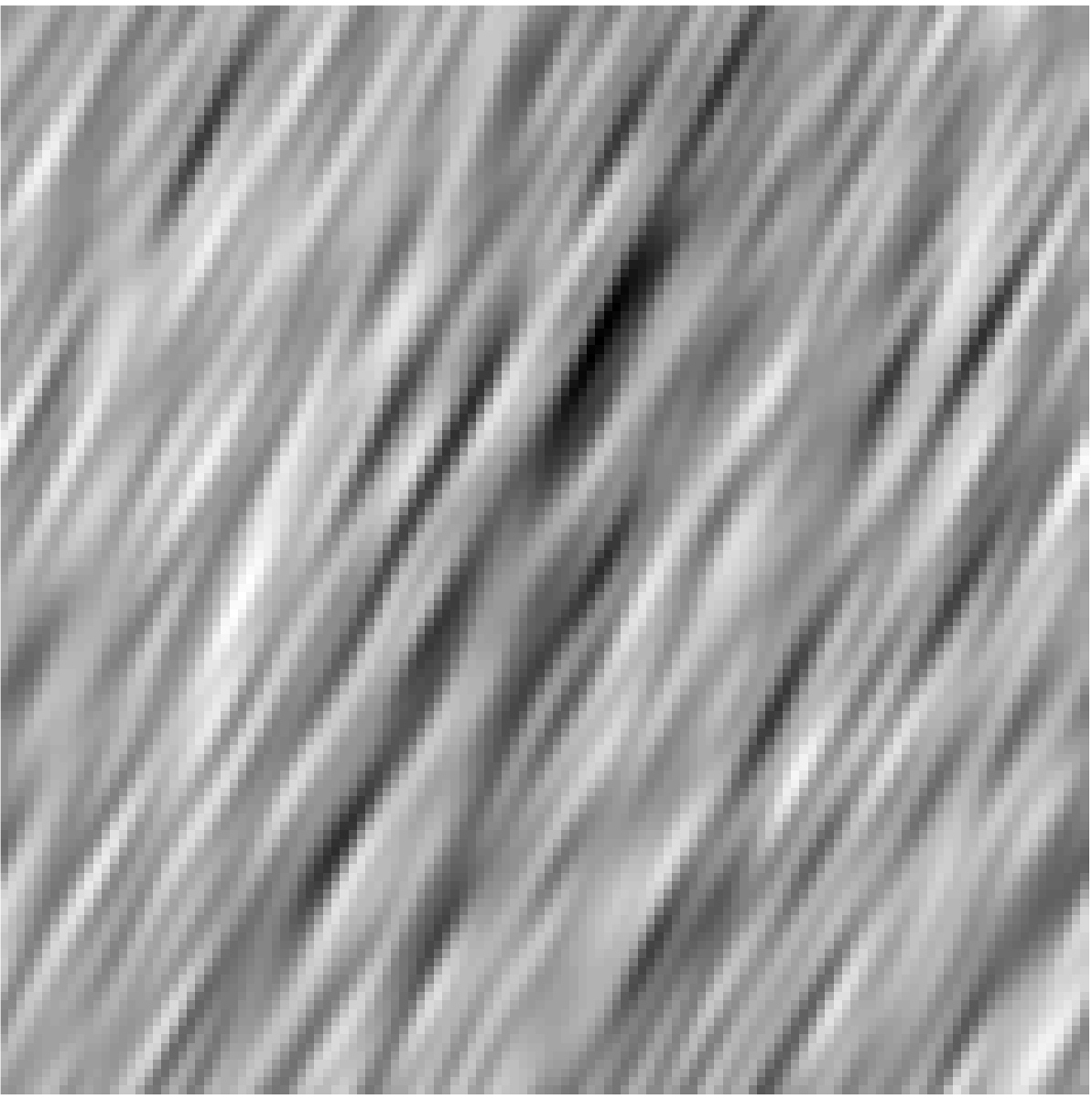}
\end{minipage}
\begin{minipage}[t]{0.15\textwidth}
\includegraphics[width=\linewidth]{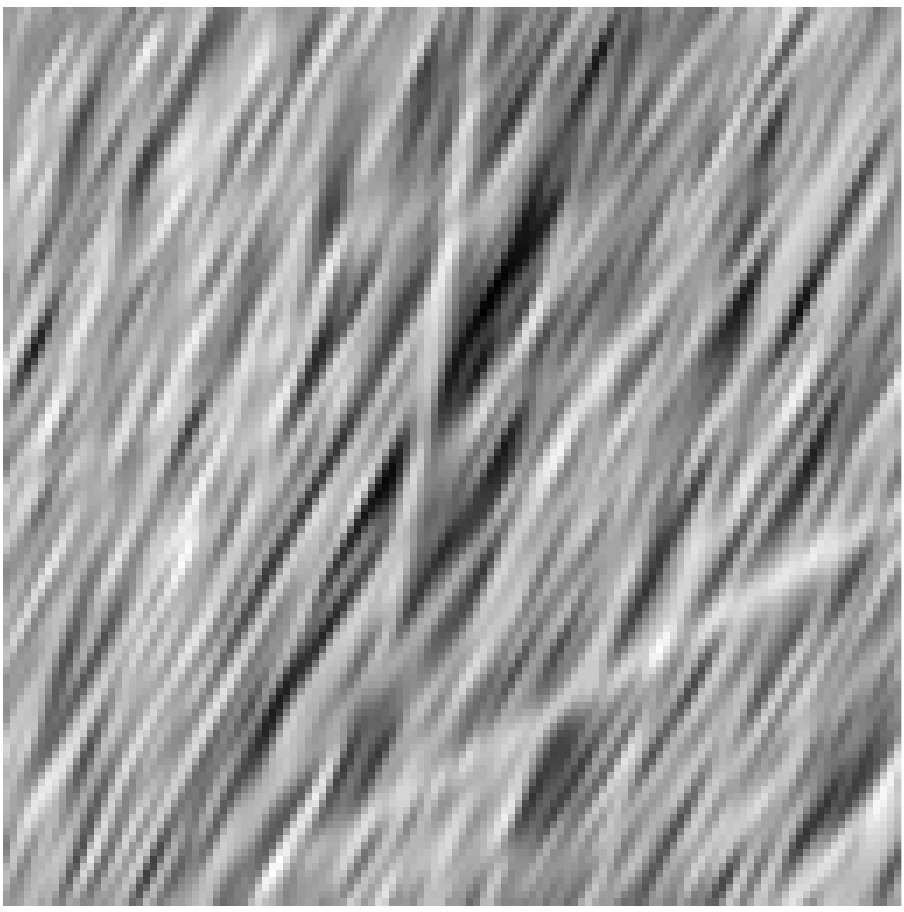}
\end{minipage}
\hspace{\fill}
\begin{minipage}[t]{0.15\textwidth}
\includegraphics[width=\linewidth]{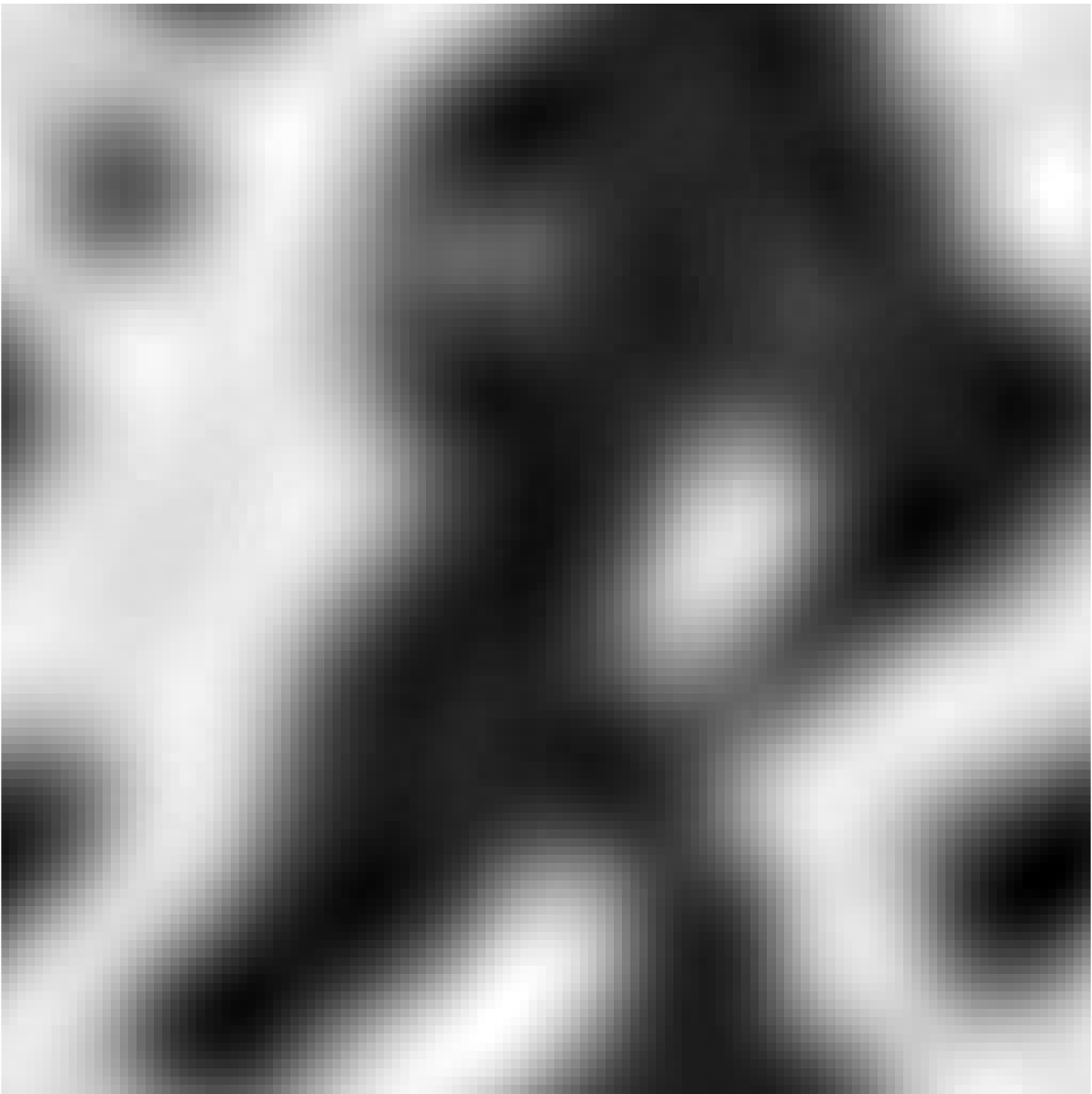}
\end{minipage}
\hspace{\fill}
\begin{minipage}[t]{0.15\textwidth}
\includegraphics[width=\linewidth]{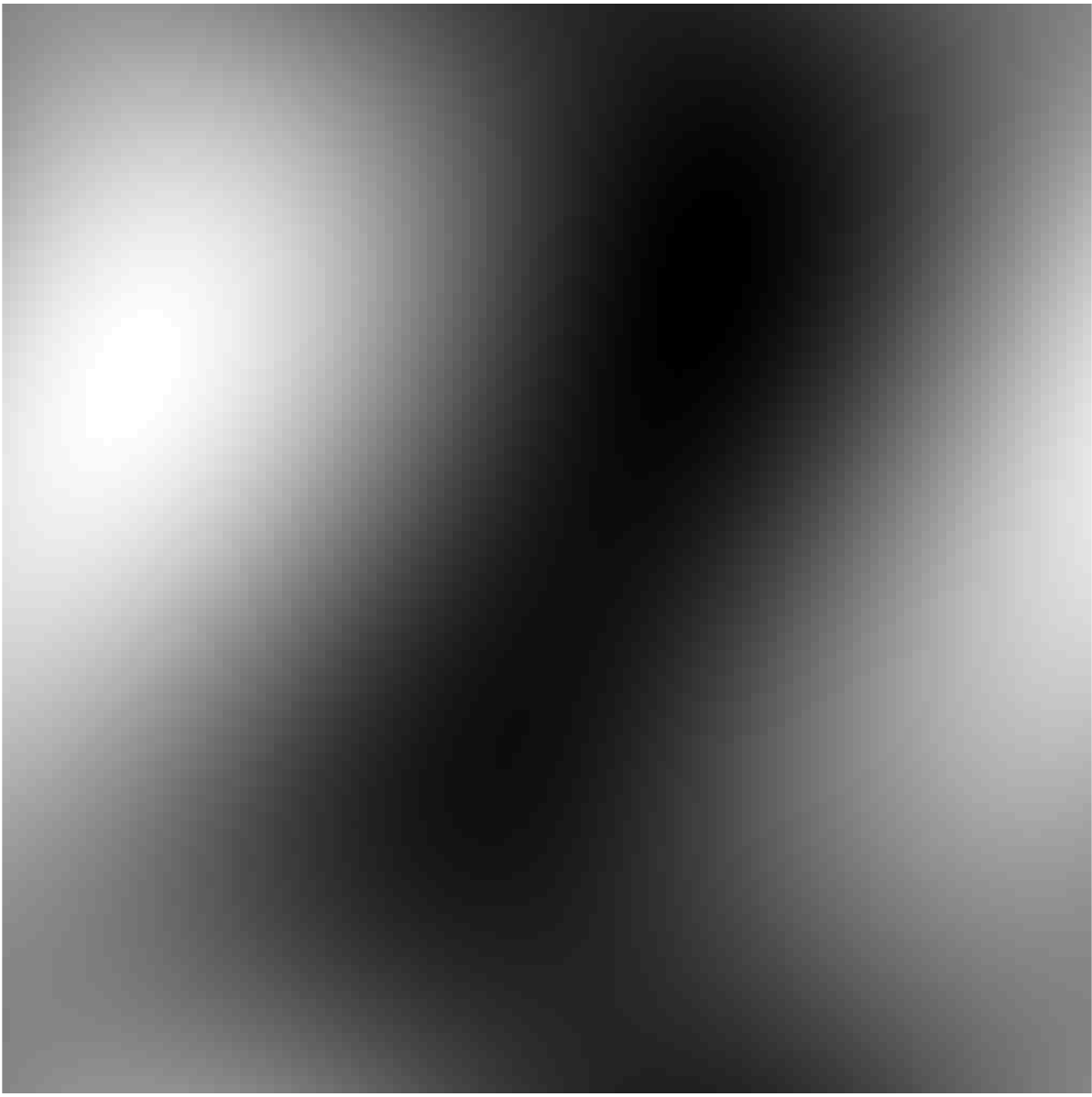}
\end{minipage}
\hspace{\fill}
\begin{minipage}[t]{0.15\textwidth}
\includegraphics[width=\linewidth]{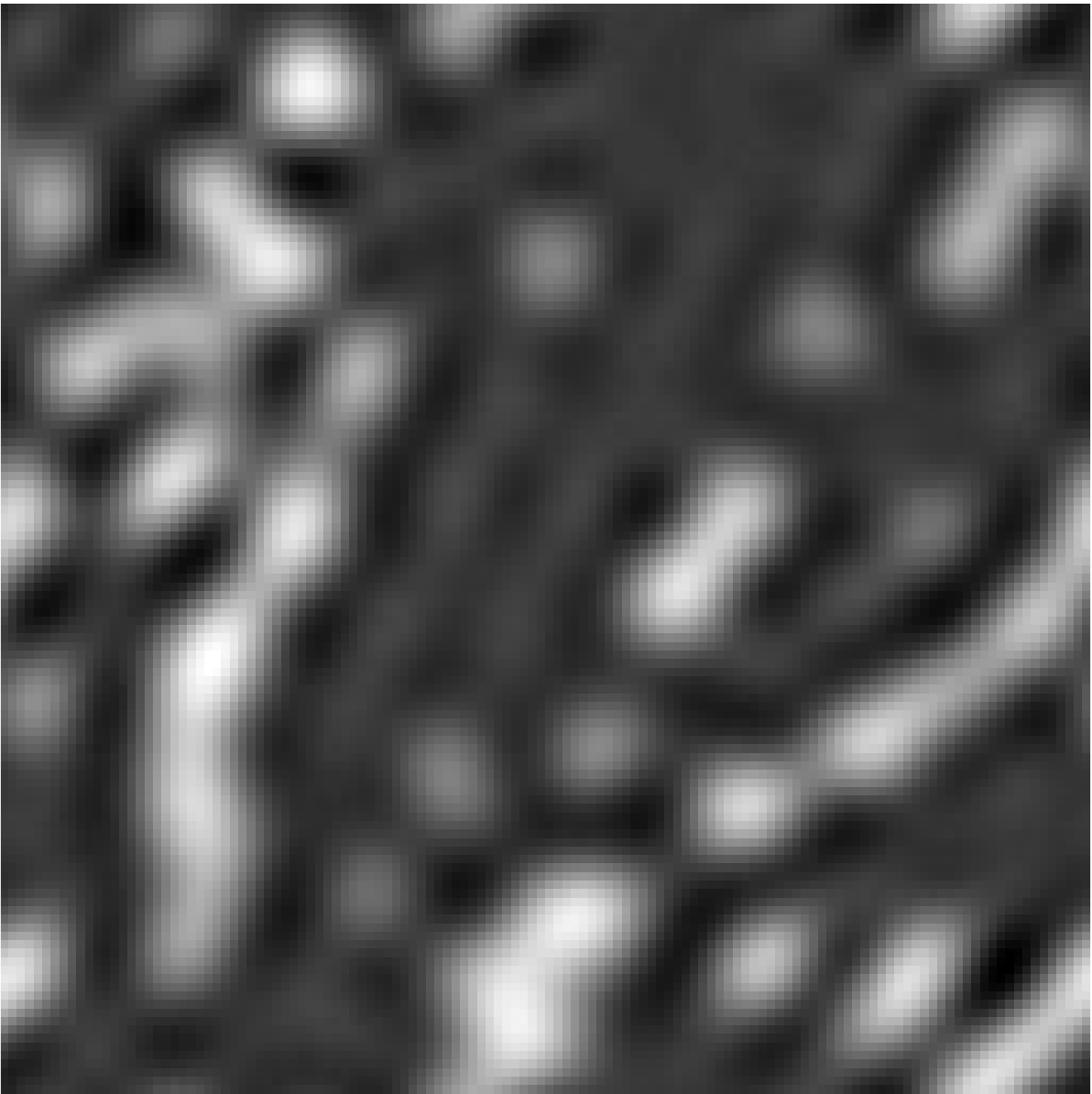}
\end{minipage}
\hspace{\fill}
\begin{minipage}[t]{0.15\textwidth}
\includegraphics[width=\linewidth]{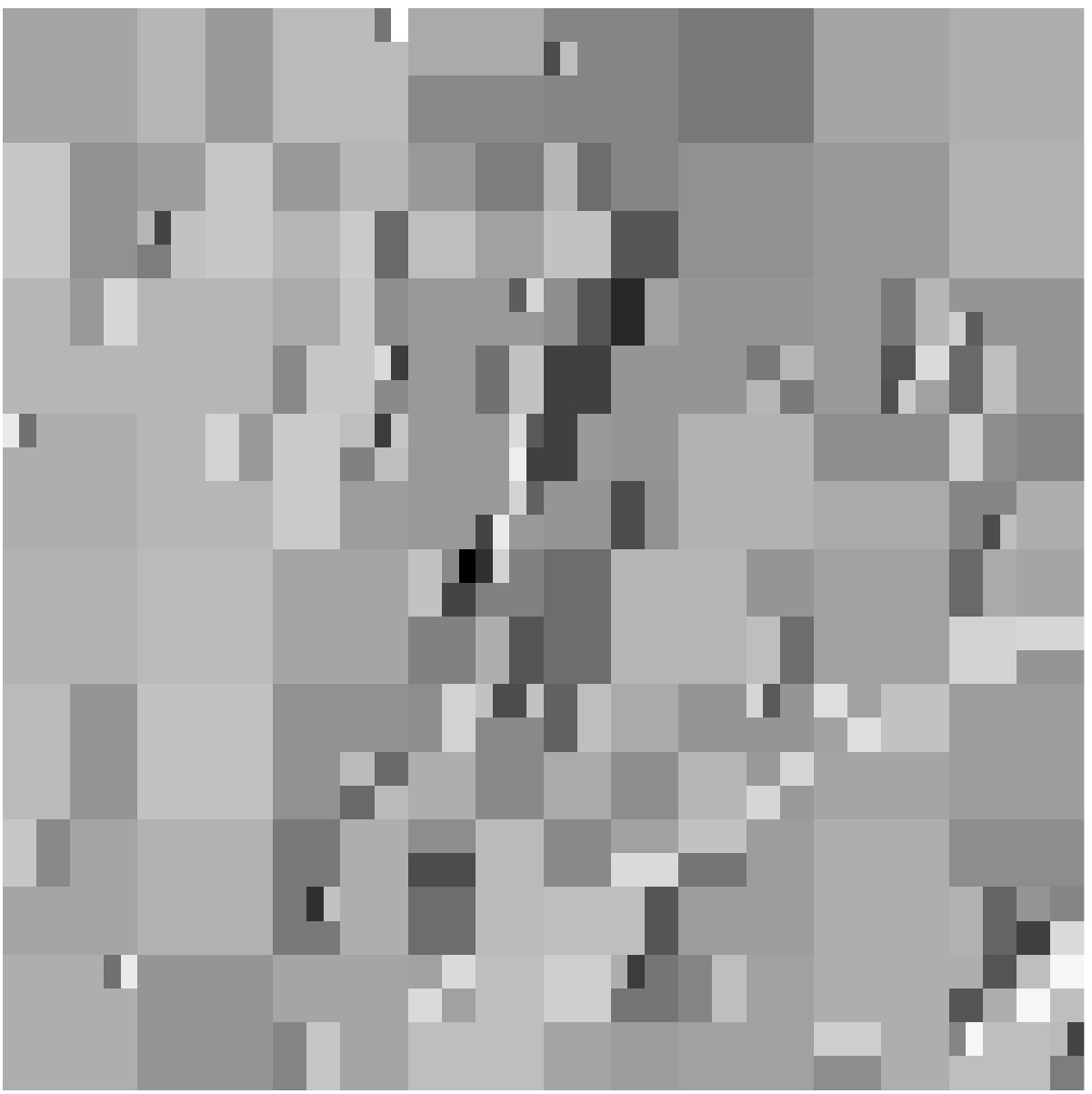}
\end{minipage}

\caption{Compressed images for straw texture.  The rows correspond to level of compression: from top to bottom, $90\%, 96\%, 98\%, 99\%$.  The columns correspond to method; from left to right: non-redundant directional Gabor frame, redundant directional Gabor frame, ShearLab, FFST, Curvelab, wavelet.}
\label{straw_approx}
\end{figure}

For our experiments, four levels of compression are considered, parameterized by compression ratios of $10,25,50,100$.  These correspond to $90\%$, $96\%$, $98\%$, $99\%$ thresholding, respectively.  Relative errors are displayed in Table \ref{tab:joint}.  We also display the compressed images produced by our experiments with straw texture in Figure \ref{straw_approx}.  Our experiments quantitatively indicate that the redundant directional Gabor systems performs optimally among the chosen methods.  Moreover, the performance is most impressive on the textures, where both the redundant and non-redundant directional Gabor methods consistently perform well, especially at high levels of compression. 

Qualitatively, our algorithm produces increasingly blurry images as compression increases.  The anisotropic methods tend to find the edges of the image quite well, but are unable to produce low error at high levels of compression.  This is not inconsistent with the the theoretical optimality of shearlets and curvelets, which are known to perform well theoretically for cartoon-like images in the asymptotic case.  For high levels of compression, we are far from the asymptotic case, so this theory no longer applies.  Moreover, the images under consideration here are strongly textural, and thus do no not fit in the cartoon-like regime.  Hence, although the multiscale anisotropic methods resolve the edges very well, they do not produce superior compression error in these experiments.  

\begin{table}[H]
\LARGE
\centering 
\resizebox{\linewidth}{!}{
\begin{tabular}{|P{3.3 cm} | P{3.8 cm}| P{2.6 cm} | P{2.6 cm} | P{2.6 cm} | P{2.6 cm} | P {2.6 cm} | P{2.6 cm}| P{2.6 cm} | P{2.6 cm} | P{3.5 cm} | P{2.6 cm} | P{2.6 cm} | P{2.6 cm} | P{2.6 cm} |}
\hline
Method & Compression Level & Straw & Rocks & Grass& Cracked Mud & Bricks & Bars & Fabric & Grate & Honeycomb  & Mandrill & Boat & Lena \\ \hline
\multirow{4} {*}{\parbox{2 cm}{\centering DGS, nonredundant}}& $90\%$ & .1099 & .0353 & .1451  & .1501 & .0851 & .0395 & .1699 & .1381 & .0361 & .1799 & .1288 & .1332 \\ 
			          & $96\%$ & .1574 &  .0561 & .2032 &  .1960 & .1021 & .0598 & .2140 & .1728 & .0541 & .2051 & .1570 & .1719 \\
			          & $98\%$ & .1921 &  .0778 & .2407 & .2225 & .1138 & .0869 & .2366 & .1964 & .0696 & .2174 & .1782 &  .2026  \\
			          & $99\%$ & .2212 &  .1042 & .2703 & .2437 & .1253 & .1299 & .2521 & .2196 & .0870 & .2279 & .2002 & .2373 \\ \hline
\multirow{4} {*}{\parbox{2 cm}{\centering DGS, redundant}} & $90\%$ & \textbf{.0453} & \textbf{.0158} & \textbf{.0597} & \textbf{.0633} & \textbf{.0412} & \textbf{.0170} & \textbf{.0692} & \textbf{.0640} & \textbf{.0148} & \textbf{.0906} & \textbf{.0619} & \textbf{.0642} \\ 
			          & $96\%$ & \textbf{.0859} &  \textbf{.0292} & \textbf{.1107} & \textbf{.1169} & \textbf{.0693} & \textbf{.0324} & \textbf{.1309} & \textbf{.1113} & \textbf{.0288} & \textbf{.1464} & \textbf{.1043} & \textbf{.1098} \\
			          & $98\%$ & \textbf{.1176} &  \textbf{.0410} & \textbf{.1525} & \textbf{.1538} & \textbf{.0865} & \textbf{.0457} & \textbf{.1695} & \textbf{.1416} & \textbf{.0408} & \textbf{.1766} & \textbf{.1310} &  \textbf{.1401}  \\
			          & $99\%$ & \textbf{.1509} &  \textbf{.0563} & \textbf{.1934} & \textbf{.1855} & \textbf{.1005} & \textbf{.0630} & \textbf{.2018} & \textbf{.1688} & \textbf{.0546} & \textbf{.1985} & \textbf{.1541} & \textbf{.1688} \\ \hline
\multirow{4}{*}{\parbox{2 cm}{\centering Shearlab}} & $90\%$ & .1191 & .0521 & .1808 & .1613 & .0912 & .0439 & .1651 & .1675 & .1003 & .1509 & .1079 & .1380 \\ 
				       & $96\%$ & .1946 & .0961 & .2542 & .2222 & .1335 & .1166 & .2364 & .2581 & .2582 & .2105 & .1682 & .2303 \\ 
				       & $98\%$ & .2784 & .1528 & .3136 & .2824 & .2482 & .2861 & .3079 & .3443 & .3974 & .2516 &  .2379 & .3202 \\
				       & $99\%$ & .6367& .5639 & .6131 &  .6250 & .6541 & .6699 & .6647 & .6675 & .6879 & .5590 & .5694 & .5264 \\ \hline
\multirow{4}{*}{\parbox{2 cm}{\centering FFST}} & $90\%$ & .1292 & .0820 & .1926 & .1742 & .0979 & .0467 & .1864 & .1925 & .1115 &  .1769 & .1300 & .1607 \\
				  & $96\%$ & .2393 & .1760 & .3025 & .2600 & .1849 & .1918 & .2778 & .3104 & .3589 & .2767 & .2392 & .3354 \\
				  & $98\%$ & .6037 & .5776 & .6061 & .5874 & .5786 & .6274 & .6085 & .6209 & .6655 & .5910 & .5230 & .5771 \\ 
				  & $99\%$ & .7860  & .7668 & .7884 & .7873 & .7961 & .8015 & .7991 & .8164 & .8413 & .7679 & .7250 & .7438 \\ \hline
\multirow{4}{*}{\parbox{2 cm}{\centering Curvelab}} & $90\%$ & .1297 & .0378 & .1613 & .1593 & .1232 & .1255 & .2013 & .2224 & .2334 & .1370 & .0979 & .1059 \\ 
					& $96\%$ & .1822 & .0613 & .2221 & .2112 & .1829 & .1874 & .2582 & .2906 & .3440 & .1898 & .1469 & .1687 \\
					& $98\%$ & .3062 & .1709 & .3136 & .3134 & .3186 & .3549 & .3584 & .3977 & .4562 & .2734 & .2243 & .2499 \\
			                 & $99\%$ & .6395 & .5963 & .6147 & .6380 & .6729 & .6959 & .6847 & .7028 & .7435 &  .5959 & .5926 & .5070  \\ \hline
\multirow{4}{*}{\parbox{2 cm}{\centering Wavelet}} & $90\%$ & .1569 & .0407 & .1607 & .1504 & .1084 & .1148 & .2019 & .2546 & .3033 &  .1528 & .0937 & .1008 \\ 
				     & $96\%$ & .2008 & .0642 & .2137 & .1991 & .1642 & .1804 & .2461 & .2967 & .3534 & .1926 & .1380 & .1609 \\
				     & $98\%$ & .2261 & .0859 & .2482 &  .2277 & .1957 & .2107 & .2677 & .3155 & .3728 & .2136 & .1682 & .2043 \\
				     & $99\%$ & .2469 & .1134 & .2771 &  .2521 & .2168 & .2331 & .2824 & .3275 & .3838 & .2320 & .1962 & .2473 \\ \hline
\end{tabular} 
}
\caption{\label{tab:joint}Compression errors.}
\end{table}

\subsection{Experiments with Denoising}

To perform denoising, all images were converted to grayscale and corrupted with Gaussian noise of mean $\mu=0$ and variance $\sigma^{2}=.1$.  The goal is to remove as much of the noise as possible, without adding artifacts.  A standard metric for denoising is the peak signal to noise ratio (PSNR):

\begin{align*}PSNR(\tilde{I},I)=-10\log_{10}(\|\tilde{I}-I\|_{2}),\end{align*}where $\tilde{I}$ is the denoised version of the noisy image, and $I$ is the original, noiseless image.  Higher levels of PSNR indicate superior denoising.  However, it is also important to consider visual evaluation for denoising, as PSNR can sometimes be misleading \cite{Morel2005}.

To evaluate the proposed methods for denoising, we considered universal hard and soft thresholding in the transform domain at a range of thresholding levels.  We considered proportional thresholding values running from .75 to .995, incrementing by steps of size .005.  For example, a thresholding value of .95 means denoising was performed by keeping only the top $(1-.95) \times 100 \% = 5 \%$ coefficients.  Lower proportional thresholds were also considered, but the results for these values were universally poorer.  Each method was then evaluated at its best thresholding level, as determined by this sweep over the proportional thresholding parameter and type of thresholding.  We then averaged results over 10 instances, in order to remove the stochastic effects of random noise generation.  The results for universal hard thresholding appear in Table \ref{tab:PSNR}; results were similar but generally poorer for universal soft thresholding.

\begin{table}[H]
\Large
\begin{center}
\resizebox{\textwidth}{!}{
\begin{tabular}{| P{3.3 cm} | P{3.3 cm} | P{3.3 cm} | P{3.3 cm} | P{3.3 cm} | P{3.3 cm} | P{3.3 cm} |}
\hline
Image & DGS, non-redundant & DGS redundant & SL PSNR & FFST PSNR & CL PSNR & Wavelet PSNR\\ \hline
Straw & 16.1568 & \textbf{16.9645} & 15.9545 & 15.8711 & 15.5610 & 15.2016 \\ \hline
Rocks &\textbf{19.6596} & 19.6498 & 18.2772 & 15.9770 & 19.3774 & 19.0008 \\ \hline
Grass & 18.0498 &  \textbf{18.1756} & 17.9382 & 17.4305 & 17.9655 & 17.9847 \\ \hline
Cracked Mud & 15.5964 & \textbf{16.1645} & 15.2524 & 15.0387 & 15.1937 & 15.2714  \\ \hline
Bricks & \textbf{21.8677} & 20.9070 & 19.7579 & 19.5595 & 18.4367 & 18.9953 \\ \hline
Bars &  \textbf{21.6415} & 21.1399 & 19.0507 & 19.2957 & 17.6218 &18.2207  \\ \hline
Fabric & 16.2813 & \textbf{16.5387} & 15.7118 & 15.6807 & 15.2957 & 15.5106  \\ \hline
Grate & 17.7373 & \textbf{17.8509} & 15.5762 & 15.6462 & 14.7761 & 14.9119 \\ \hline
Honeycomb  & \textbf{19.9771} & 19.6670 & 14.6194 & 15.2500 & 12.8580 & 11.7106\\ \hline
Mandrill & 16.5017 & \textbf{16.5743} & 16.4203 & 15.1830 & 16.5653 & 16.2576 \\ \hline
Boat & 17.9514 & 18.0641 & 17.6947 & 16.4199 & \textbf{18.2324} & 18.0268 \\ \hline
Lena & 17.3064 & \textbf{17.6526} & 16.5380 & 15.3970 & 17.1781 & 16.7578 \\ \hline
\hline
\end{tabular} 
}
\end{center}
\caption{\label{tab:PSNR}Best-case thresholding PSNR for denoising experiments.}
\end{table}

\begin{figure}[H]
\centering
\includegraphics[width=.2\textwidth]{x1_1_03.pdf}
\includegraphics[width=.2\textwidth]{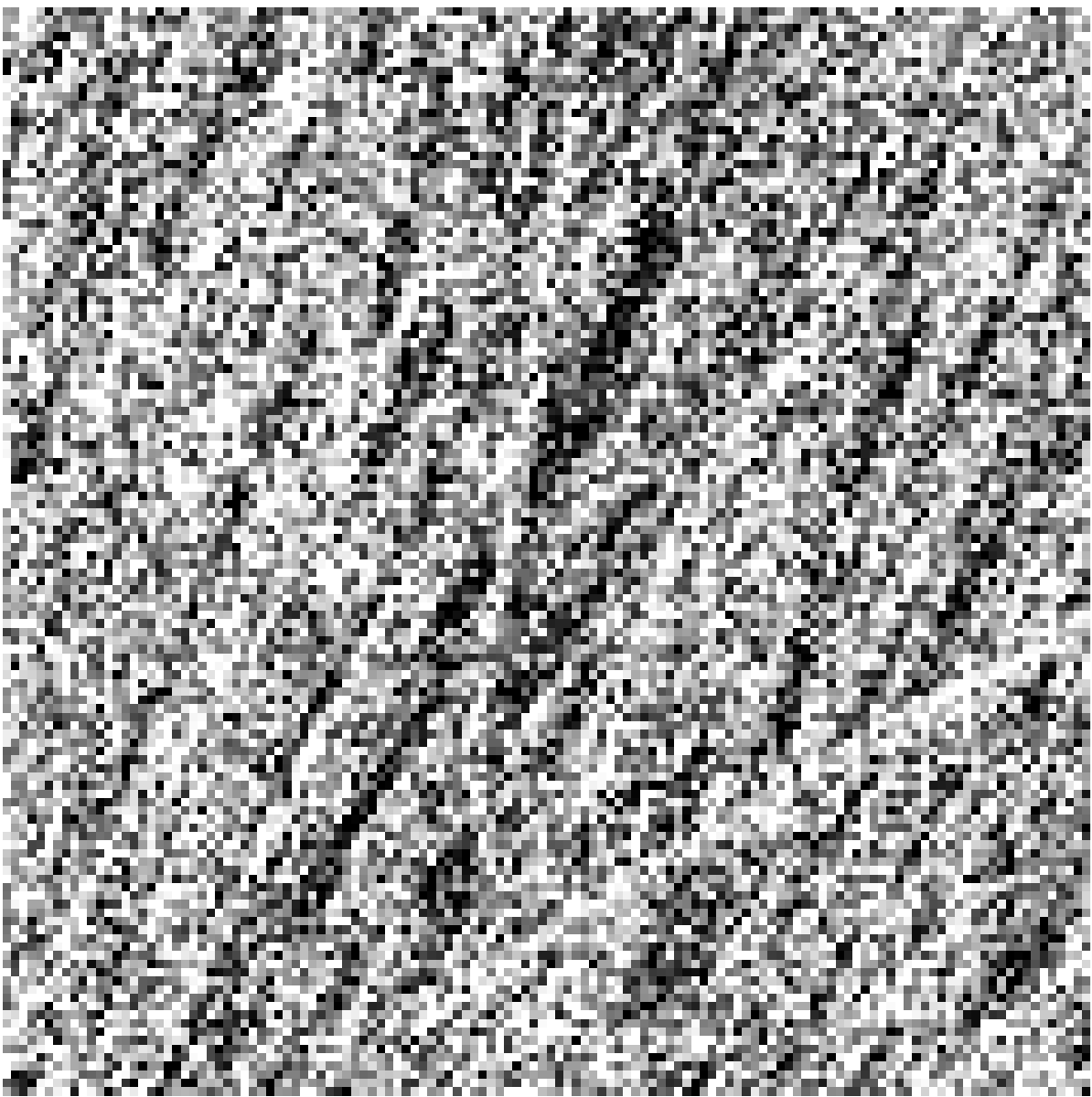}
\includegraphics[width=.2\textwidth]{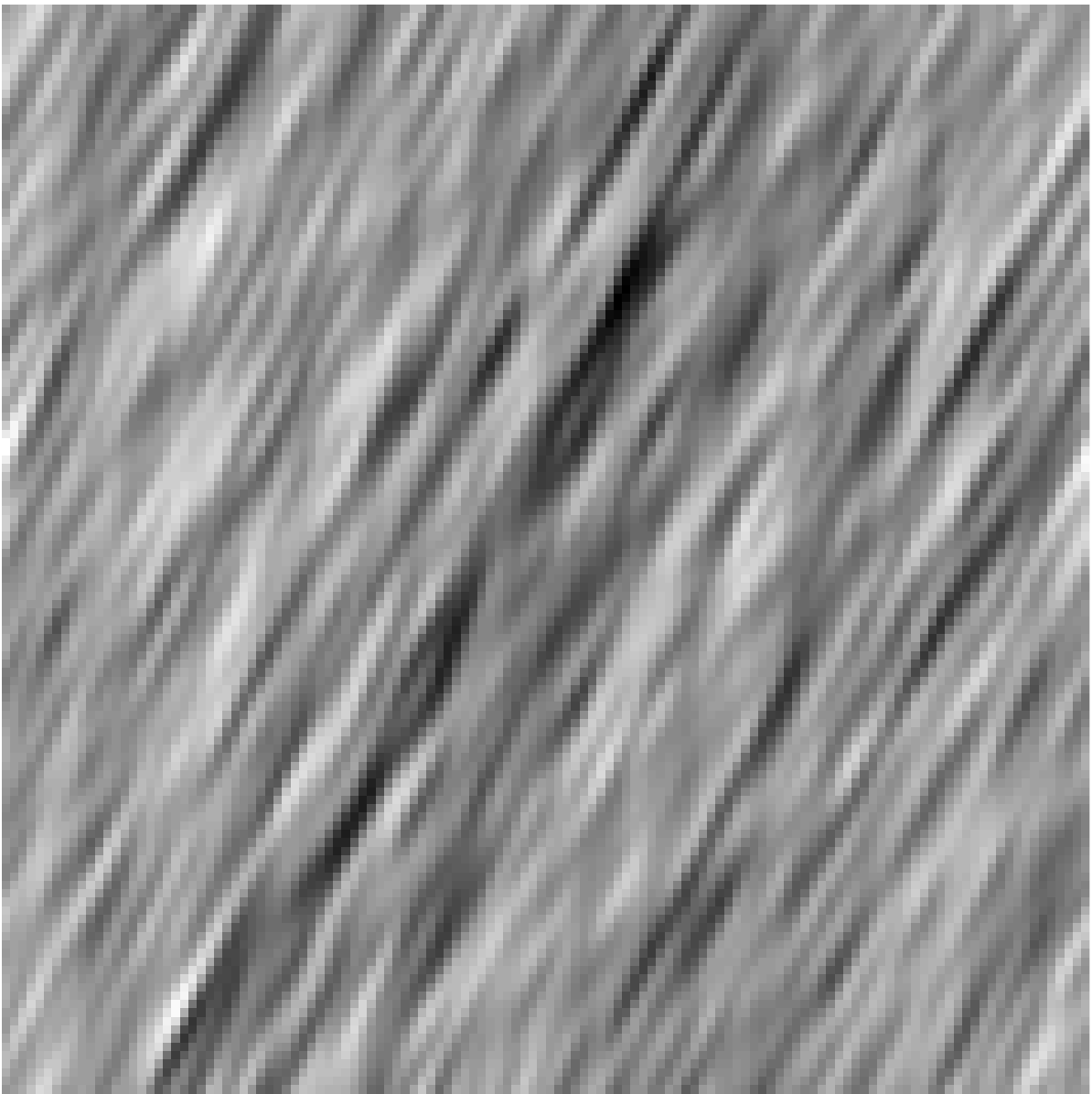}
\includegraphics[width=.2\textwidth]{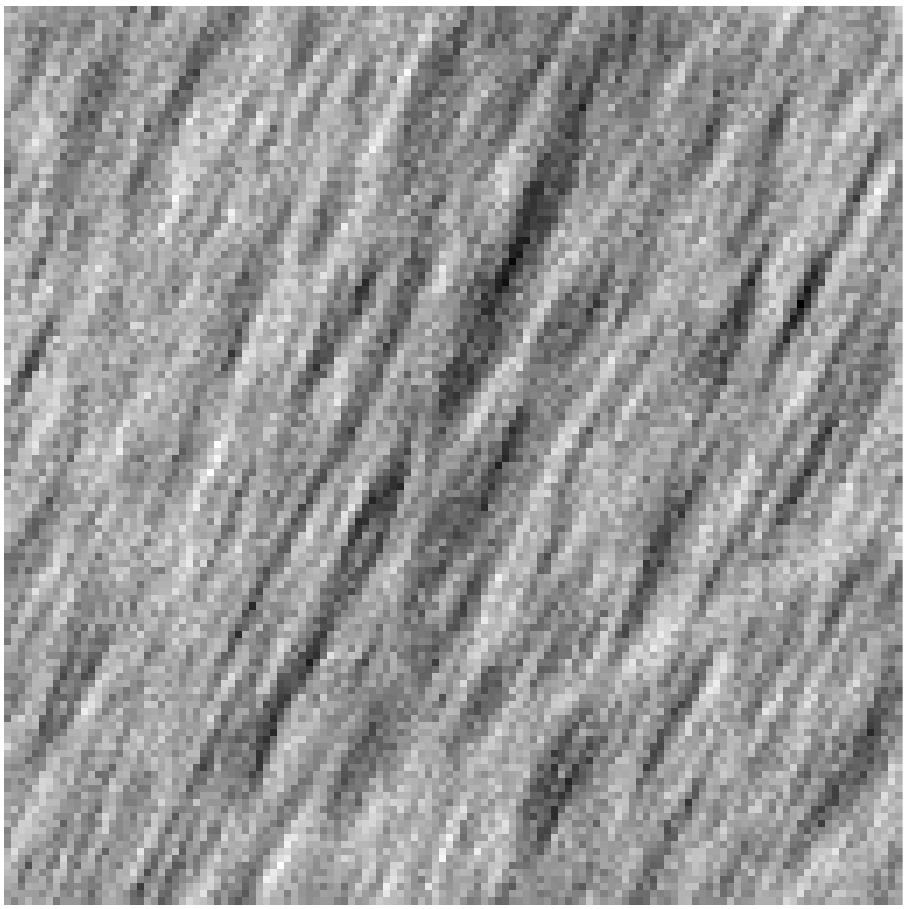}
\includegraphics[width=.2\textwidth]{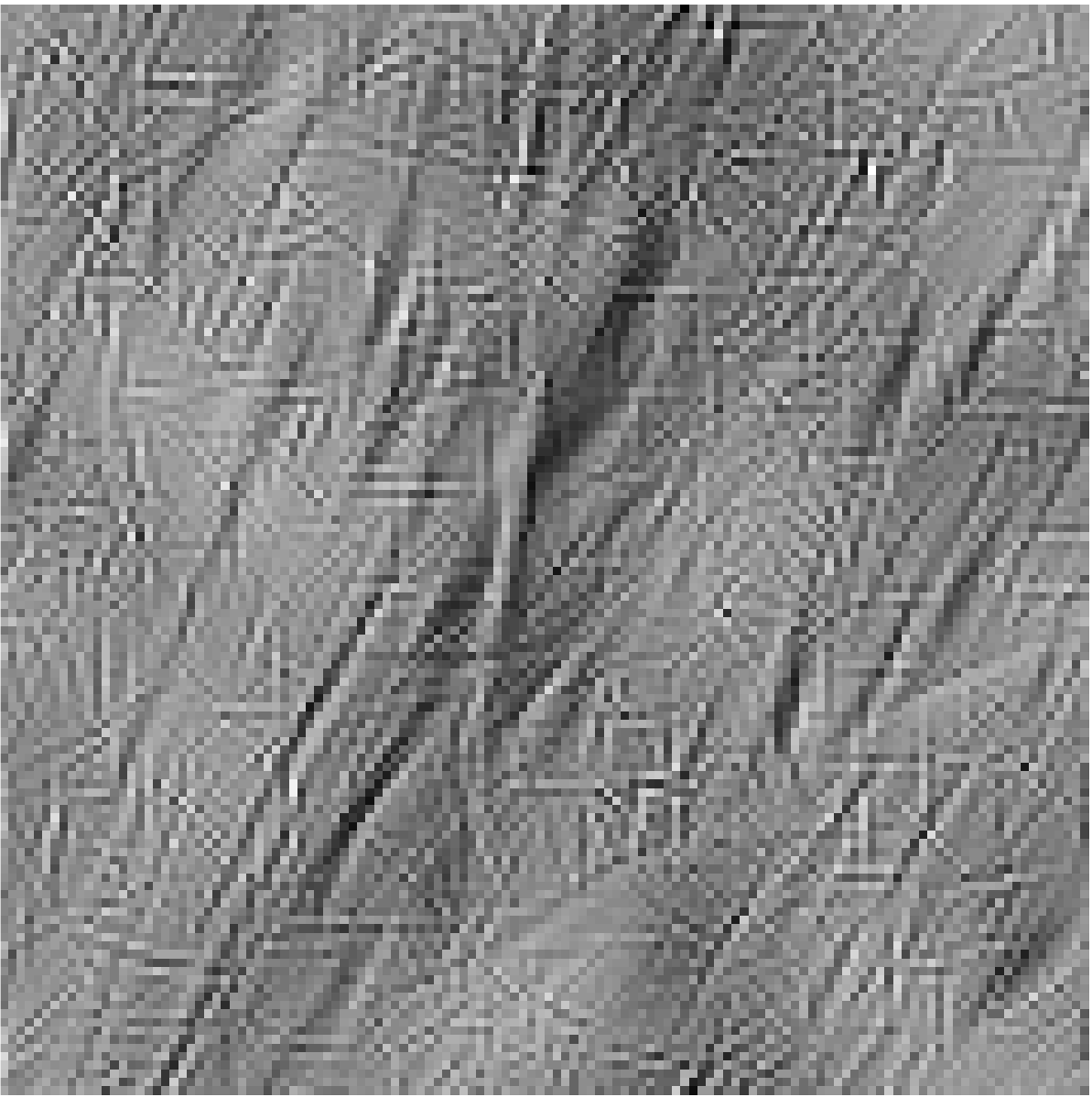}
\includegraphics[width=.2\textwidth]{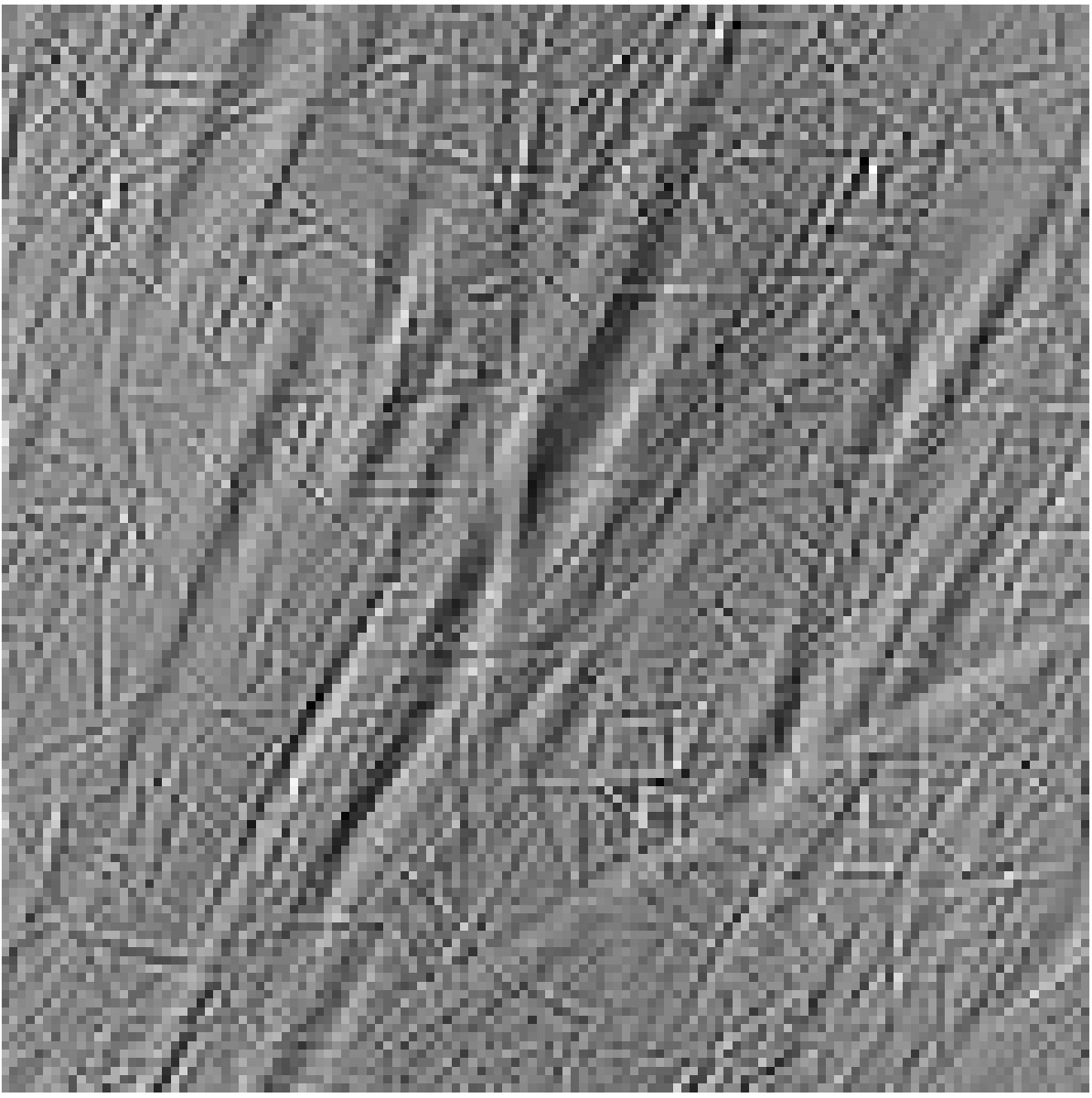}
\includegraphics[width=.2\textwidth]{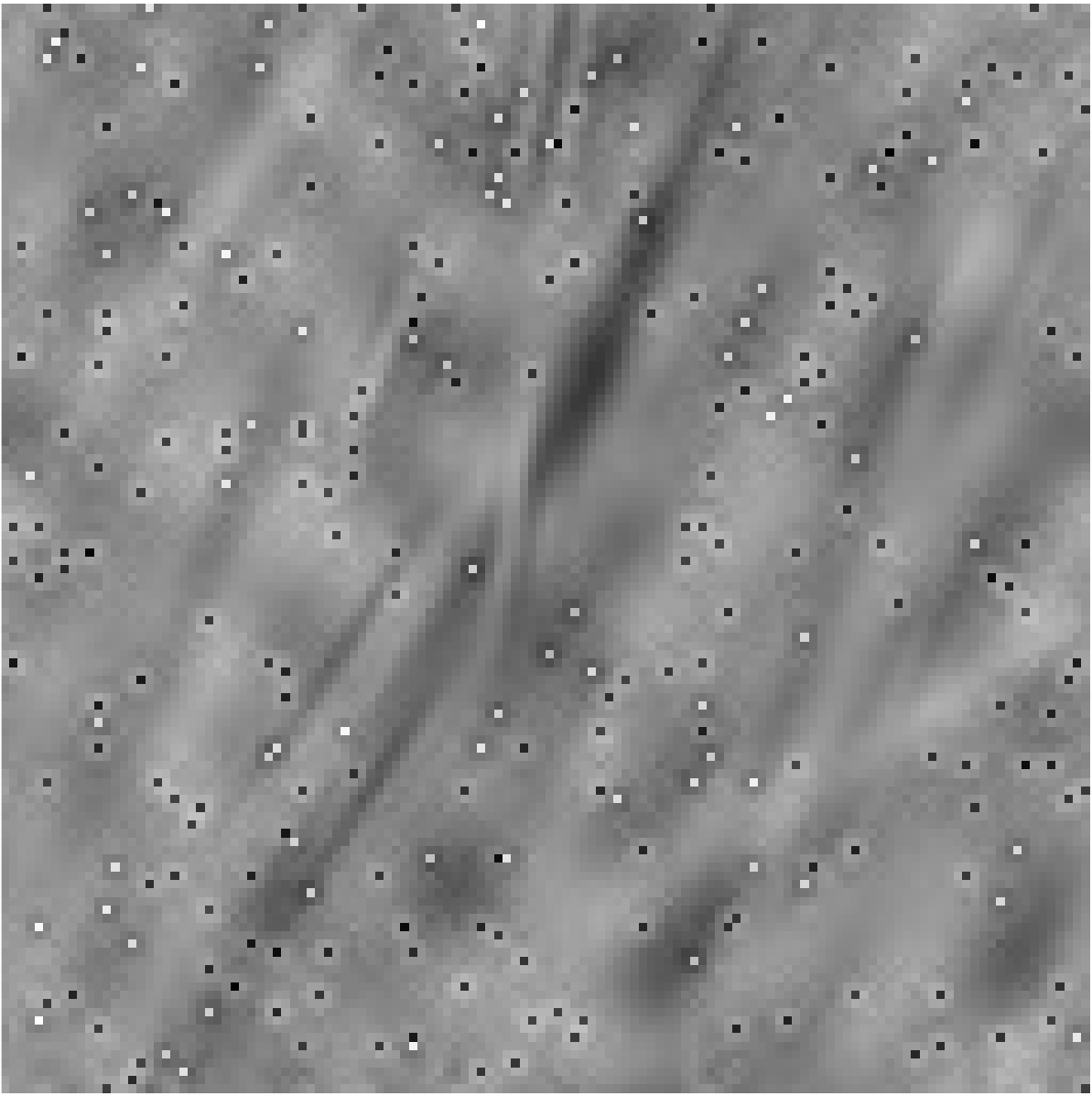}
\includegraphics[width=.2\textwidth]{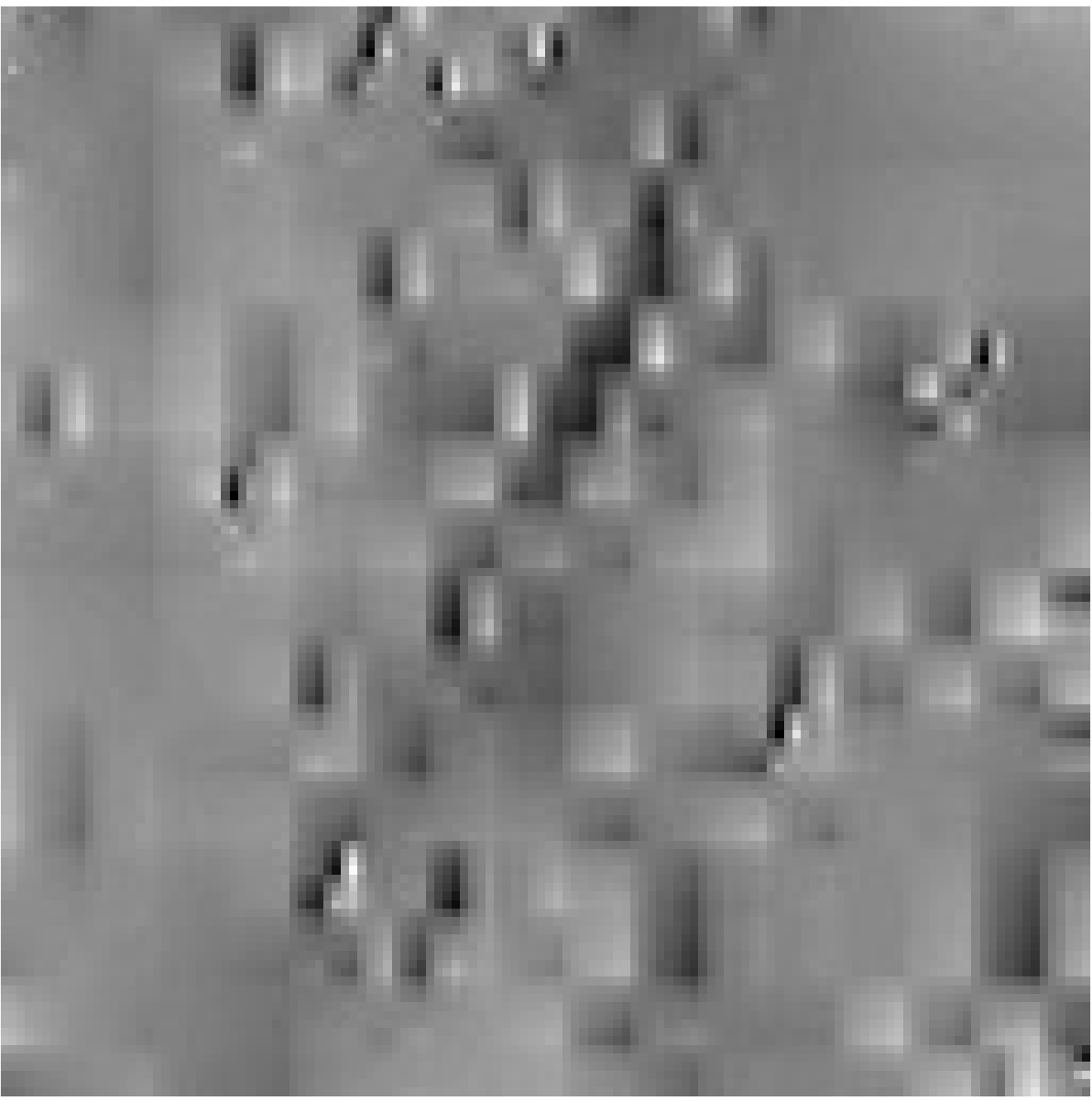}
\caption{Denoising results for straw.  First row, left to right: original, noisy, denoised with non-redundant directional Gabor system, denoised with redundant directional Gabor system.  Second row, left to right: denoised with shearlab shearlets, denoised with FFST shearlets, denoised with curvelab curvelets, denoised with wavelets.}
\label{denoised_straw}
\end{figure}

As in the compression experiments, directional Gabor systems perform competitively, particularly in the textural images.  For textures, directional Gabor systems are clearly optimal in all examples.  However, for the natural images, results are more ambiguous.   While directional Gabor systems are effective for Lena, the results are poorer for mandrill and boat.  This is reasonable, given that our method does not have the multiscale anisotropic properties of shearlets and curvelets that make them near-optimal for cartoon-like images.  Visual results for the straw texture appear in Figure \ref{denoised_straw}, which illustrates interesting properties of the directional Gabor frame.  Our method provides a relatively smooth denoised image, with few substantial artifacts.  The multiresolution methods all produce substantial denoising artifacts that affect the PSNR negatively, and also reduce visual quality.  However, our method does not resolve many of the edges as well as the anisotropic multiresolution methods.  For example, in Figure \ref{denoised_straw}, some of the most prominent edges are better resolved by shearlets and curvelets.  We note that there are other denoising regimes, beyond hard and soft universal thresholding, that can in some cases produce denoised images with better PSNR and visual quality.  

\subsection{Summary and Conclusions}

The experiments demonstrate competitive performance of directional Gabor systems against well-known redundant anisotropic frames, as well as non-redundant wavelets.  The performance of the directional Gabor system is particularly notable on textures, where it consistently outperforms other methods.   This is an important case in image processing, as many state-of-the-art frames, such as shearlets and curvelets, are not known to be near-optimal for textures.  We note that in its current prototype stage, our algorithm runs more slowly than FFST, ShearLab, or CurveLab.  A future research direction is to consider existing fast implementations of frames using ridge functions \cite{ridgelets_code}, in order to improve our prototype.

\subsection{Acknowledgments}

This work was supported in part by a research grant from the Defense
Threat Reduction Agency, HDTRA1-13-1-0015 ``Harmonic Analysis Methods for Autonomous Radiological Search: A Data Driven Approach", and by an Army Research Office grant W911NF1610008 ``Nonlinear and Probabilistic Analysis with Frames".  The authors also thank the USC-SIPI Image Database for collecting the images used for the numerical analysis.  We would also like to thank the reviewers for their helpful comments.

\clearpage
\bibliography{JMMcitations_WC5-17.bib}

\end{document}